\documentclass[12pt]{article}
\usepackage[margin=1in]{geometry}
\usepackage[utf8]{inputenc}
\usepackage{amssymb}
\usepackage{graphicx}
\usepackage{amsmath}
\usepackage{amsthm}
\usepackage[usenames,dvipsnames]{xcolor}
\definecolor{LightGray}{gray}{0.95}
\usepackage[colorinlistoftodos,color=LightGray]{todonotes}
\usepackage[numbers,square,comma,sort&compress]{natbib}
\usepackage{url}
\usepackage{hyperref}
\hypersetup{
  colorlinks   = true, 
  urlcolor     = Blue, 
  linkcolor    = Cerulean, 
  citecolor   = SeaGreen 
}

\usepackage[normalem]{ulem}
\usepackage{charter}
\usepackage{parskip}
\usepackage{subcaption}

\providecommand{\R}{\mathbb{R}}

\theoremstyle{plain}

\newtheorem{example}{Example}[section]

\theoremstyle{definition}
\newtheorem{definition}{Definition}
\theoremstyle{remark}

\newcommand{\captionfonts}{\small\it}
\makeatletter  % Allow the use of @ in command names
\long\def\@makecaption#1#2{%
  \vskip\abovecaptionskip
  \sbox\@tempboxa{{\captionfonts #1: #2}}%
  \ifdim \wd\@tempboxa >\hsize
    {\captionfonts #1: #2\par}
  \else
    \hbox to\hsize{\hfil\box\@tempboxa\hfil}%
  \fi
  \vskip\belowcaptionskip}
\makeatother   % Cancel the effect of \makeatletter

\newcommand{\beginsupplement}{%
        \setcounter{table}{0}
        \renewcommand{\thetable}{A\arabic{table}}%
        \setcounter{figure}{0}
        \renewcommand{\thefigure}{A\arabic{figure}}%
     }

\newcommand{\missingcite}[1]{\textcolor{ForestGreen}{$[\bf cite]$}}

\title{\Large The underlying connections between identifiability, active subspaces, and parameter space dimension reduction}

\author{Andrew F. Brouwer$^{*,\dagger}$ and Marisa C. Eisenberg$^{*,\dagger,\ddagger}$\\
\small $^*$ co-corresponding, brouweaf@umich.edu, marisae@umich.edu. \\
\small $^\dagger$ Department of Epidemiology, $^\ddagger$ Departments of Complex Systems and Mathematics, \\
\small University of Michigan, Ann Arbor}

\date{}

\begin{document}
\maketitle

\begin{abstract}
The interactions between parameters, model structure, and outputs can determine what inferences, predictions, and control strategies are possible for a given system. Parameter space reduction and parameter estimation---and, more generally, understanding the shape of the information contained in models with observational structure---are thus essential for many questions in mathematical modeling and uncertainty quantification. As such, different disciplines have developed methods in parallel for approaching the questions in their field. Many of these approaches, including identifiability, sloppiness, and active subspaces, use related ideas to address questions of parameter dimension reduction, parameter estimation, and robustness of inferences and quantities of interest. 

In this paper, we show that active subspace methods have intrinsic connections to methods from sensitivity analysis and identifiability, and indeed that it is possible to frame each approach in a unified framework. A particular form of the Fisher information matrix (FIM), which we denote the sensitivity FIM, is fundamental to all three approaches---active subspaces, identifiability, and sloppiness. Through a series of examples and case studies, we illustrate the properties of the sensitivity FIM in several contexts. These initial examples show that the interplay between local and global and linear and non-linear strongly impact the insights each approach can generate. These observations underline that one's approach to parameter dimension reduction should be driven by the scientific question and also open the door to using tools from the other approaches to generate useful insights.
\end{abstract}

\section*{Introduction}

Both parameter space dimension reduction and parameter identifiability are, fundamentally, a pursuit of the underlying structure of a map from an input space to an output space. This pursuit is an essential aspect of mathematical modeling, where a model is the map of interest. Indeed, different disciplines have developed methods in parallel for approaching dimension reduction, each with their own emphasis depending on the important questions of the parent field. Although much dimension reduction work has been done in the context of dynamical systems, many of the concepts and techniques apply to a wide range of models.

Identifiability, which emphasizes questions of parameter estimation (i.e., which parameters or parameter combinations can be uniquely determined from observed data),  primarily grew out of statistics and engineering, starting in the 1940s and '50s, with applications particularly focused in pharmacokinetics \cite{geary1941inherent,reiersol1950identifiability, koopmans1950identification, koopmans1949identification}. Since then, periods of renewed and more generalized interest have followed in the '70s and after 2000, particularly with advent of the differential algebra method for identifiability of dynamical systems~\cite{Bellman1970,Rothenberg1971,Reid1977,Pohjanpalo1978,Cobelli1980,Jacquez1985,Jacquez1990,Jacquez1991,Ljung1994,Audoly2001,Saccomani2001,Saccomani2003,Meshkat2009,Raue2009,Meshkat2011,Meshkat2012,Eisenberg2014}.) These methods often include parameter reduction approaches, often by either combining or fixing estimated parameters to ensure model identifiability.

More recently, in the early 2000s, the concept of model sloppiness was developed by researchers investigating dynamical systems in biology and physics, with the goal of developing reduced models with nearly the same dynamical properties by using differential geometry methods, like the manifold boundary approximation method (MBAM)~\cite{Brown2003,Gutenkunst2007,Transtrum2011,Transtrum2014,Tonsing2014,Transtrum2015,transtrum2016bridging}. Sloppiness concepts and techniques are closely related to the parameter reduction approaches seen in the identifiability literature and seek to identify sloppy (insensitive) and stiff (sensitive) directions in parameter space. Moreover, the reduced models identified by the methods should be able to generate the same overall input--output behavior as the original (i.e., given the same inputs, it generates the same outputs), thereby seeming to indicate that the original model was unidentifiable and, potentially, that the reduced model is identifiable (or closer to it). These connections have been explored in several papers, which have shown that, while connections to identifiability may not be one-to-one, sloppiness is a closely related concept~\cite{Chis2016, dufresne2016geometry}.

Active subspaces is a relatively new approach  that came out of uncertainty quantification and that seeks to reduce the number of parameters needed to approximate a quantity of interest, particularly in the case of large models~\cite{Smith2014,constantine2015active}. Larger models often exhibit challenges that may make more standard identifiability and sloppiness approaches difficult to implement, particularly when running the model is computationally intensive. Similar to both identifiability and sloppiness, active subspaces methods focus on understanding input--output relationships---how the model output changes as a function of its inputs or parameters---and on developing reduced models which generate a similar input--output structure as the original model. In the case of active subspaces, these ideas are framed around the idea of active (sensitive) directions in parameter space, versus inactive (insensitive) directions in parameter space. The inactive directions would seem to naturally correspond to compensation between parameters, potentially in an identifiable combination, while active directions might correspond to changing the value of the identifiable combination itself. In this paper we will examine this potential connection further.
 
Each of these methods has close connections to ideas of parameter sensitivities, which will form the basis of the links we will draw out in this paper. Indeed, some links along these lines have been drawn out in varying levels of detail for all three methods \cite{constantine2015active,Jacquez1985,Jacquez1990,Jacquez1991, Transtrum2015}. Unidentifiability, sloppiness, and inactive directions in parameter space are each manifestations of insensitivity of the model output to changes along some direction in parameter space. Sensitivity analysis methods often focus on determining which individual parameters are sensitive or insensitive, but the generalization of this idea to a multi-parameter case can lead one to any of the three approaches mentioned here.

For more comprehensive reviews of these approaches to dimension reduction and techniques popular in each field, we refer the reader to \cite{constantine2015active, Miao, chis2011structural,eisenberg2017confidence, Gutenkunst2007,Transtrum2015}. Here, we examine some of the underlying connections between identifiability and parameter reduction methods, in particular focusing on active subspaces and sloppiness. We show that some of the main objects and concepts in active subspaces and other parameter space reduction approaches can be framed in terms of a commonly used form of the Fisher Information Matrix which we will term the sensitivity Fisher Information Matrix (sFIM). While much of the material we present is a review or reframing of existing approaches, we hope that the translation dictionary developed here will facilitate cross-talk between these similar approaches. Finally, we illustrate how these concepts and approaches may interact in practice with a series of examples and case studies.

\section*{Framework and notation}
We begin by setting up the framework and notation we will use throughout. We will consider primarily either algebraic equations or ordinary differential equations (ODEs) of the form:
\begin{equation} \label{eq:model}
    \begin{aligned}
    \dot{x} &=w(x,t,\theta),\\
    y &= v(x,\theta),
    \end{aligned}
\end{equation}
where $t$ is time, $w$ and $v$ are functions, and $\theta$ represents the (vector of) parameters, which may in some cases include initial conditions, input variables, or other quantities affecting the model behavior. Here, $x$ is the (unobserved) state variable vector, and $y$ represents the measured  (observed) outputs or quantities of interest. 
In many cases, one might also have known inputs or forcing functions which drive the model and be included in $f$---these would typically be denoted $u$. For non-differential, algebraic models, we will also use the same notation of state variables ($x$), observed variables ($y$), and parameters ($\theta$). More generally, our notation will follow the following conventions:
\begin{description}
\item[$\theta$] -- Parameters, input variables, initial conditions or other varied quantities. These are often denoted $\boldsymbol{x}$ in the active subspace literature \cite{constantine2015active} and $p$ or $\theta$ in the identifiability literature. Here, $n$ is the length of the parameter vector, $\bf \theta = \{\theta_1, \dots, \theta_n\}$.
\item[$x$] -- Model state variables (unobserved), a vector.
\item[$q$, $y$] --  Model output (observed) or quantity or quantities of interest (QOI). A QOI can often be viewed as the output of the model, and is a function of $x$, $\theta$, and potentially time or other independent variables. In the active subspace literature, this is typically viewed as a scalar \cite{constantine2015active}, and is denoted $q$. In the identifiability and parameter estimation literature, the QOI might be a scalar in the form of the likelihood or sum of squares, or potentially a vector of model measurements or an observed trajectory of some function of the model variables (e.g., one of the variables scaled by a constant). If the QOI is a vector of measurements or an observed trajectory, it would more typically be denoted $y$ in the identifiability literature. We will thus use $q$ to represent a scalar QOI and $y$ to represent a vector or continuous QOI.
\item[$f$] -- Model map from the model parameters $\theta$ to the model output $q$ or $y$, $f:\Theta\subset\R^n \rightarrow \R^m$. For scalar $q$, $m=1$. In an identifiability context, $f$ is the model map used when evaluating injectivity, which the input--output equations represent implicitly (described further in the next section). 
\end{description}

\section*{Concepts in identifiability and parameter space reduction }

A model parameter is said to be \textit{identifiable} if it can be uniquely determined from the model output, and a model is identifiable if all of its parameters are identifiable. If a parameter is not identifiable, then the model output is either \textit{insensitive} to that parameter, or the parameter is part of an \textit{identifiable parameter combination}, meaning that, while the value of the parameter itself is not fixed by the model output, some function of it and other parameters is. For example, in the model $y=(m_1+m_2)x+b$ with $(x,y)$ pairs as the observed model output, parameter $b$ is identifiable and, while parameters $m_1$ and $m_2$ are individually unidentifiable, the sum $m_1+m_2$ is an identifiable parameter combination.

A common distinction in examining identifiability is between \textit{structural} versus \textit{practical} identifiability. Structural identifiability focuses purely on identifiability issues inherent to the model structure (such as in the linear example described above), while practical identifiability considers the estimation issues that come with real data (such as error, number or timing of samples taken, etc.). (In some cases, these two categories are denoted \textit{identifiability} and \textit{estimability}.) In the ODE case, structural identifiability is often framed  as a best-case scenario wherein the data are assumed to be known completely (i.e., smooth, noise-free, and continuously sampled), although one can also consider structural identifiability when particular measurement times are specified. Structural identifiability of a model is a necessary condition but not sufficient condition for parameter estimation with real-world data \cite{Cobelli1980}, since failure to recover parameters in the ideal case implies failure in the imperfect, i.e., real-world data case as well. More formally, we can define structural identifiability as follows:
\medskip
\begin{definition} An individual parameter $\theta_i$ in $\theta$ (Eq.~\eqref{eq:model}) is \emph{globally (also termed uniquely) structurally identifiable} if, for almost all values $\theta_i^*$ and initial conditions, the observation of an output ($y=y^*$) uniquely determines the parameter value $\theta_i$ ($\theta_i=\theta_i^*$), i.e., if only one value of $\theta_i$ could have resulted in the observed output. Similarly, a parameter $\theta_i$ is said to be \emph{locally (also termed non-uniquely) structurally identifiable} if there are a finite number of parameter values which can generate the observed output. 
\end{definition}

Similarly, a model is said to be globally (respectively locally) structurally identifiable if every parameter is globally (at least locally) structurally identifiable.  In subsequent sections, ``structurally identifiable'' will be understood to mean ``globally structurally identifiable'' unless otherwise indicated. If a model is not structurally identifiable, it is \textit{unidentifiable}, and there exists a set of identifiable combinations of parameters that represents the parametric information available in the data (except in degenerate cases where the model is reducible or has insensitive parameters)~\cite{Cobelli1980}. Such a set is not unique; any set of combinations that generates the same field is an equivalent set of identifiable combinations, e.g.,\ $\{\theta_1\theta_2, \theta_3/\theta_2\}$ and $\{\theta_1\theta_2,\theta_1\theta_3\}$ are equivalent sets of identifiable parameter combinations.

Many different analytical approaches to structural identifiability have been developed \cite{Chappell1998, Cobelli1980, Pohjanpalo1978, Vajda1989, Miao, chis2011structural}. However, analytical methods for identifiability can be  computationally intensive, making applications beyond relatively simple models challenging~\cite{Saccomani2003, Vajda1989, Pohjanpalo1978, chis2011structural}. 
Of these techniques, differential algebra has gained significant traction and has been the source of a range of recent advances the field of identifiability \cite{mahdi2014structural, meshkat2014identifiable, meshkat2015identifiability, brouwer2016systematic, haffke2015validation, harrington2016differential, saccomani2016structural, merkt2015higher, bearup2013input}. By contrast, while most numerical approaches to identifiability provide only local (rather than global) information about the parameters, they are often more computationally tractable~\cite{Hengl2007}. Many of these methods can be used to address both structural and practical identifiability, often by using simulated data (either without noise or with a range of different noise assumptions, depending on whether structural or practical identifiability is considered) \cite{Raue2009, Eisenberg2014}.

Both structural and practical identifiability  are often used in to inform parameter space reduction. One simplistic approach is to fix the values of individual parameters until the identifiable parameter combinations uniquely determine the remaining parameters, e.g., in the example with identifiable combinations $\theta_1\theta_2$ and $\theta_1\theta_3$, fixing $\theta_1=\theta_1^*$ would allow $\theta_2$ and $\theta_3$ to be uniquely determined from the identifiable parameters. 
Another, more elegant approach is to restructure and reparameterize the model in terms of its identifiable combinations, yielding an identifiable form of the model \cite{Chappell1998,Brouwer2017b}. In the structural identifiability case, this is often termed a \textit{identifiable reparameterization}. For example, a simple identifiable reparameterization would be to take the $y = (m_1 + m_2)x + b$ example above and define $m = m_1 + m_2$ so that our model is now the identifiable equation $y = mx + b$, with two parameters identifiable from $(x,y)$ data, $m$ and $b$. An example of a practical parameter reduction is the linear approximation of a Hill function~\cite{Holmberg1982}, where $y = \frac{V x}{x + K}$. If all $(x,y)$ data is in the region where $x\ll K$, then this model becomes practically unidentifiable and can be approximated by the identifiable linear model $y = m x$, where $m = V/K$. These examples are extremely simple compared to most models used in practice, but they illustrate the concepts. 

Parameter space dimension reduction based on structural identifiability takes advantage of the inherent structure of the model by reparameterizing only in terms of the structurally identifiable parameter combinations; the dimension of the parameter space is reduced but no information is lost. 
Viewing the model as a map from parameter space to output space, we can consider the fibers of this map (i.e., the sets of all input values that correspond to each output value or trajectory). If the fibers contain only one or finitely many elements (i.e., only one or finitely many parameter values can generate a given output), then the model is identifiable, or equivalently, cannot be structurally reduced in dimension. However, if the fibers contain infinitely many elements, the model is unidentifiable---or put in a dimension reduction framework, these fibers represent opportunities for dimension reduction, as they all yield the same output. The dimension reduction process can then be viewed as collapsing/taking a set of representatives of the equivalence classes generated by the model map from parameters to output. These reduced dimension models would be termed identifiable reparameterizations of the model in an identifiability context.

Practical dimension reduction, whether based on practical identifiability approaches or more general numerical approaches to the problem, identifies lower-dimensional models with output that is nearly indistinguishable, but not necessarily equivalent, to the output of the original, higher-dimensional model. In the parameter estimation context, this might occur because the collected data cannot distinguish (for some level of significance) between similar output trajectories associated with different parts of parameter space. Such problems can arise, for example, from excessive noise or measurement error or because of insufficiently generic measurement times---for instance, if one measures the value of a periodic function only once per period, the amount of information in the data for a periodic model is limited. Of course, this may be framed as a problem or an opportunity depending on the question at hand. As for structural identifiability, practical identifiability approaches to dimension reduction can find practically identifiable parameter combinations, although there may not exist reparameterizations of the original model in terms of these practical parameter combinations. The scientific (and philosophical) implications of practical dimension reduction, particularly in biological models, have been a focus the sloppiness literature among others~\cite{transtrum2016bridging, villaverde2017dynamical}. Practical dimension reduction is also used, as in the active subspace literature, to find computationally tractable approximations to to computationally intensive models. The focus in this context is usually not on the scientific implications but rather on the effective outcomes.

Dimension reduction depends not only the model but also on the output considered---parameters that are identifiable for one kind of observed output may not be for another. This observed output might be a single quantity of interest (QOI), as is typical in the active subspaces literature; among the mathematical biology/dynamical systems literature, on the other hand, the output is typically a trajectory measured over time. If one observes a trajectory over time, then each data point could be considered its own QOI. Alternatively, one can aggregate the fit of the model output to all data points simultaneously in one cost function, as in the case of maximum likelihood estimation. The dimension reduction techniques will only identify parameter combinations that are common to all QOIs considered. Hence, which QOIs to consider in one's analysis should be driven by one's question. In a periodic model, for example, do you want to find the parameter space that closely matches the observed output, or do you simply want to match the period and amplitude?

\section*{The sensitivity Fisher Information Matrix}\label{sec:sFIM}

Next, we introduce the sensitivity-matrix formulation of the Fisher Information Matrix (FIM), which will underpin our development of the connections between parameter space reduction methods.  This formulation of the FIM has useful identifiability properties \cite{Rothenberg1971} and has---like sensitivity analysis more generally---a long history of use in identifiability and parameter space reduction \cite{Rothenberg1971,Jacquez1985,Jacquez1990,Jacquez1991,Cobelli1980,Cintron-Arias2009,capaldi2012parameter,banks2007sensitivity,yue2006insights,rodriguez2006hybrid,balsa2010iterative}. 

\medskip 
\begin{definition}
For a vector-valued function $f(\theta): \Theta\subset\mathbb{R}^n\to\mathbb{R}^m$, the \textit{sensitivity Fisher Information Matrix (sFIM)}, denoted $F(f;\theta)$, is the symmetric, $n\times n$ matrix whose $(i,j)$\textsuperscript{th} element is given by 
\begin{equation}
F_{ij}(f;\theta) = \sum_{k=1}^m \left(\frac{\partial f_k}{\partial\theta_i}\right)\left(\frac{\partial f_k}{\partial\theta_j}\right).
\end{equation}
\end{definition}

\textbf{Remark.} Although $F$ is often simply called the Fisher Information Matrix in the literature, we want to distinguish between this object and the expected Fisher Information Matrix; we compare and contrast these objects in a later section.

The entries of $F$ are the sensitivity coefficients \cite{Smith2014, distefano2015dynamic} of our quantity of interest $f(\theta)$, i.e., the partial derivatives of $f$ with respect to each parameter. Sensitivity coefficients are the core objects of local sensitivity analysis, and they link naturally to questions of parameter space reduction and identifiability---if a parameter is insensitive, it is (practically or structurally) unidentifiable, and the model can be reduced. (Note that the reverse is not true: individual parameters may be highly sensitive but also unidentifiable).

For a univariate (i.e., $m=1$) QOI, $q = f(\theta)$ , as in the active subspaces literature, $F(f)$ can be conveniently written as 
\begin{equation}
F(f,\theta) = \left(\nabla f\right)\left(\nabla f\right)^T
\end{equation} 
where $\nabla f$ is the column gradient vector of sensitivities
\begin{equation}
\nabla f= \begin{bmatrix}
\frac{\partial f}{\partial \theta_1}\\ \vdots \\ \frac{\partial f}{\partial \theta_n}
\end{bmatrix}.
\end{equation}

For a multivariate QOI,  $y=f(\theta)$, it is convenient to use the Jacobian 
\begin{equation}
\chi = J(f) = \begin{bmatrix} \frac{\partial f_1}{\partial\theta_1} & \cdots & \frac{\partial f_1}{\partial\theta_n}\\ \vdots & \dots & \vdots \\ \frac{\partial f_m}{\partial\theta_1} & \cdots & \frac{\partial f_m}{\partial\theta_n},\\
\end{bmatrix}
\end{equation}
which is often called the \textit{sensitivity matrix} in this context~\cite{Cobelli1980}, and write
\begin{equation}
F(f)=\chi^T\chi.
\end{equation}

\textbf{Remark.} Depending on the context, an analytic formula for $f$ may or may not be available, and so the derivatives are often calculated numerically. Constantine~\cite{constantine2015active} discusses some practical considerations in gradient calculation.

\textbf{Remark.} We note that the forms given for $\chi$ and $F$ assume a vector of discrete QOIs forming $y$, while many definitions of structural identifiability use the full trajectory of the model as the output or QOI (i.e., with complete, continuous temporal and/or spatial information for the model). In such cases, we can typically approximate the full trajectory by taking very frequent samples (as would be generated in most numerical solvers), allowing us to numerically evaluate local structural identifiability. Alternatively, we may also define the output as being measured at specific times.

To understand why $F$ is an important object, we first consider $\chi$. The linear approximation of the change in the output $f$ as a function of the change in the parameters $\theta$ can be written as
\begin{equation}
\Delta f \approx \chi \Delta\theta.
\end{equation}
In the early identifiability literature, $f$ was said to be ``sensitivity identifiable'' if $\Delta\theta$ was (locally) recoverable from $\Delta f$~\cite{Reid1977,Cobelli1980}. More generally, it has been shown that the rank of $\chi$---or, equivalently, $F$---is the number of locally identifiable combinations, so that the model is locally identifiable when $\chi$ or $F$ has full rank \cite{Rothenberg1971, Jacquez1985, Jacquez1990, distefano2015dynamic}. Here, it is useful to think of $\chi$ as a map from $\mathbb{R}^n$ to $\mathbb{R}^m$. The nullspace of $\chi$ at any point $\theta$ in parameter space gives the linearization of the structures along which parameters can move without changing $f$ (indicating the identifiable combinations), so characterizing the nullspace of $\chi$ is one approach to the goal of dimension reduction. Alternatively, we could characterize the coimage of $\chi$, that is, the quotient space $\mathbb{R}^n/\ker{\chi}$ (also called the orthogonal complement of the kernel or the row space in the language of linear algebra). The rank of this space---which is at most $m$---is the number of identifiable parameter combinations, and finding a basis for this space identifies the linearizations of the parameter combinations. Hence, from the perspective of parameter reduction, we hope that $\chi$ does not have full rank. However, from the perspective of parameter estimation, we hope that it does. 

In practice, it is easier  to work with the map $F=\chi^T\chi$, which has the same nullspace and coimage as $\chi$ but also is square, symmetric, positive semi-definite, and, if full rank, invertible and positive definite. As noted above, if $F$ has full rank at $\theta$, then we say that $f$ is \textit{(locally) structurally identifiable} at $\theta$. In many situations, $F$ is nearly rank deficient, with one or more eigenvalues $\lambda$ close to zero. Near rank deficiency is an indication of practical unidentifiability and an opportunity for practical dimension reduction. The sloppiness literature has called models sloppy when $\lambda_{\max}/\lambda_{\min}$ is large. In this situation, certain directions in parameter space---corresponding the eigenvectors of the small eigenvalues---do not greatly affect the value of $f$, at least relative to the change along the directions of the eigenvectors of the large eigenvalues. The notion of large vs small eigenvalues is ill-defined, which is a common complaint about the notion of sloppiness (e.g.~\cite{Chis2016}). 

The sensitivity FIM $F(f;\theta)$, is a local object, as it depends on $\theta$. It can be full rank in some regions of parameter space but nearly rank deficient in others; that is, practical identifiable parameter combinations in one region might be resolved into individually identifiable parameters elsewhere. For example, when measuring the sum of two periodic functions, the amplitudes of two functions may be difficult to distinguish if their periods are similar but easy to distinguish if their periods are very different. 

\subsection*{The expected and sensitivity FIMs}
Although we developed the sFIM above in the parameter sensitivity context, it is closely related to (and can be viewed as a special case of) a similar object encountered in the parameter estimation context---namely, the expected Fisher Information Matrix.
In parameter estimation, the fit of a model to the available data $z$ is usually measured by a cost function, often a statistical likelihood $\mathcal{L}(z,\theta)$. Least-squares fitting falls into this category since, it is equivalent to maximum likelihood assuming Gaussian measurement error. In this context, we do not consider the fit to each point individually but rather to all points as a whole. Here, we are interested in shape of information near some parameter vector, typically the maximum likelihood estimate,
\begin{equation}
\hat\theta=\arg\max_{\theta} \mathcal{L}(z;\theta) = \arg\min_{\theta} (- \log\mathcal{L}(z;\theta))
\end{equation}

The expected Fisher Information matrix is an important information-theoretic object associated with a likelihood function $\mathcal{L}(z,\theta)$ and is given by
\begin{equation}
\begin{aligned}
I(\theta)&=\int\left[\left(\nabla \log \mathcal{L}(z;\theta)\right)  \left(\nabla \log \mathcal{L}(z;\theta)\right)^T\right] \mathcal{L}(z;\theta) \ dz,\\
&= E_z[\left(\nabla \log \mathcal{L}(z;\theta)\right)  \left(\nabla \log \mathcal{L}(z;\theta)\right)^T].
\end{aligned}
\end{equation}

Because $I(\theta)$ is integrated over the data $z$, it is dependent only on the parameters $\theta$. The sensitivity FIM $F(\log \mathcal{L},\theta)$  is a special case of expected FIM $I(\theta)$ when $z$ has Gaussian error with mean zero and variance one~\cite{Cobelli1980}.

Given sufficient regularity, the entries of $I(\theta)$ can be written
\begin{equation}
I_{ij}(\theta)%&=-E\left[ \left(\frac{\partial}{\partial \theta_i} \log\mathcal{L}(z;\theta)\right)\left(\frac{\partial}{\partial \theta_i} \log\mathcal{L}(z;\theta)\right)  \right]\\
= - E_z\left[\frac{\partial^2}{\partial \theta_i\partial\theta_j} \log\mathcal{L}(z;\theta)\right].
\end{equation}

Thus, the negative Hessian matrix of $\log\mathcal{L}$ evaluated at the maximum likelihood estimate $\hat\theta$ (also called the observed information matrix), is sometimes used to assess parameter sensitivity and identifiability~\cite{Little2010}.

\section*{The underlying connection between active subspaces and the sensitivity FIM}
With the sFIM introduced, we now present the main objects used in active subspaces, framed around the sFIM. Unlike the local sFIM, the active subspace approach is more interested in the behavior of a function over large regions of parameter space. The active subspace approach thus integrates $F$ over the parameter space with respect to some density $\rho$ on $\theta$ (often uniform or Gaussian), to create a global object, denoted $C$ in \cite{constantine2015active}. In our, more general notation, we can write $C$ as
\begin{equation}
C=\int F(f;\theta) \rho(\theta)\ d\theta.
\end{equation}

This matrix $C$, then, is \textit{average sensitivity FIM} over parameter space and is a global object. Unlike $F$, which has rank at most $m$, $C$ is not similarly constrained and can have rank up to $n$. Active subspace analysis using $C$ will capture unidentifiable parameters that are linear (or nearly so) over the whole space but will not find more non-linear parameter combinations, as we will see in the examples.

\subsection*{Eigendecomposition of local and average sensitivity FIMs}

Because local and average sensitivity FIMs $F$ and $C$ are symmetric, real matrices, they have orthogonal eigendecompositions 
\begin{equation}
Q\Lambda Q^T.
\end{equation}

The active subspace and sloppiness literatures have used the eigendecomposition of what we call the average and local sensitivity FIMs, respectively, to designate the directions spanned by eigenvectors corresponding to large eigenvalues as \textit{active}/\textit{stiff} and those corresponding to zero or small eigenvalues as \textit{inactive}/\textit{sloppy}. One could also call these directions \textit{identifiable} and \textit{unidentifiable} (whether practically or structurally so), although we will use the active subspace notation here. Again, the designation of eigenvalues as large or small has largely been relatively ad-hoc and driven by the questions of the investigators.

Denote $Q_a$ as the matrix whose columns are the active eigenvectors of $C$ and $Q_i$ as the matrix whose columns are the inactive eiegenvectors. Then, $Q_a\theta$, which is called the active variables in the active subspace literature, are the linearized identifiable parameter combinations. Taking advantage of the decomposition~\cite{constantine2015active}
\begin{equation}
\theta= I\theta = QQ^T\theta = Q_aQ_a^T\theta+Q_iQ_i^T\theta,
\end{equation}

one can approximate $f$ as 
\begin{equation}
g(\theta):= f(Q_aQ_a^T\theta) \approx f(\theta),
\end{equation}
though, in practice, creating a response surface approximation for $f$ through  regression on the active directions may be a more robust and less computationally intensive choice~\cite{constantine2015active}.

If all identifiable combinations are linear, then $F$ and $C$ will have the same eigenvectors. If there are non-linear identifiable combinations, then the two matrices will have different eigenspaces and, accordingly, different approximations of $f$. The approach one uses should be directed by the scientific question.

\subsection*{The cost function as a link between scalar and vector QOIs}
Both of the sensitivity FIM objects $F$ and $C$ introduced thus far can be considered in either the scalar or vector case, 
although researchers using active subspaces have reported difficulty or the need for further examination of using active subspace methods for vector-valued QOIs~\cite{constantinewebsite, constantine2015active}.

However, the framing of $C$ in terms of $F$ suggests one possible option for a scalar summary QOI to be used when we are considering a relatively local analysis---the cost function (e.g. least squares, or other norms and likelihood functions). In the case of the local sFIM $F$, there is a natural link between the vector QOI where $y$ is a vector of measurements and a cost function as the analogous scalar QOI, as follows: let $F_y(\theta)$ be the sFIM for a vector of measurements $y$, and suppose we evaluate $F_y$ at a point $\hat{\theta}$. Then, let us consider $q(\theta) = \sum_1^m(y(\theta) - y(\hat{\theta}))^2$ to be the least squares cost function using $y(\hat{\theta})$ as the `data' (although we note a range of other cost functions would also work). The QOI $q$ is zero at $\theta = \hat{\theta}$, and will remain (roughly) zero if we perturb $\theta$ in an insensitive/inactive/unidentifiable direction for $y$. Conversely, it will increase as we move $\theta$ in a sensitive/active/identifiable direction for $y$. Thus, when evaluated at $\hat{\theta}$, $F_q(\hat{\theta})$, should broadly be sensitive/insensitive in the same/similar directions as $F_y(\hat{\theta})$.

Then, if we are considering examining $C$ using a vector QOI in a region around some nominal set of parameters, the cost function or likelihood may provide a useful way to generate a scalar, summary QOI, although it does then tie one's QOI to a specific local area of parameter space. A similar approach could also be applied if one is working in a parameter estimation context, where the $y(\hat{\theta})$ might be replaced with the data set used for estimation.

\section*{Visual tools for parameter space reduction}
There are several visual tools  used to assess the parameter identifiability or opportunities for parameter space reduction. \textit{Parameter profile plots}, developed to assess identifiability, is the gold standard for determining practical identifiability, and is not based on the sFIM. \textit{Sufficient summary plots}, on the other hand, are used to visualize results of active subspace analysis and are thus based on the average sFIM.

\subsection*{Parameter profiles}
The profile likelihood is a  visual tool in assessing parameter identifiability~\cite{Raue2009,murphy2000profile,venzon1988method} that we will generalize here as the parameter profile.
Conceptually, one tries to identify ridges and other structures in the response surface. This approach `profiles' a single parameter $\theta_i$ by fixing the value of $\theta_i$ across a range of values and fitting all remaining parameters for each fixed value of $\theta_i$ either to data or to a model trajectory. The optimal cost function value at each value of $\theta_i$ constitutes the likelihood profile for the fixed parameter.

The first step in profiling a parameter is to select a point in parameter space $\hat\theta$.  In some contexts, this point will be the maximum likelihood estimate, though this may be generalized for broader parameter reduction questions to simply be a set of nominal parameters around which we plan to profile. Then, to profile a parameter $\theta_i$, we fix $\theta_i$ at a value $\theta_i^*$. Define $\theta_{j\not=i}=\{\theta_j|j\not=i\}$ and 
\begin{equation}
c(\theta_{j\not=i}):=(f(\theta_i^*,\theta_{j\not=i})-f(\hat\theta))^2.
\end{equation}

\textbf{Remark.} This definition of a cost function is useful for profiling parameters in the dimension reduction context. In the parameter estimation context, on the other hand, the cost function will be an assessment of fit to the data, typically a likelihood, and this profile is called a profile likelihood.

The cost function $c$ is a map $\mathbb{R}^{n-1}\to\mathbb{R}^+\cup\{0\}$. Although we've defined $c$ here to correspond to minmization of the $\mathcal{L}^2$ norm, other metrics might be appropriate depending on the context.

Let 
\begin{equation}
\theta_{j\not=i}^*:=\text{argmin } c(\theta_{j\not=i}).
\end{equation}

\begin{figure}[t]
\centering
\includegraphics[width=\textwidth]{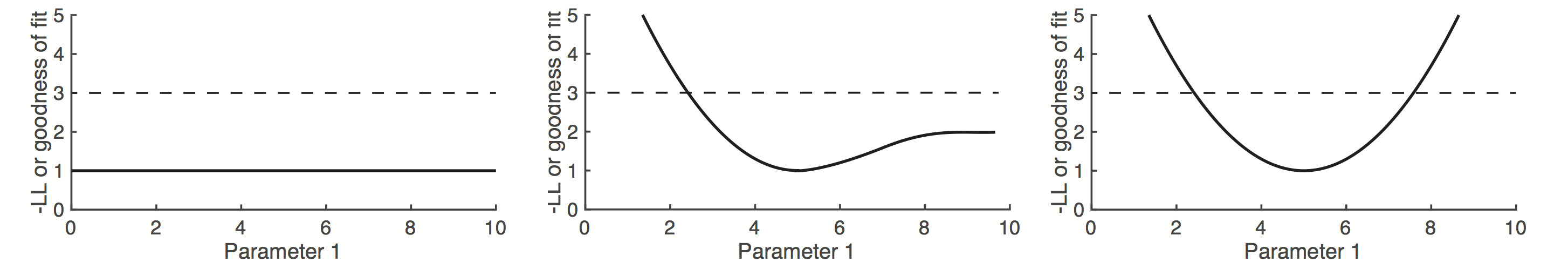}\\
\includegraphics[width=\textwidth]{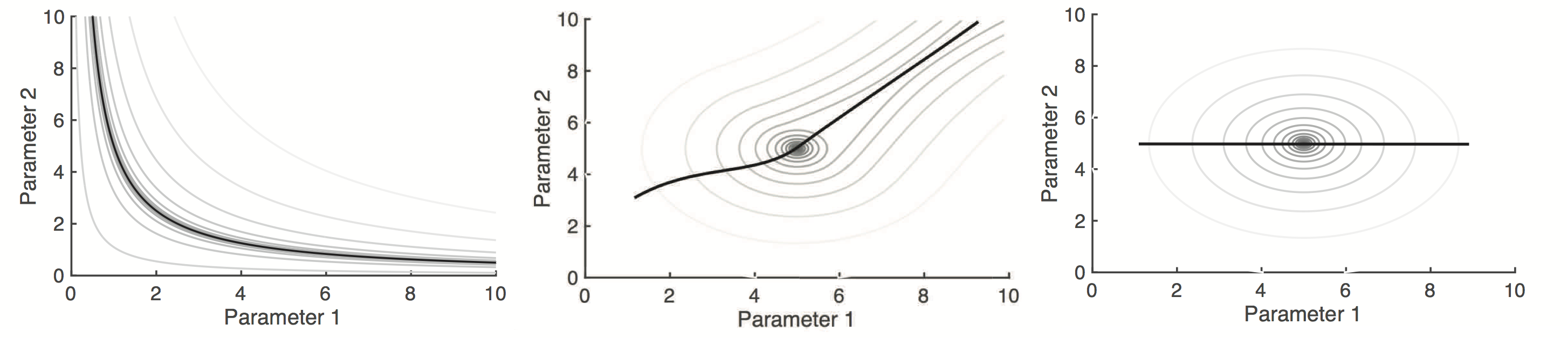}
\caption{Example profile likelihoods and parameter relationship contour plots for structural unidentifiabilty, practical unidentifiability, and identifiability. Top row: profile likelihoods (solid line) showing negative log likelihood (-LL) or goodness of fit values. Dashed line shows a confidence interval threshold (e.g. for 95\% confidence bounds), which is infinite, finite on one side, or bounded, respectively (left to right). Bottom row: corresponding parameter relationship plots to the top row. Contours indicate goodness of fit, with the best fit value of $p_2$ (Parameter 2) for each fixed value of $p_1$ (Parameter 1) shown as a black line. The structurally unidentifiable case illustrates a combination of the form $p_1 p_2$.}
\label{fig:proflike}
\end{figure}

That is, in the dimension reduction context, we are looking for the values of all parameters---excluding $\theta_i$, which is fixed at $\theta_i^*$---that make the function as close as possible to $f(\hat\theta)$. The plot of $c(\theta_{j\not=i}^*)$ vs $\theta_i^*$ is called a parameter profile. The shape of this profile is informative. As shown in Figure~\ref{fig:proflike}, if the shape is concave (i.e., trough-shaped), the parameter is practically identifiable. If the profile is flat or flat on one side, the parameter is at least practically (and potentially structurally) unidentifiable. In the practical identifiability literature, one defines a threshold value $\Delta$ and a confidence interval for $\theta_j$,
\begin{equation}
\{\theta_i^* |~c(\theta_{j\not=i}^*)<\Delta\}.
\end{equation}

The parameter is said to be practically identifiable if the confidence interval is finite~\cite{Raue2009,Brouwer2017b}; whether or not the parameter is practically identifiable can depend on the threshold chosen (much like the designation of active or inactive subspaces depends on the eigenvalue cut-off). When the cost function is a relative negative log-likelihood, $2\Delta$ is given by the chi-squared distribution $\chi^2(1-\alpha,n)$ where $\alpha$ is the level of significance ($\alpha$=0.05 for 95\% confidence intervals) and $n$ is the number of parameters~\cite{Raue2009}, but the choice of $\Delta$ in other contexts is more heuristic.

Perhaps more useful in considering parameter space reduction, often one uses plots of $\theta_{j\not=i}^*$ vs $\theta_i^*$ to identify the relationship between parameters in identifiable parameter combinations~\cite{Eisenberg2014} (bottom row of Figure~\ref{fig:proflike}). By examining how the estimates of the remaining parameters change as we profile a particular parameter $\theta_j$, we can trace out the form of the identifiable combinations (as developed in \cite{Raue2009}). For example, if we have a combination $\theta_1 + \theta_2$, then when profiling $\theta_1$, we would expect $\theta_2$ to change in a compensatory way which preserves the sum $\theta_1 + \theta_2$, as this will preserve the value of identifiable combination at $\hat\theta$ (or fit to the data).  However, as noted in \cite{Raue2009, Eisenberg2014}, this approach is ill-conditioned when there are multiple parameters in a combination or multiple combinations. Any extra degree of freedom in a combination allows the fitted parameters to compensate for one another and avoid tracing out the form of the identifiable combination with the profiled parameter. For example,  when profiling $\theta_1$, if our combination is $\theta_1 + \theta_2 + \theta_3$ then there are infinitely many ways that $\theta_2$ and $\theta_3$ can compensate to maintain the sum $\theta_1+\theta_2+\theta_3$ (and thus maintain the same fit to the data). The resulting profiled parameter relationships are often noisy or arbitrary, as there is a range of values for the unidentifiable parameters which will yield the same output in the profile. These issues can be addressed by restricting the set of parameters used for profiling to maintain appropriate degrees of freedom, such as using an sFIM-based approach \cite{Eisenberg2014}, as well as other methods \cite{Hengl2007}. Ultimately, one can plot the value of $f$ against the identifiable parameter combinations, when they are determined.

\subsection*{Sufficient summary plots}
A related figure is the \textit{sufficient summary plot}, which is a plot of a single QOI versus one row of $Q_a^T\theta$, that is, an informative linear combination of parameters/input variables~\cite{constantine2015active}. This plot can be conceptualized as a rotation of a surface plot of the function $q=f(\theta)$ to reveal a (potentially) lower dimensional structure by viewing it edge-on.

In practice, one may generate points by sampling from the parameter space (e.g., Latin hypercube sampling or using the density $\rho(\theta)$) and computing $q$ at each point~\cite{constantine2015active}. Because each row of $Q_a^T$ is a linearized identifiable combination, the relationship between $Q_a^T$ and $q$ should be nearly one-dimensional. A sufficient summary plot is typically used as validation tool to confirm that the active subspace adequately captures the desired variation in the data. The sufficient summary plot can also be generalized to include multiple rows of $Q_a^T$ by using 3D plots or heatmaps. In general, if there is a sizable eigenvalue gap after the first eigenvalue, one might expect a single-row sufficient summary plot to capture a univariate trend, while if the gap comes after the second eigenvalue, a two-row plot may be useful.

\section*{Examples and case studies}
We illustrate the definitions and techniques with three simple, analytic examples where we begin to explore the strengths and weaknesses of the local and global techniques for linear and non-linear identifiable combinations. Then, we consider two real-world case-studies to highlight the importance of tailoring the technique to the scientific question. Code for each of these examples and case studies is provided on Github (\url{https://github.com/epimath/sFIM-param-reduction}).

\subsection*{Example 1: Linear identifiable combination}
The first example,
\begin{equation}
f(\theta_1,\theta_2)=\exp(\theta_1+\theta_2),
\end{equation}
has a linear, structural identifiable parameter combination. This two parameter function is univariate, so \textit{a priori} its inputs cannot be uniquely determined from its output. Because the identifiable parameter combination $\theta_1+\theta_2$ is linear, we will be able to reconstruct it with these linear techniques.

Here,
\begin{equation}
\chi=\begin{bmatrix} \exp(\theta_1+\theta_2) & \exp(\theta_1+\theta_2)
\end{bmatrix}
\end{equation}
and 
\begin{equation}
F=\begin{bmatrix} \exp(2(\theta_1+\theta_2)) & \exp(2(\theta_1+\theta_2))\\ \exp(2(\theta_1+\theta_2)) & \exp(2(\theta_1+\theta_2))
\end{bmatrix}.
\end{equation}

\begin{figure}[h!]
\centering
	\begin{subfigure}[b]{0.45\textwidth}
		\includegraphics[width=1\textwidth]{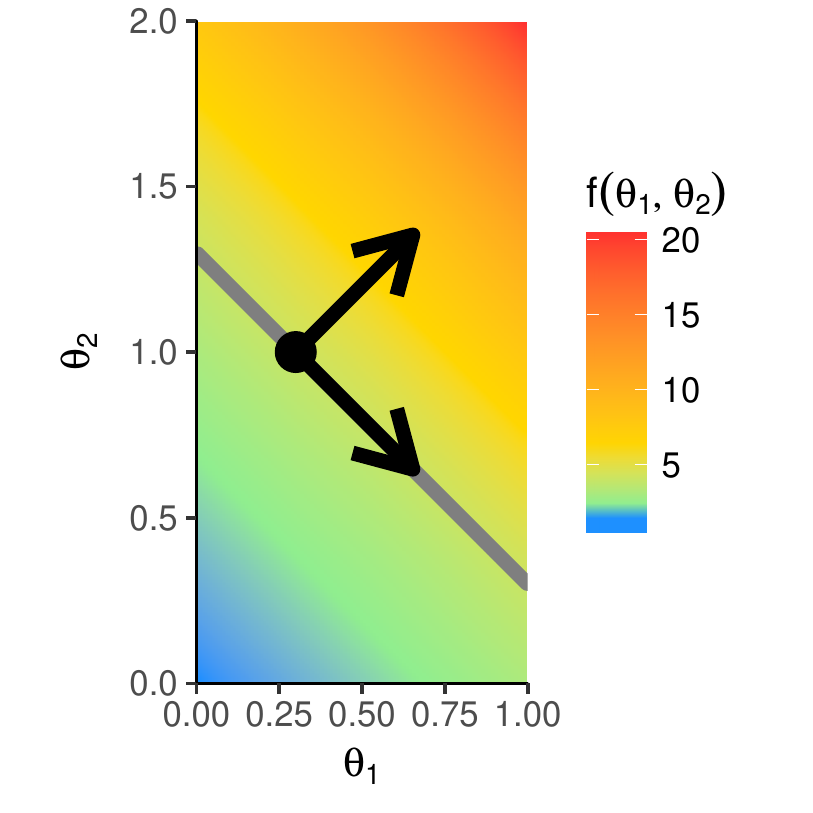}
		\caption{}
	\end{subfigure}
	\begin{subfigure}[b]{0.45\textwidth}
		\includegraphics[width=1\textwidth]{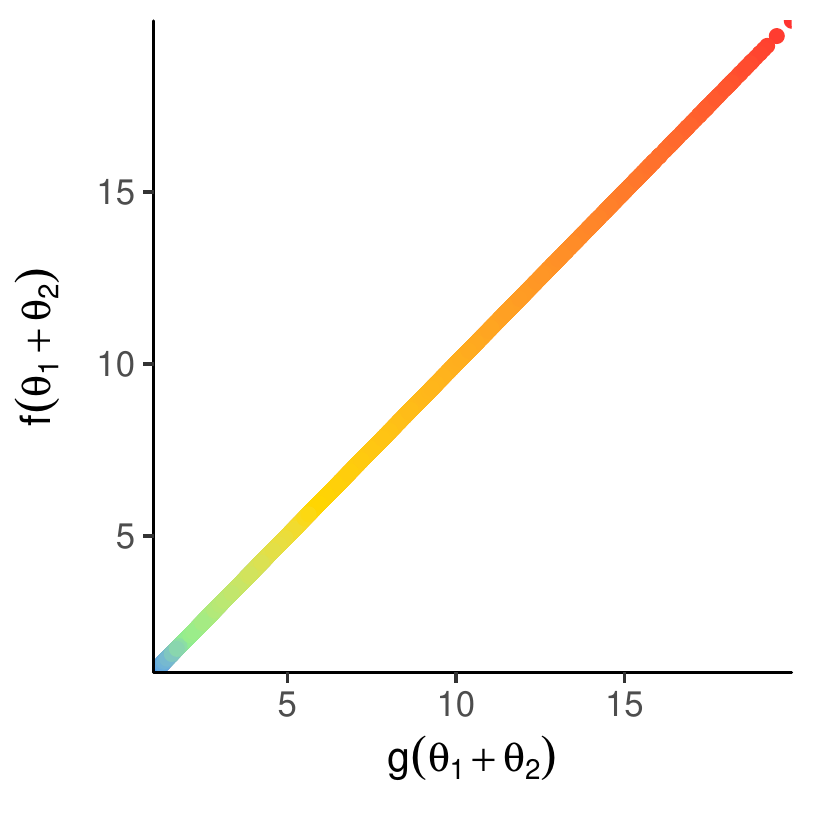}
		\caption{}
	\end{subfigure}
\caption{(a) Heat map of $f(\theta_1,\theta_2)=\exp(\theta_1+\theta_2)$ with the eigenvector directions of the sensitivity FIM evaluated at (0.3,1.0). The gray contour is the set of points $\{(\theta_1,\theta_2): f(\theta_1,\theta_2)=f(0.3,1.0)\}$. If we were to profile on $\theta_1$, we would find that this contour is also $\{(\theta_1^*,\text{argmin }  c(\theta_2)\}$, and since $\min c(\theta_2)=0$ along this profile, the parameter profile would be flat. (b) The approximation $g(\theta_1,\theta_2)$ is identical to $f(\theta_1,\theta_2)$ because the identifiable combination is linear. }
\label{Exp1}
\end{figure}

Because $f$ is univariate, the local sensitivity FIM $F$ can have rank at most 1. Indeed, the eigenvalues of $F$ are $\lambda_1=2\exp(2(\theta_1+\theta_2))$ and $\lambda_2=0$, with eigenvectors 
\begin{equation}
\nu_1=\frac{1}{\sqrt{2}}\begin{bmatrix} 1\\ 1 \end{bmatrix}, \quad \nu_2 =\frac{1}{\sqrt{2}}\begin{bmatrix} 1\\ -1 \end{bmatrix},
\end{equation}
as seen in Figure~\ref{Exp1}a. The rank deficiency of $F$ indicates that $f$ is not structurally identifiable, and $\nu_1$ correctly identifies the identifiable combination $\theta_1+\theta_2$.

To find the active subspaces, we must define a parameter region and density. Let us take $\theta_1$ and $\theta_2$ uniformly distributed on [0,1]$\times$[0,2]. The choice of domain here is meant to remove symmetry that could result in non-generalizable results.

Then
\begin{equation}
\begin{aligned}
C&=\int_0^2\int_0^1 \begin{bmatrix} \exp(2(\theta_1+\theta_2)) & \exp(2(\theta_1+\theta_2))\\ \exp(2(\theta_1+\theta_2)) & \exp(2(\theta_1+\theta_2))
\end{bmatrix} \,d\theta_1\,d\theta_2,\\
&=\frac{1}{4}(e^2-1)(e^2+1) \begin{bmatrix}  1 & 1\\ 1 & 1
\end{bmatrix}.
\end{aligned}
\end{equation}

Like the (local) $F$, the average sensitivity FIM $C$ is not full rank. Because the identifiable parameter combination is linear, the local and global techniques identify the same directions; the stiff/active/identifiable direction corresponds to the parameter combination $\theta_1+\theta_2$, given by $\nu_1$, while kernel is spanned by $\nu_2$, the sloppy/inactive/unidentifiable direction, corresponding to compensation between $\theta_1$ and $\theta_2$. Moreover, the approximation,
\begin{equation}
\begin{aligned}
g(\theta_1,\theta_2)&= f\left(Q_aQ_a^T\begin{bmatrix} \theta_1\\\theta_2\end{bmatrix}\right),\\
&= f\left(\frac{1}{2}\begin{bmatrix} 1\\1\end{bmatrix}\begin{bmatrix} 1 & 1\end{bmatrix}\begin{bmatrix} \theta_1\\\theta_2\end{bmatrix}\right),
\end{aligned}
\end{equation}
is in fact equal to $f(\theta_1,\theta_2)$ (Figure~\ref{Exp1}b).

Here, with a linear parameter combination, the active subspace gives the same answer as the local identifiability analysis. Both identify the linear, structural parameter combination, and we can successfully make a low-rank approximation of $f$.

\clearpage
\subsection*{Example 2: Non-linear identifiable combination}

Now, consider an example with a structural, non-linear identifiable parameter combination,
\begin{equation}
f(\theta_1,\theta_2)=\exp(\theta_1\theta_2).
\end{equation}

Again, this function is necessarily not identifiable because the number of parameters ($n=2$) is greater than the number of outputs  ($m=1$).

Here,
\begin{equation}
\chi=\begin{bmatrix} \theta_2\exp(\theta_1\theta_2) & \theta_1\exp(\theta_1\theta_2)
\end{bmatrix}
\end{equation}
and 
\begin{equation}
F=\begin{bmatrix} \theta_2^2\exp(2\theta_1\theta_2) & \theta_1\theta_2\exp(2\theta_1\theta_2)\\ \theta_1\theta_2\exp(2\theta_1\theta_2) & \theta_1^2\exp(2\theta_1\theta_2)
\end{bmatrix}.
\end{equation}

\begin{figure}%[h!]
\centering
	\begin{subfigure}[b]{0.45\textwidth}
		\includegraphics[width=1\textwidth]{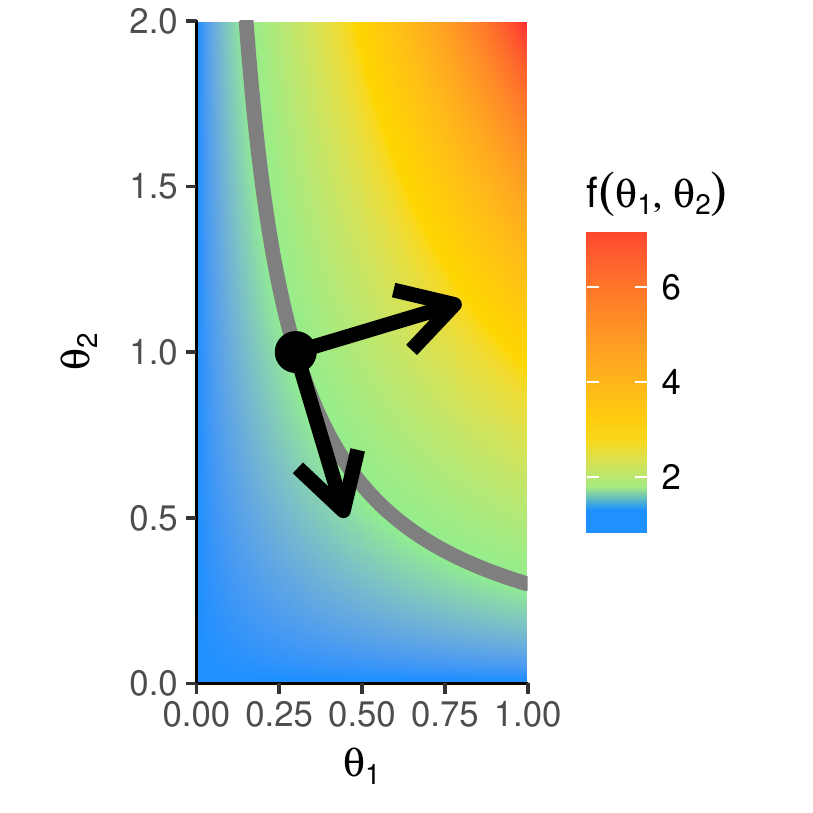}
		\caption{}
	\end{subfigure}
	\begin{subfigure}[b]{0.45\textwidth}
		\includegraphics[width=1\textwidth]{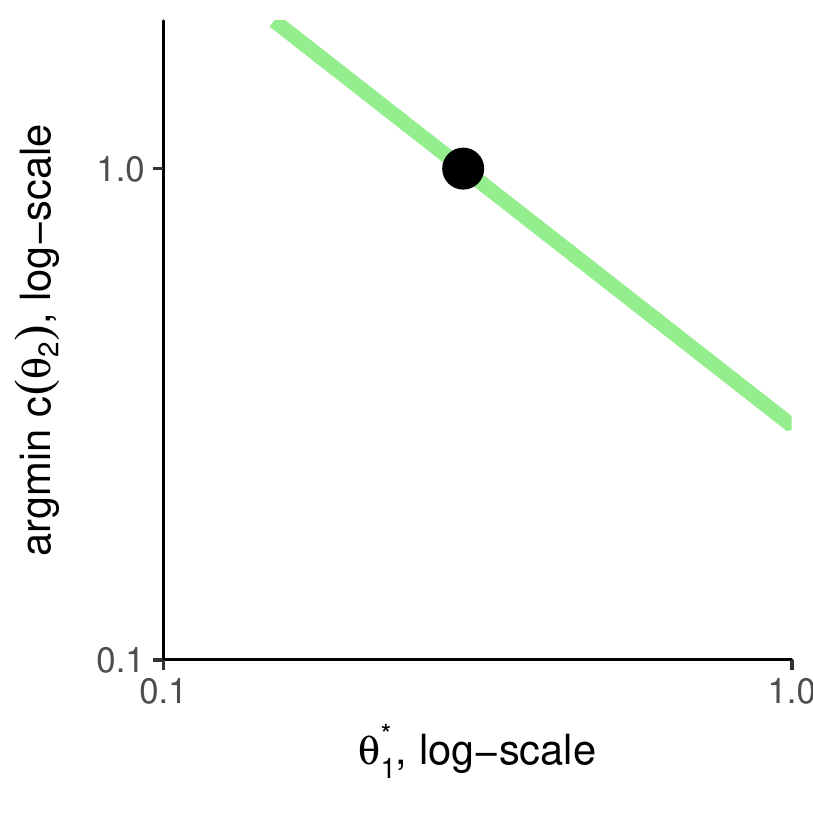}
		\caption{}
	\end{subfigure}
    \begin{subfigure}[b]{0.45\textwidth}
		\includegraphics[width=1\textwidth]{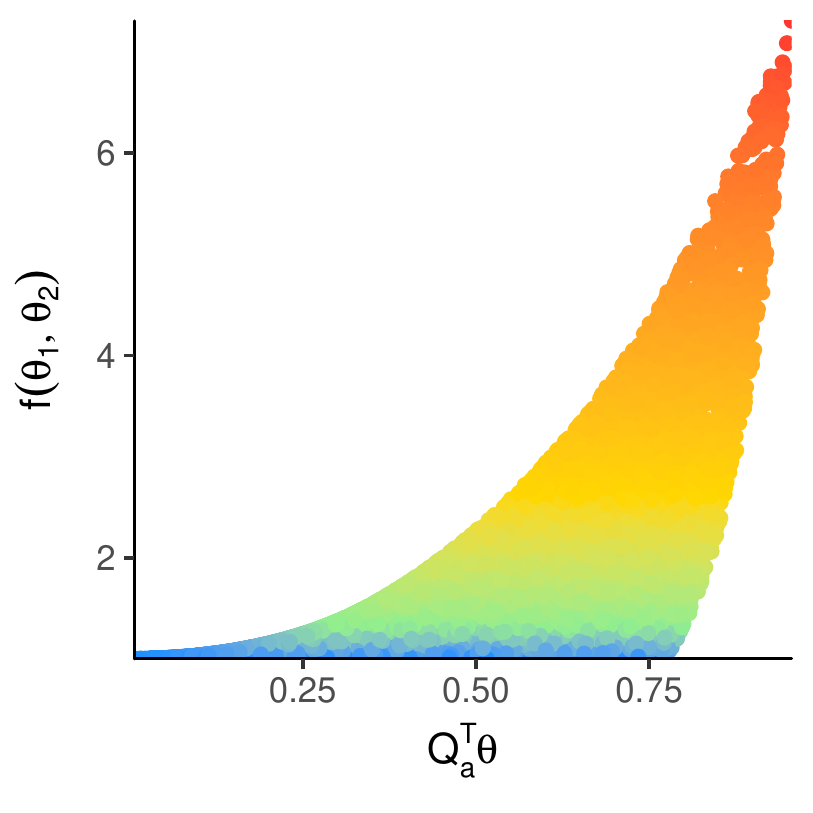}
		\caption{}
	\end{subfigure}
       \begin{subfigure}[b]{0.45\textwidth}
		\includegraphics[width=1\textwidth]{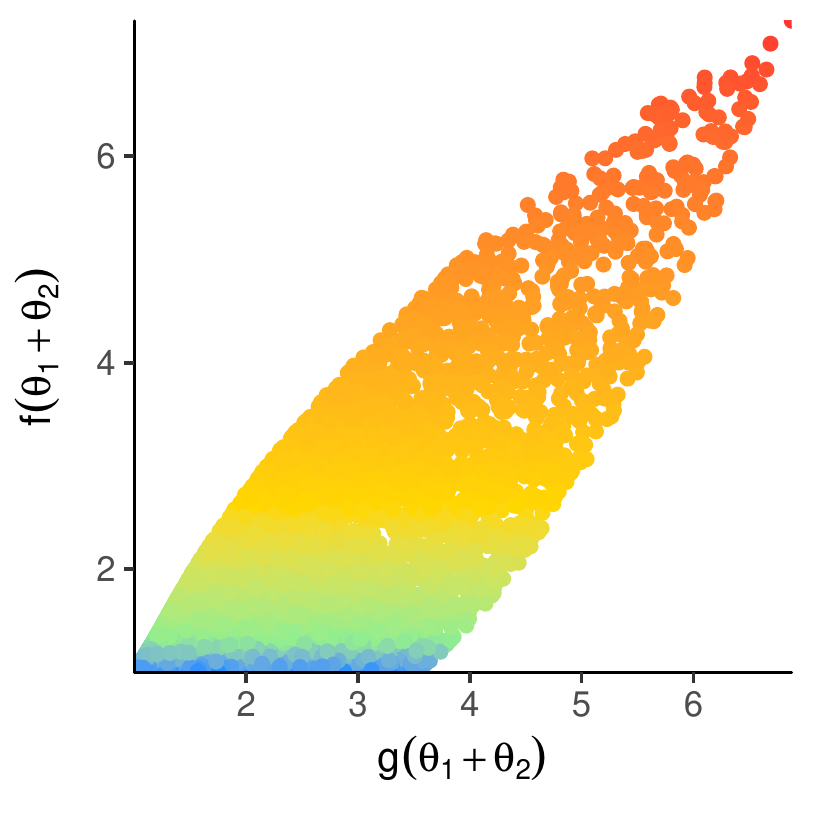}
		\caption{}
	\end{subfigure}
\caption{(a) Heat map of $f(\theta_1,\theta_2)=\exp(\theta_1+\theta_2)$ with the eigenvector directions of the sensitivity FIM evaluated at (0.3,1.0). The gray contour is the set of points $\{(\theta_1,\theta_2): f(\theta_1,\theta_2)=f(0.3,1.0)\}$. As in Example 1, if we were to profile on $\theta_1$, we would find that this contour is also $\{(\theta_1^*,\text{argmin }  c(\theta_2)\}$, and since $\min c(\theta_2)=0$ along this profile, the likelihood profile would be flat. (b) We plot this contour  on the log-log scale, and the linear relationship indicates that the parameter combination is a product. The value of $\theta_2$ minimizes the squared difference of $f(\theta_1^*,\theta_2)$ and $f(0.3,1)$. (c) Sufficient-summary plot of the first eigenvector of the average sensitivity FIM. A nearly linear relationship would indicate that a lower dimensional structure is present. (d) The approximation $g(\theta_1,\theta_2)$ deviates from $f(\theta_1,\theta_2)$ because the identifiable combination is non-linear.}
\label{Exp2}
\end{figure}

Again, for any $(\theta_1,\theta_2)$, the local sensitivity $F$ is not full rank. It has eigenvalues $\lambda_1=(\theta_1^2+\theta_2^2)\exp(2\theta_1\theta_2)$ and $\lambda_2=$0 with corresponding eigenvectors

\begin{equation}
\nu_1=\begin{bmatrix} \theta_2\\ \theta_1 \end{bmatrix}, \quad \nu_2 =\begin{bmatrix} \theta_1\\ -\theta_2 \end{bmatrix},
\end{equation}
as seen in Figure~\ref{Exp2}a.

For any particular $\left(\hat\theta_1,\hat\theta_2\right)$, this analysis will suggest that $(\hat\theta_2/\hat\theta_1)\theta_1 + \theta_2$ is a (linearized) identifiable combination. In this situation, we can recover the true identifiable combination through profiling. That is, we fix $\theta_1$ at a series of values $\theta_1^*$ and determine the value of $\theta_2$ such that

\begin{equation}
c(\theta_2)=(f(\theta_1^*,\theta_2)-f(\hat\theta_1,\hat\theta_2))^2
\end{equation}

is minimized. In this case, profiling reveals a linear relationship between $\theta_1^*$ and $\arg\min c(\theta_2)$ when plotted on a log-log scale (Figure~\ref{Exp2}b),  demonstrating that 

\begin{equation}
\arg\min c(\theta_2) \propto \frac{1}{\theta_1^*},
\end{equation}

i.e. $\theta_1\theta_2$ is an identifiable combination.

Now, we consider whether we can create a one-dimensional global approximation for $f$ using the average sensitivity FIM. We again assume that $\theta_1$ and $\theta_2$ are uniformly distributed on [0,1]$\times$[0,2]. Then
\begin{equation}
\begin{aligned}
C&=\int_0^2\int_0^1 \begin{bmatrix} \theta_2^2\exp(2\theta_1\theta_2) & \theta_1\theta_2\exp(2\theta_1\theta_2)\\ \theta_1\theta_2\exp(2\theta_1\theta_2) & \theta_1^2\exp(2\theta_1\theta_2)
\end{bmatrix} \,d\theta_1\,d\theta_2,\\
&=\begin{bmatrix}  4.89983 & 8.9827\\ 8.9827 & 19.5993
\end{bmatrix}.
\end{aligned}
\end{equation}

The average sensitivity FIM $C$ has the eigendecomposition

\begin{align}
C&= \begin{bmatrix} 0.428222 & -0.903673 \\ 0.903673 & 0.428222
\end{bmatrix}\begin{bmatrix} 23.8559 & 0 \\ 0 & 0.64321
\end{bmatrix} \begin{bmatrix} 0.428222 & 0.903673 \\ -0.903673 & 0.428222
\end{bmatrix}
\end{align}

In this calculation, $C$ is full rank, but there is a small eigenvalue gap between $\lambda_1=$23.8559 and $\lambda_2=$0.64321, raising the possibility of a one-dimensional approximation. The active subspace is the span of 
\begin{equation}
 Q_a=\begin{bmatrix} 0.428222  \\ 0.903673
\end{bmatrix}. 
\end{equation}
The sufficient summary plot for the active subspace (Figure~\ref{Exp2}c) indicates that, although there is some sort of structure, most of the variance is not captured by the active subspace alone. Indeed, the approximation
\begin{equation}
\begin{aligned}
g(\theta_1,\theta_2)&=f\left(Q_aQ_a^T\begin{bmatrix}\theta_1  \\ \theta_2
\end{bmatrix}\right)\\
&=f(0.183374 \theta_1 + 0.386973 \theta_2, 2.11029\,(0.183374 \theta_1 + 0.386973 \theta_2)),
\end{aligned}
\end{equation}
deviates a great deal from $f$ (Figure~\ref{Exp2}d).

In contrast to Example 1, both the local and global sensitivity FIM analyses had difficulty because we are applying linear methods to a non-linear problem.
Profiling can determine the form of the non-linear combination, locally, but the global analysis was unable to find the low-dimensional approximation.

\clearpage
\subsection*{Example 3: Nearly rank deficient sFIM}

Now, we consider an example with a nearly rank-deficient sensitivity FIM. Let
\begin{equation}
f(\theta_1,\theta_2)=(\theta_1+\log(1+\theta_2),~\theta_1+\theta_2).
\end{equation}

This function has a two-dimensional output, so, unlike the previous examples, it could be identifiable (and in fact is structurally identifiable because $f$ is injective). But, because $\log(1+\theta_2)\approx\theta_2$ for $\theta_2$ near zero, we expect to find that $\theta_1+\theta_2$ is a practically identifiable combination when $\theta_2$ is small.

Computing,
\begin{equation}
\chi=\begin{bmatrix} 1 & 1 \\ \frac{1}{1+\theta_2} & 1
\end{bmatrix},
\end{equation}

and 
\begin{equation}
F=\begin{bmatrix}1 + \frac{1}{(1+\theta_2)^2} & 1+\frac{1}{1+\theta_2}\\ 1+\frac{1}{1+\theta_2} & 2\\
\end{bmatrix}.
\end{equation}

The ratio of the eigenvalues of this matrix varies widely depending on the value of $\theta_2$. Let us consider two points $(\theta_1,\theta_2)=$(1.95,0.05) and (0.05,1.95). At the first point, the eigenvalues are 
$\lambda_1=$3.91 and $\lambda_2=$5.80E-4, with eigenvectors 
\begin{equation}
\nu_1=\begin{bmatrix} 0.699\\ 0.715 \end{bmatrix}, \quad \nu_2 =\begin{bmatrix} -0.715\\ 0.699 \end{bmatrix}.
\end{equation}

The matrix is nearly rank deficient at this point, and the eigenvectors suggest that we can only practically identify $\theta_1+\theta_2$ from the value of $f$ at this point. At (0.05,1.95), on the other hand, the eigenvalues are much closer,  with eigenvalues $\lambda_1=$2.96 and $\lambda_2=$0.15, indicating that the values of both $\theta_1$ and $\theta_2$ can be identified.

As a global measure, active subspaces will not identify local opportunities for dimension reduction. On the other hand, making system simplifications and creating a function approximation based on the local practical identifiability near the first point would not be useful for global approximation.

\pagebreak
\subsection*{Case study: Cell cycle model}

In this case study, we consider a cell cycle model developed by Gerard and Goldbeter~\cite{Gerard2011}: 
\begin{align}
\begin{split}
\frac{d\textrm{Md}}{dt} =& v_\textrm{sd}\cdot\left(\frac{\textrm{GF}}{K_\textrm{gf}+\textrm{GF}}\right) - V_\textrm{dd}\cdot\left(\frac{\textrm{MD}}{K_\textrm{dd}+\textrm{Md}}\right),\\
\frac{d\textrm{E2F}}{dt} =&V_\textrm{le2f}\cdot\left(\frac{(\textrm{E2F}_\textrm{tot}-\textrm{E2F})}{K_\textrm{le2f}+(\textrm{E2F}_\textrm{tot}-\textrm{E2F})}\right)\cdot(\textrm{Md}+\textrm{Me})\\&-V_\textrm{2e2f}\cdot\left(\frac{\textrm{E2F}}{K_\textrm{2e2f}+\textrm{E2F}}\right)\cdot \textrm{Ma}\\
\frac{d\textrm{Me}}{dt}=& v_\textrm{se}\cdot\textrm{E2F}- V_\textrm{de}\cdot\textrm{Ma}\cdot\left(\frac{\textrm{Me}}{K_\textrm{de}+\textrm{Me}}\right),\\
\frac{d\textrm{Ma}}{dt}=&v_\textrm{sa}\cdot\textrm{E2F}-V_\textrm{da}\cdot\textrm{Cdc20}\cdot\left(\frac{\textrm{Ma}}{K_\textrm{da}+\textrm{Ma}}\right),\\
\frac{d\textrm{Mb}}{dt}=&v_\textrm{sb}\cdot\textrm{Ma}-V_\textrm{db}\cdot \textrm{Cdc20}\cdot \left(\frac{\textrm{Mb}}{K_\textrm{dv}+\textrm{Mb}}\right),\\
\frac{d\textrm{Cdc20}}{dt}=&V_{1cdc20}\cdot\textrm{Mb}\cdot\left(\frac{(\textrm{Cdc20}_\textrm{tot}-\textrm{Cdc20})}{K_\textrm{1cdc20}+(\textrm{Cdc20}_\textrm{tot}-\textrm{Cdc20})}\right)\\
&-V_\textrm{2cdc20}\cdot\left(\frac{\textrm{Cdc20}}{K_\textrm{2Cdc20}+\textrm{Cdc20}}\right).
\end{split}
\end{align}

\begin{figure}
\centering
    \begin{subfigure}[b]{0.48\textwidth}
    \includegraphics[width=1\textwidth]{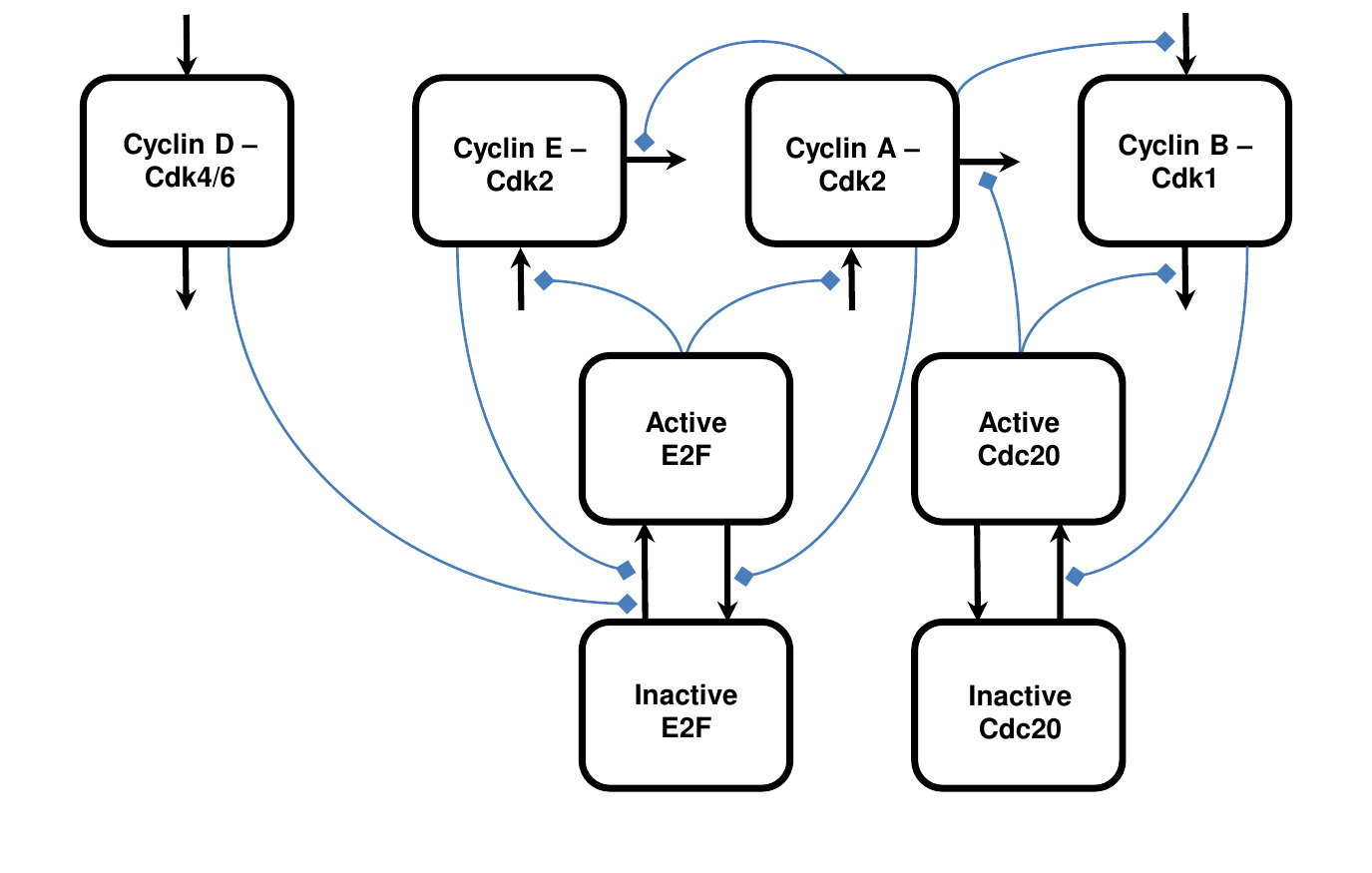}
		\caption{}
    \end{subfigure}
	\begin{subfigure}[b]{0.48\textwidth}
    \includegraphics[height=3in]{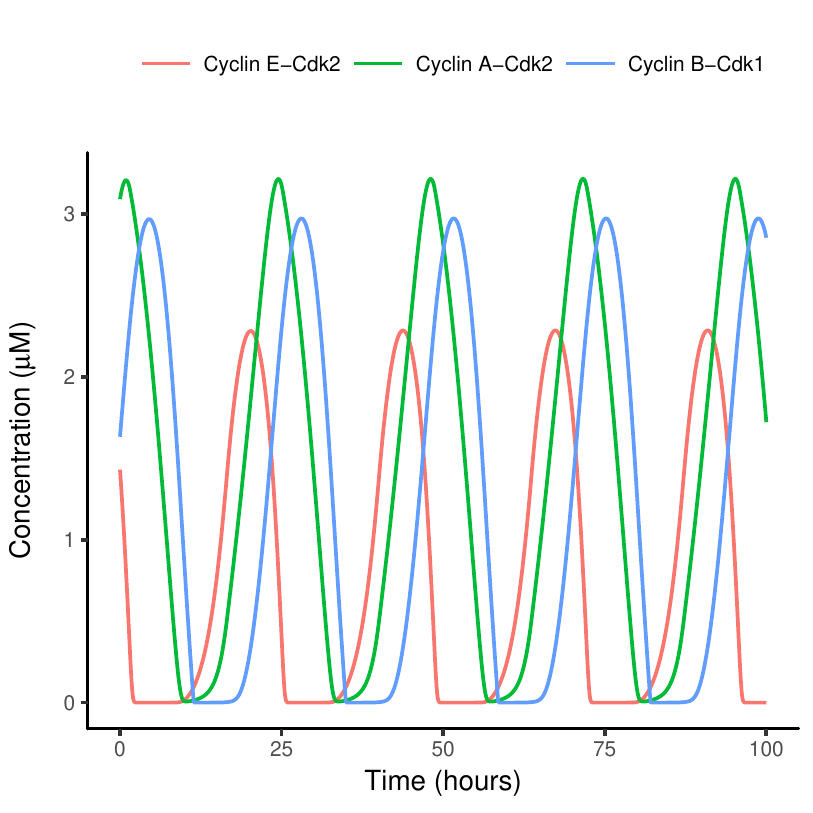}
		\caption{}
	\end{subfigure}

\caption{(a) A skeleton model of the cell cycle developed by Gerard and Goldbeter~\cite{Gerard2011}. (b) The cell progresses from G1 to S to G2 To M phases as cyclins E, A, and B---in complex with their cyclin dependent kinases (Cdk)---wax and wane periodically.} 
\label{cell_cycle}
\end{figure}

Reduced from their original model of thirty-nine variables, this skeleton model qualitatively reproduces cell cycle behavior in six variables and twenty-four parameters (Figure~\ref{cell_cycle}; see \cite{Gerard2011} for variable and parameter definitions). Cyclins are a family of proteins that, in complex with cyclin-dependent kinases (Cdk), drive a cell through the G1, S, G2, and M phases of the cell cycle. Transcription factor E2F and protein Cdc20 help regulate this cell cycle progression. Here, we consider one quantity of interest---the period of the cell cycle---and we ask whether we can find a lower dimensional structure in parameter space that predicts it. 

\begin{figure}[h!]
\centering
	\begin{subfigure}[b]{0.15\textwidth}	\includegraphics[height=3in]{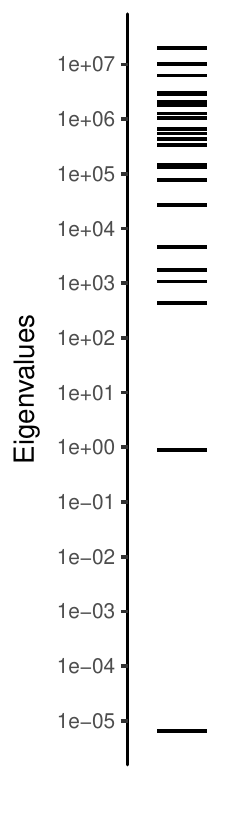}
		\caption{}
	\end{subfigure}
	\begin{subfigure}[b]{0.48\textwidth}
    \includegraphics[height=3in]{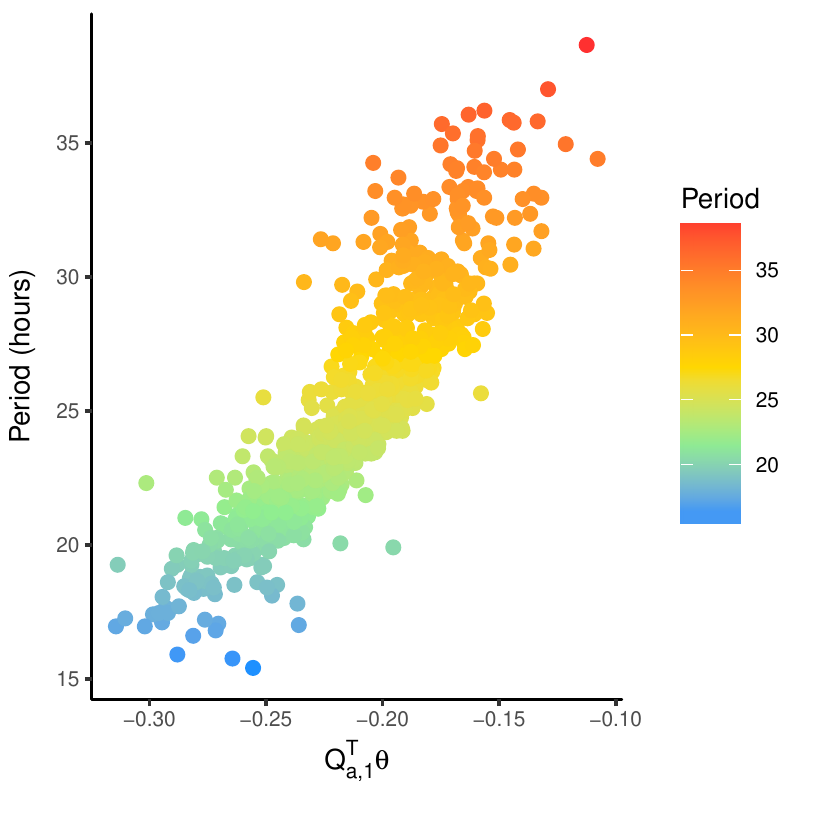}
		\caption{}
	\end{subfigure}
    	\begin{subfigure}[b]{0.48\textwidth}		\includegraphics[height=3in]{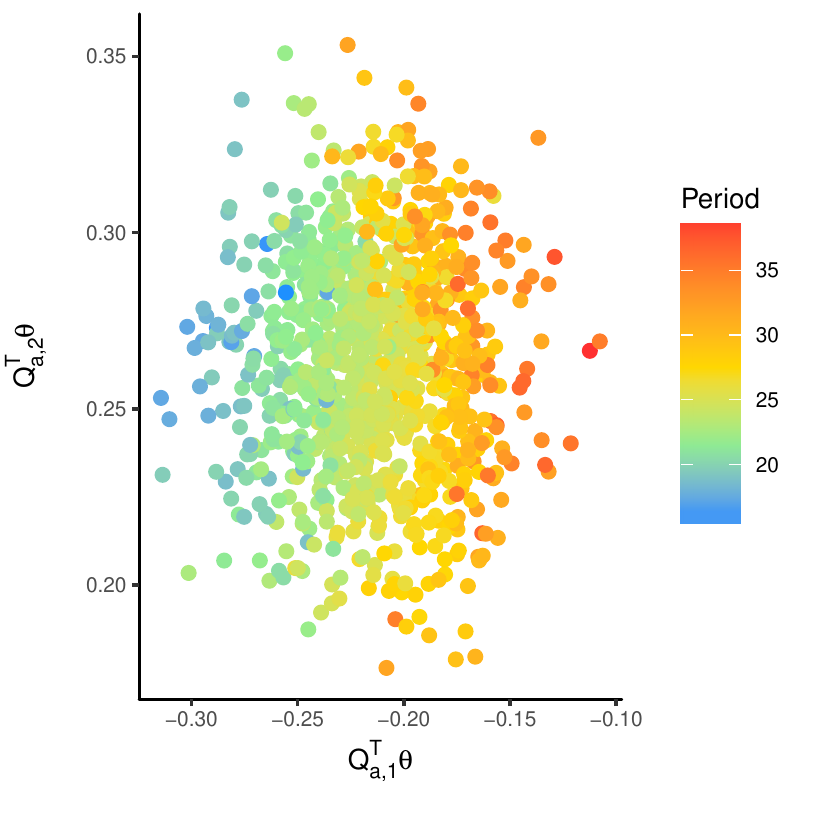}
		\caption{}
	\end{subfigure}
        \begin{subfigure}[b]{0.48\textwidth}		\includegraphics[height=3in]{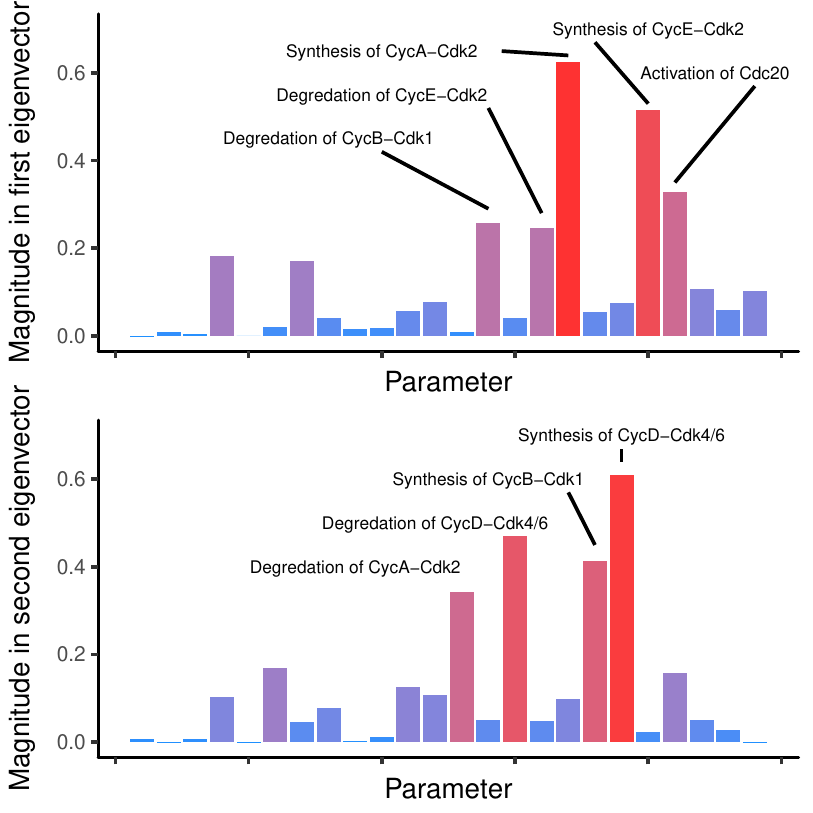}
		\caption{}
    \end{subfigure}
\caption{(a) Global sensitivity FIM analysis of the period of the cell cycle model finds that eigenvalues cluster near the top of the ladder, suggesting that only a few of the variables do not impact the period. Nevertheless, the sufficient-summary plot of the first (b) and first two eigenvectors (c) demonstrate that the period is determined by a low-dimensional structure in parameter space. (d) Analysis of the parameter loadings of the first two eigenvectors highlight the importance of the synthesis and degredation of the cyclin--Cdk complexes. }
\label{cell_cycle_results}
\end{figure}

Although our question is a global one, the model does not exhibit periodic dynamics everywhere in parameter space. Hence, we restrict our analysis to parameter values between 50--150\% of the default parameters of Gerard and Goldbeter~\cite{Gerard2011}. We sample 1,000 points from this restricted parameter space and calculate the period and estimate the gradient at each point. We compute the global sensitivity FIM $C$. Although the eigenvalues of $C$ do not display a large eigenvalue gap at the top of the eigenvalue ladder (Figure~\ref{cell_cycle_results}(a)), sufficient-summary plots of the first (Figure~\ref{cell_cycle_results}(b)) and first two (Figure~\ref{cell_cycle_results}(c)) eigenvectors demonstrate that most of the variance in the period is controlled by a subset of parameters. Considering the parameter loadings of these eigenvectors (Figure~\ref{cell_cycle_results}(d)), we see that the synthesis and degradation of the cyclin-Cdk complexes have the greatest effect on the period.

Exact computation of identifiable combinations for observing output trajectories becomes increasing computationally intensive for even moderately sized models. Moreover, there is no clear way to formulate the period as a rational function of an output trajectory, as would be required for differential algebra method. Parameter profiling here would be computationally intensive, and it may be difficult glean useful information about the  likely complex practical identifiable combinations without an involved analysis fixing different combinations of parameters.
Local sensitivity FIM may be useful, but, as a first pass it will be insufficient for understanding which parameters are important over a wide range of parameters. Global sensitivity FIM analysis (active subspaces) is useful in this instance because we i)  have  a single quantity of interest that cannot be analytically expressed as a function of the output trajectories and ii) are  interested in determining which parameters would be needed to develop a low-dimensional, computationally fast approximation to the period that does not require solving a system of ODEs or numerically estimating the period.

\clearpage
\subsection*{Case study: infectious disease transmission model}
\begin{figure}
\centering
\includegraphics[width=1\textwidth]{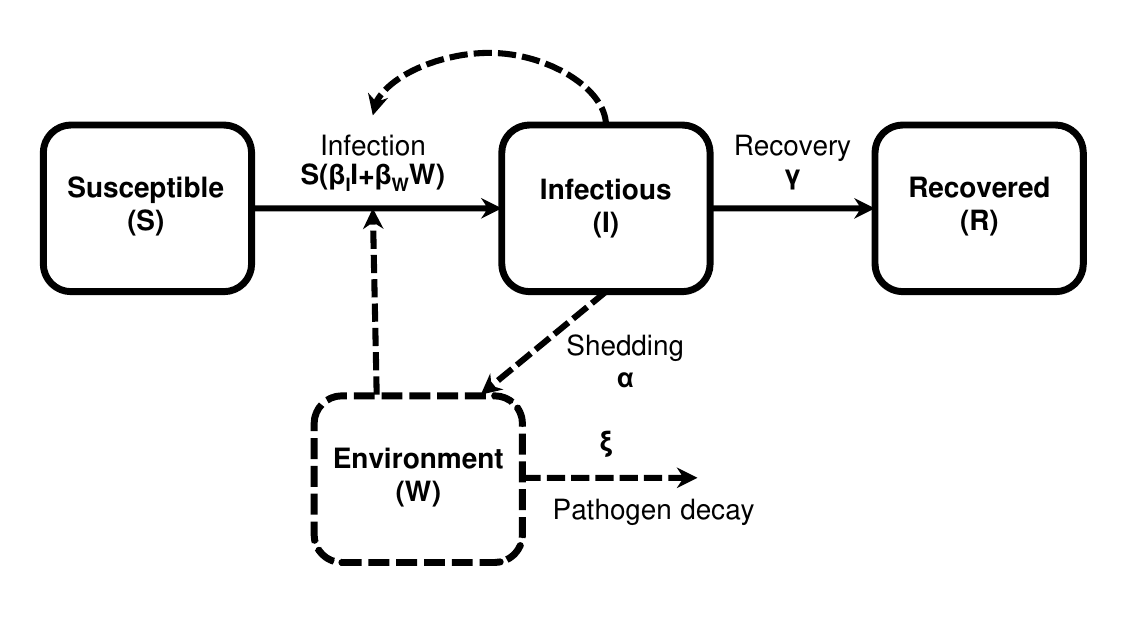}
\caption{The Susceptible, Infectious, Water, Recovered (SIWR) model of infectious disease transmission with direct and indirect transmission pathways.}
\label{SIWR}
\end{figure}
We finally illustrate sensitivity FIM techniques on a ordinary differential equation model of infectious disease transmission with two pathways: a direct (person-to-person) route and an indirect (environmental) route, shown in Figure~\ref{SIWR}. This example will allow us to illustrate two points: how different identifiable combinations (active/inactive directions) may be present in different parts of parameter space, and how one might examine vectors of QOIs in the form of a time series. 

The model we consider is often referred to as the SIWR model as an initialism of the compartments: susceptible, infectious, water, and recovered~\cite{Tien2010,Eisenberg2013a}. Infectious people infect susceptible people directly with rate $\beta_I$, recover from infection at rate $\gamma$, and pathogens in the water infect susceptible people with rate $\beta_W$ and decay at rate $\xi$ (we note also that $W$ has been re-scaled by the pathogen shedding rate $\alpha$ and decay rate $\xi$, which is why $\alpha$ does not appear in the equations, see~\cite{Tien2010,Eisenberg2013a}). The units for all human compartments ($S$, $I$, and $R$) are assumed to be as fractions of the total population at risk. A schematic of this system shown in Figure~\ref{SIWR} and the system of equations is
\begin{align}\label{eq:siwr}
\begin{split}
\dot S &= -S(\beta_I I + \beta_W W),\\
\dot I &= S(\beta_I I + \beta_W W) - \gamma I,\\
\dot R &= \gamma I,\\
\dot W &= \xi (I-W).
\end{split}
\end{align}
For our QOI/output, we take in this case the vector of measured cases over time, $y = \kappa I$, where $\kappa$ is the reporting rate multiplied by the size of the at-risk population. Following Eisenberg et al.~\cite{Eisenberg2013a}, we define $k = 1/\kappa$ and work with $k$ (so that $y = I/k$), as it is often easier to estimate ($k$ is bounded generally between 0 and 1, while $\kappa$ is typically a large number ranging from the hundreds to the millions). 

Using the differential algebra method, Eisenberg et al.~\cite{Eisenberg2013a} previously showed that the scaled form of the model given in Eq.~\eqref{eq:siwr} is structurally identifiable, and used the model to demonstrate that both a direct and indirect pathway was needed to explain the observed transmission dynamics for the 2006 cholera epidemic in Angola. However, while the model is structurally identifiable, they also observed that there are often issues of practical unidentifiability between $\beta_W$ and $\xi$~\cite{Eisenberg2013a} when noisy data is considered. Further, in the limit as $\xi\rightarrow\infty$, $\beta_W$ and $\beta_I$ become indistinguishable, forming an identifiable combination $\beta_W+\beta_I$. As this happens, the model may also become insensitive to relatively small changes in $\xi$. In practice, these issues mean that depending on where one is in parameter space and the data quality (level of noise/variance, frequency of samples), the model may be practically identifiable, or unidentifiable in different ways depending on which dependencies between $\xi$, $\beta_W$, and $\beta_I$ dominate.

\begin{figure}[t]
\centering
	\begin{subfigure}[b]{0.32\textwidth}
    \includegraphics[height=2.2in]{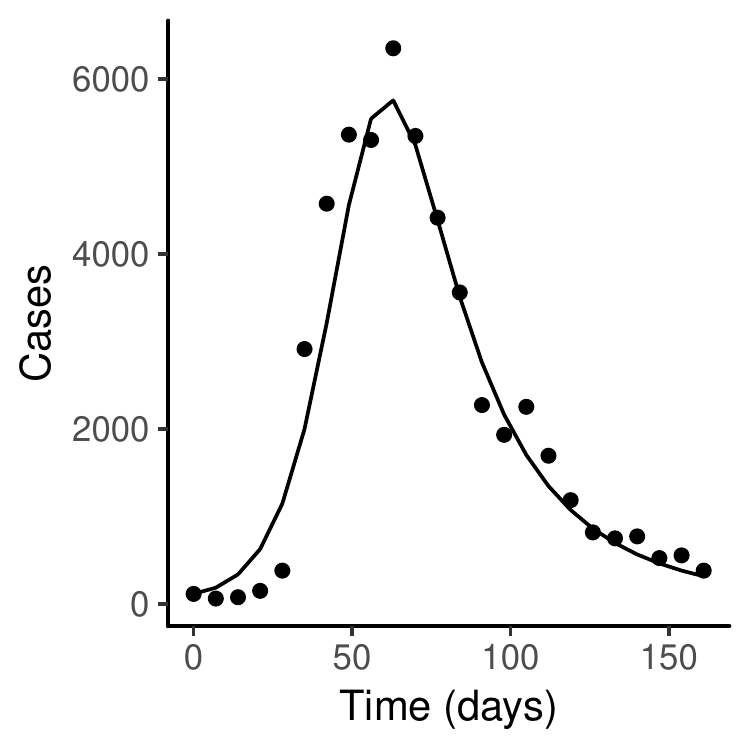}
		\caption{}
	\end{subfigure}
	\begin{subfigure}[b]{0.33\textwidth}
    \includegraphics[height=2.2in]{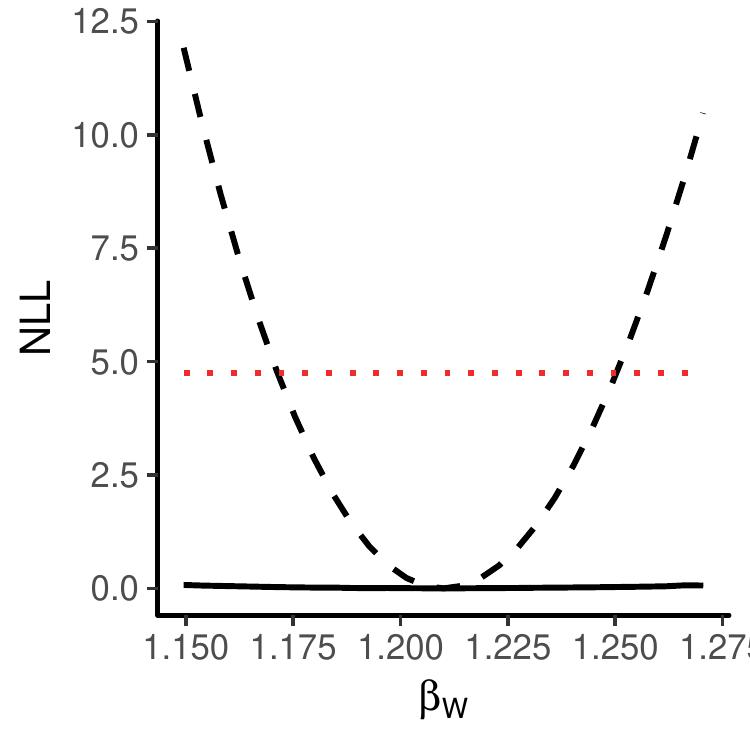}
		\caption{}
	\end{subfigure}
    	\begin{subfigure}[b]{0.33\textwidth}		
	\includegraphics[height=2.2in]{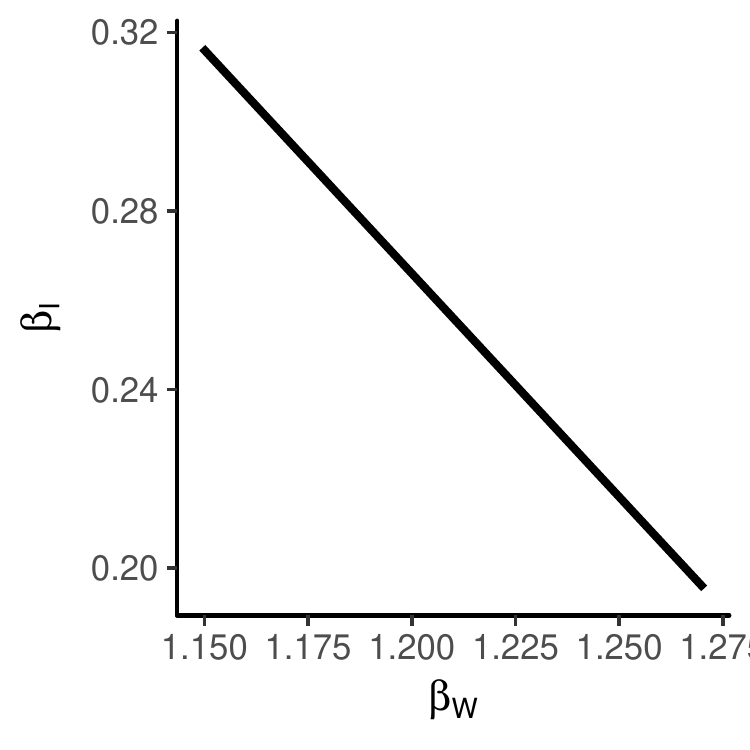}
		\caption{}
	\end{subfigure}
\caption{\textbf{(a)} Model fit to data from the 2006 Angola epidemic, using parameters estimated in \cite{Eisenberg2013a}. \textbf{(b)} Profile likelihoods for $\beta_W$ using simulated, noise-free data generated in two parameter regimes: dashed line, the parameters used to generate panel (a) (see example text for values), and solid line, the same parameters but with $\xi$ multiplied by 5000. The red, dotted line indicates the threshold for the 95\% confidence interval. Once $\xi$ is large, the model becomes unidentifiable (flat profile with an infinite confidence interval). \textbf{(c)} Corresponding parameter relationship plot, showing how $\beta_I$ varies as $\beta_W$ is profiled in the unidentifiable case for panel (b). As $\xi$ becomes large, the two transmission pathways can no longer be distinguished, with the total $\beta_I + \beta_W$ forming an identifiable combination. In active subspaces parlance, the direction along the line would be considered inactive/sloppy, while the normal direction to the line (corresponding to changing the value of the identifiable combination) would be considered active/stiff.}
\label{SIWR_profiles}
\end{figure}

We implemented the model using the parameter estimates given in~\cite{Eisenberg2013a} for the 2006 Angola epidemic (shown in Figure~\ref{SIWR_profiles}a, with parameter set $\theta = \{ \beta_I=0.256, \beta_W=1.21, \xi = 0.00756, k = 1.1212e-5\}$). As in~\cite{Eisenberg2013a}, we let $\gamma = 0.25$ be a fixed (not estimated or varied) value. We let the initial conditions be determined by the data, taking $I(0) = y(0)k$, and set the remainder of the population to be susceptible. To simplify our example somewhat, we use simulated data without noise (this can also potentially be useful in illustrating the more general aim of parameter reduction from a time series QOI, even if not working with data).

To illustrate the issues of unidentifiability as $\xi \rightarrow \infty$, we profiled the parameters in both regimes---we first profile using the above parameters as our $\hat{\theta}$ (i.e. the model trajectory $y(\hat{\theta})$ is treated as the `data' in the cost function and used for profiling). We then profiled the model with our nominal value for  $\xi$ 5000 times higher, equal to $37.8$. In both cases, we took 20 data points, spread evenly between the start and approximate end of the epidemic. As shown in Figure~\ref{SIWR_profiles}b, the model becomes unidentifiable, with an combination that has begun to approximate $\beta_W+\beta_I$ (Figure~\ref{SIWR_profiles}c). Indeed, as $\xi$ increases, the rank of the sFIM falls from full rank of 4, to 3, and finally to 2 once $\xi$ is approximately four orders of magnitude larger, indicating increasing dependencies between parameters. 

Next, we consider how we might apply more general sFIM techniques beyond rank, given that we are considering a vector of QOIs (the time series $y$). We can evaluate the eigenvalues and eigenvectors of $F$ and $C$ in the vector case, however the vector form of our QOI $y$ makes examining sufficient summary plots (which are typically defined for scalar QOI) more complicated. However, as we are partially localizing our analyses, we will define a single summary QOI, $q$, to be the cost function minimizing the sum of squares ($\mathcal{L}^2$ norm) between the current parameters and our trajectory at a set of nominal parameter values at the center of our parameter regime, denoted $\hat{\theta}$. 
This allows us to examine which directions in parameter space tend to maintain the same goodness of fit or the same model behavior as the nominal values (inactive/sloppy) versus tend to alter the model behavior/fit. However, we note that once away from $\hat{\theta}$ there will be many ways to attain the same cost function value that do not necessarily represent the same model behavior. Still, this QOI may allow us to use sufficient summary plots to examine the overall impact of the parameters across the time series. 

We restrict our analysis to parameter values between 50-150\% of $\hat{\theta}$, sampling 500 points using Latin hypercube sampling (LHS). We note that while it is often preferable for sufficient summary plots to rescale the parameters to be within $(-1,1)$ \cite{constantine2015active}, we found that the identifiability relationships in the eigenvectors were better visualized if we took symmetric ranges but left the parameters unscaled (this was likely due to the fact that the identifiable combination sum between the two parameters is more easily captured without any multiplicative scaling on the parameters). 

\begin{figure}%[h!]
\centering
$F_y$\\
	\begin{subfigure}[b]{0.09\textwidth}	
	\includegraphics[height=1.5in]{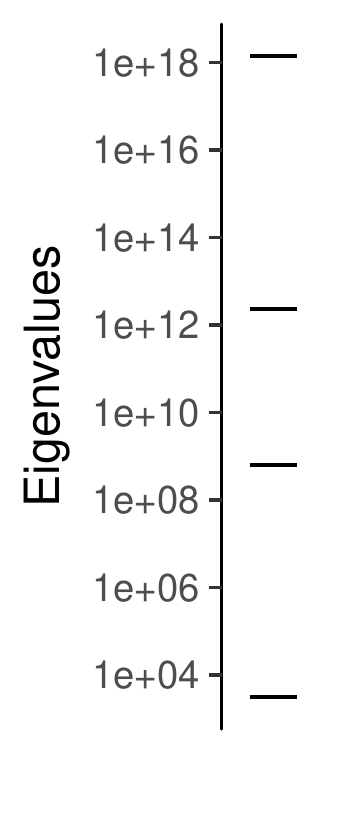}
%		\caption{}
	\end{subfigure}
	\begin{subfigure}[b]{0.22\textwidth}	
	\includegraphics[width = \textwidth]{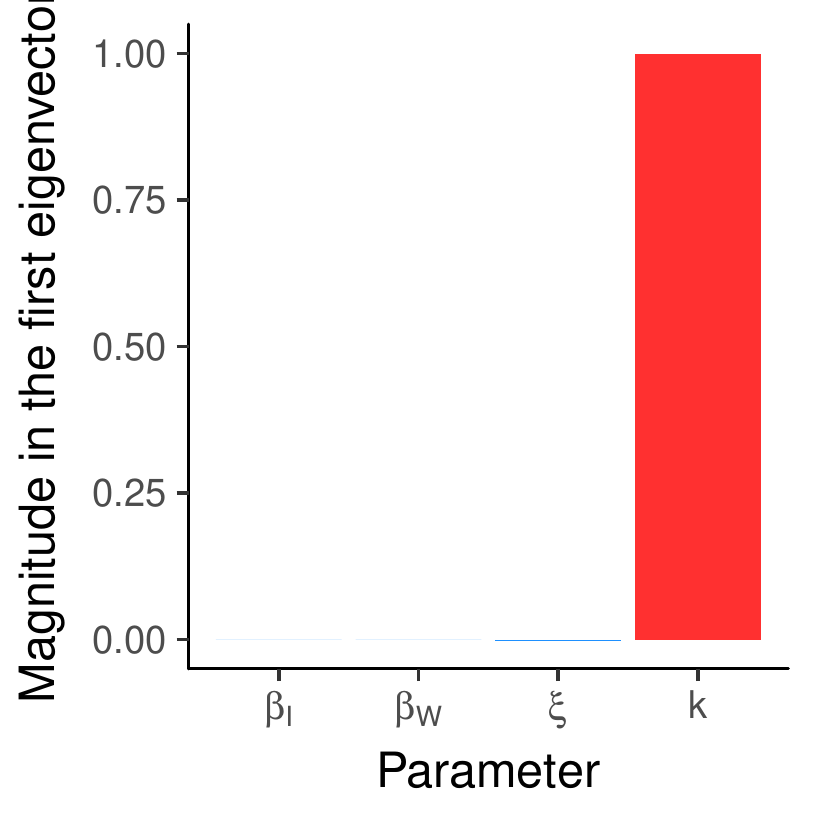}
	\end{subfigure}
	\begin{subfigure}[b]{0.22\textwidth}
	\includegraphics[width = \textwidth]{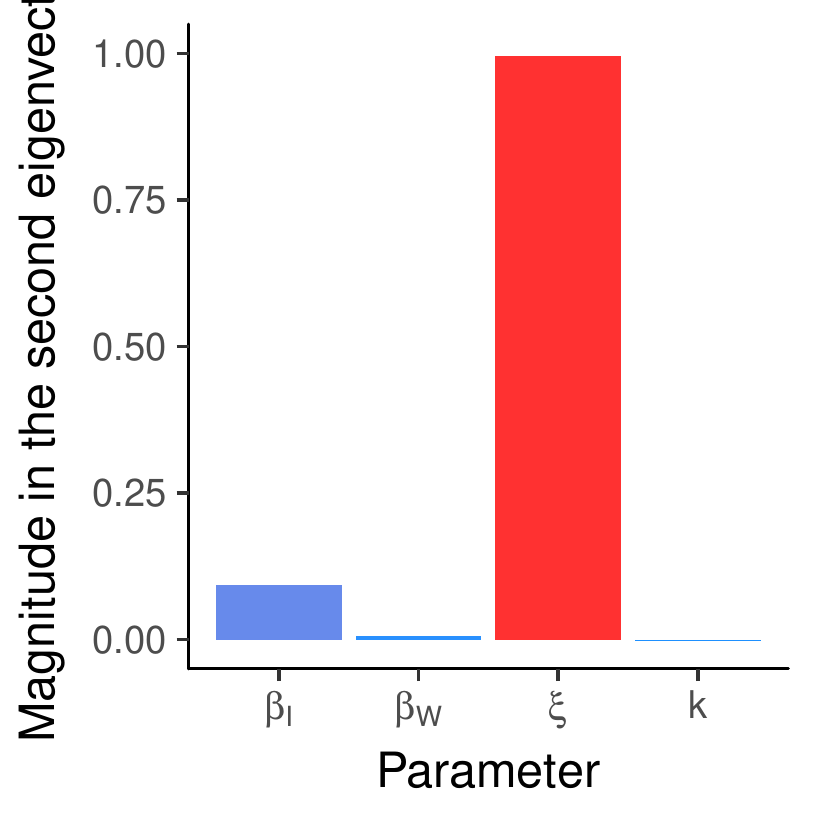}
%		\caption{}
	\end{subfigure}
	\begin{subfigure}[b]{0.22\textwidth}	
	\includegraphics[width = \textwidth]{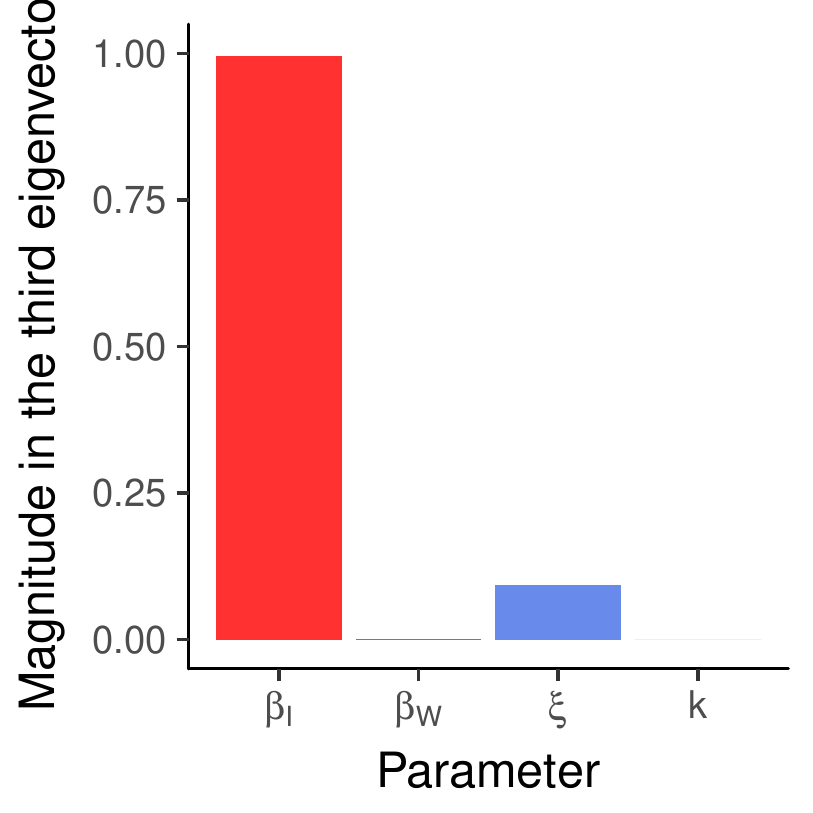}
	\end{subfigure}
	\begin{subfigure}[b]{0.22\textwidth}	
	\includegraphics[width = \textwidth]{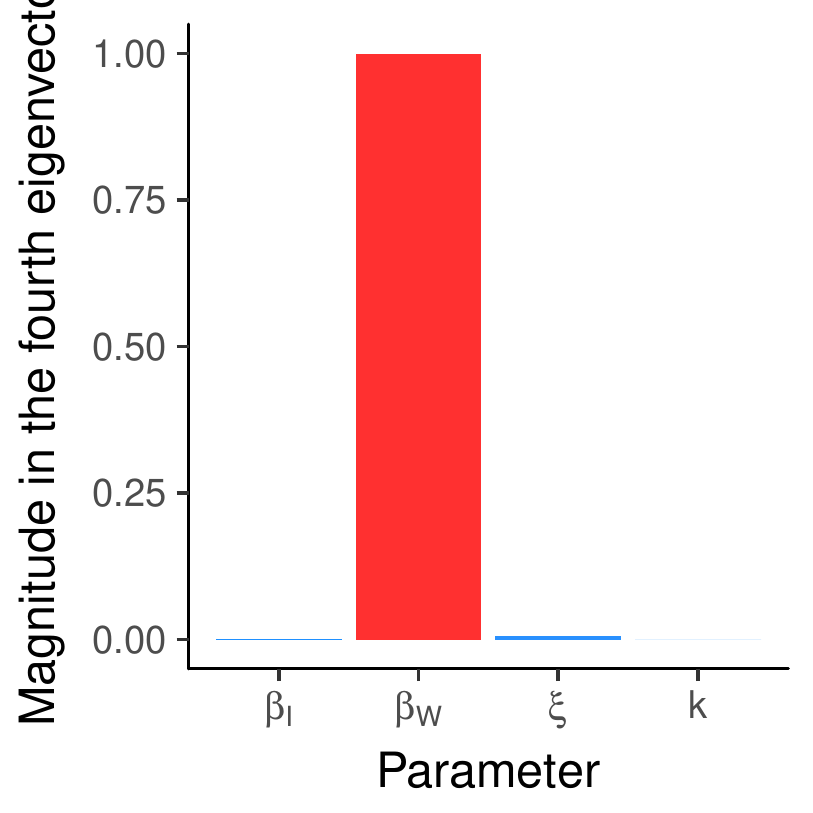}
%		\caption{}
	\end{subfigure}
\\
$F_q$\\
	\begin{subfigure}[b]{0.09\textwidth}	
	\includegraphics[height=1.5in]{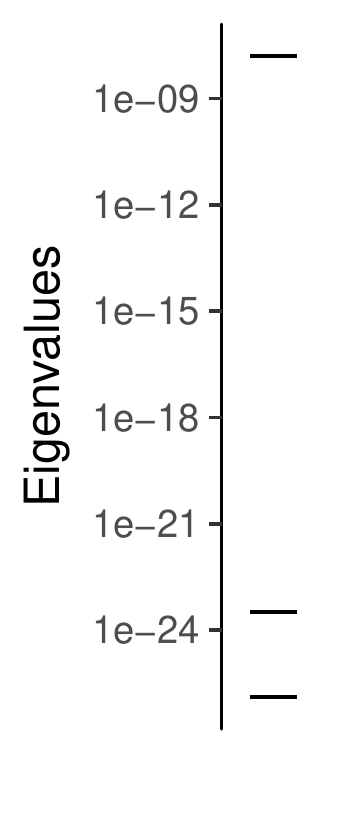}
%		\caption{}
	\end{subfigure}
	\begin{subfigure}[b]{0.22\textwidth}	
	\includegraphics[width = \textwidth]{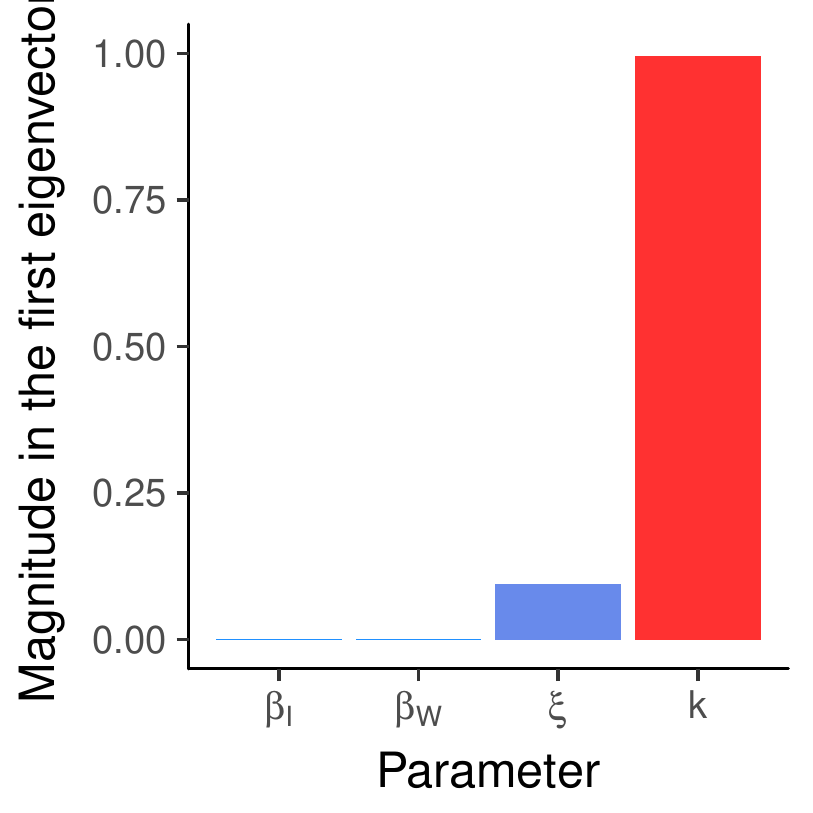}
	\end{subfigure}
	\begin{subfigure}[b]{0.22\textwidth}
	\includegraphics[width = \textwidth]{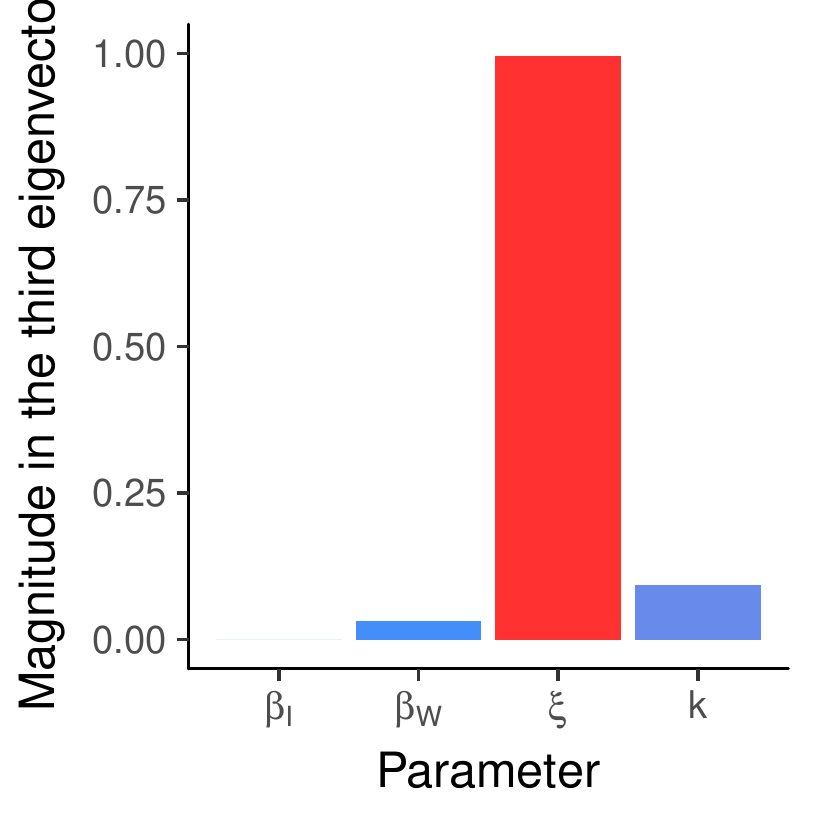}
%		\caption{}
	\end{subfigure}
	\begin{subfigure}[b]{0.22\textwidth}	
	\includegraphics[width = \textwidth]{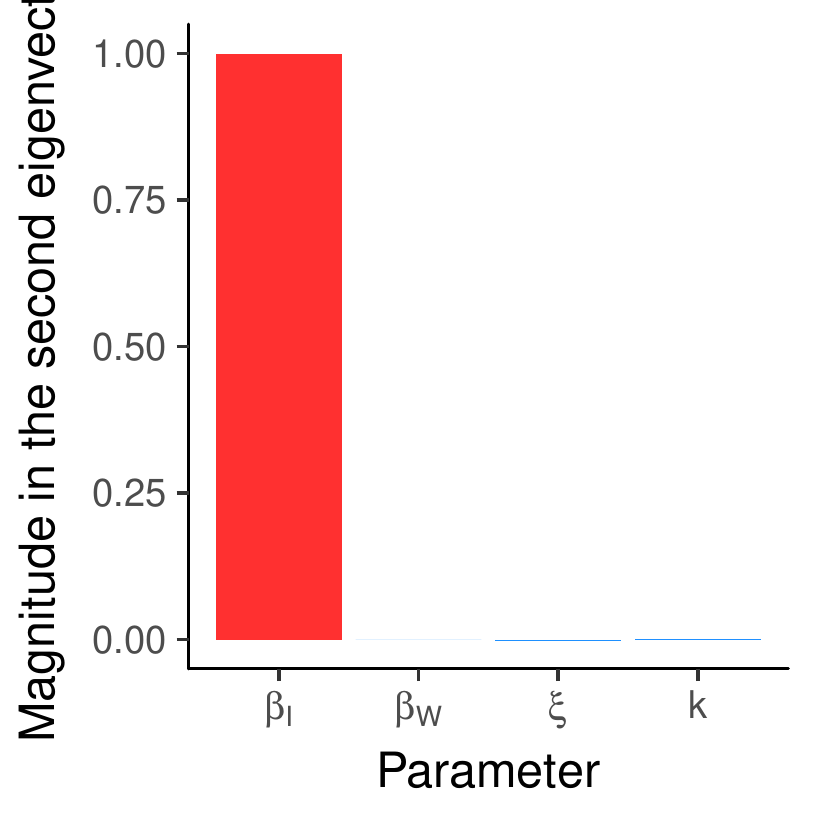}
	\end{subfigure}
	\begin{subfigure}[b]{0.22\textwidth}	
	\includegraphics[width = \textwidth]{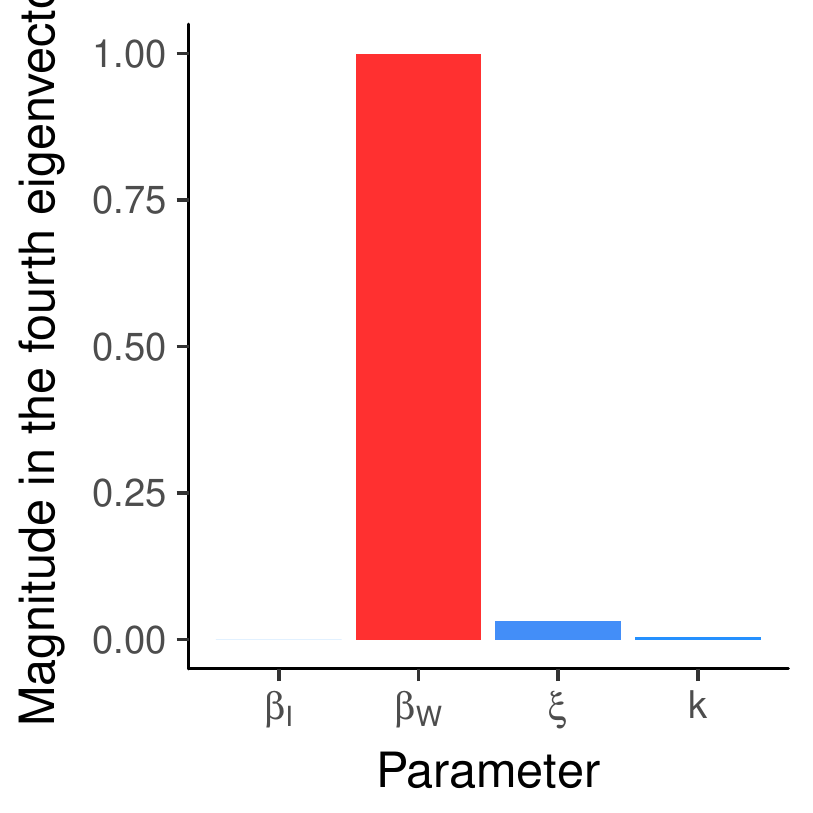}
%		\caption{}
	\end{subfigure}
	\\
$C_y$\\
	\begin{subfigure}[b]{0.09\textwidth}	
	\includegraphics[height=1.5in]{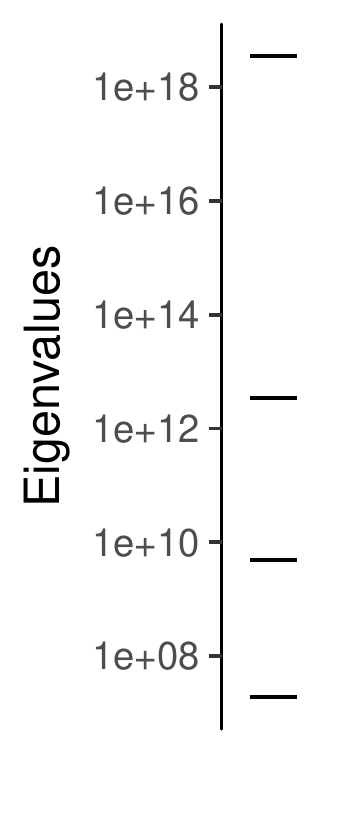}
%		\caption{}
	\end{subfigure}
	\begin{subfigure}[b]{0.22\textwidth}	
	\includegraphics[width = \textwidth]{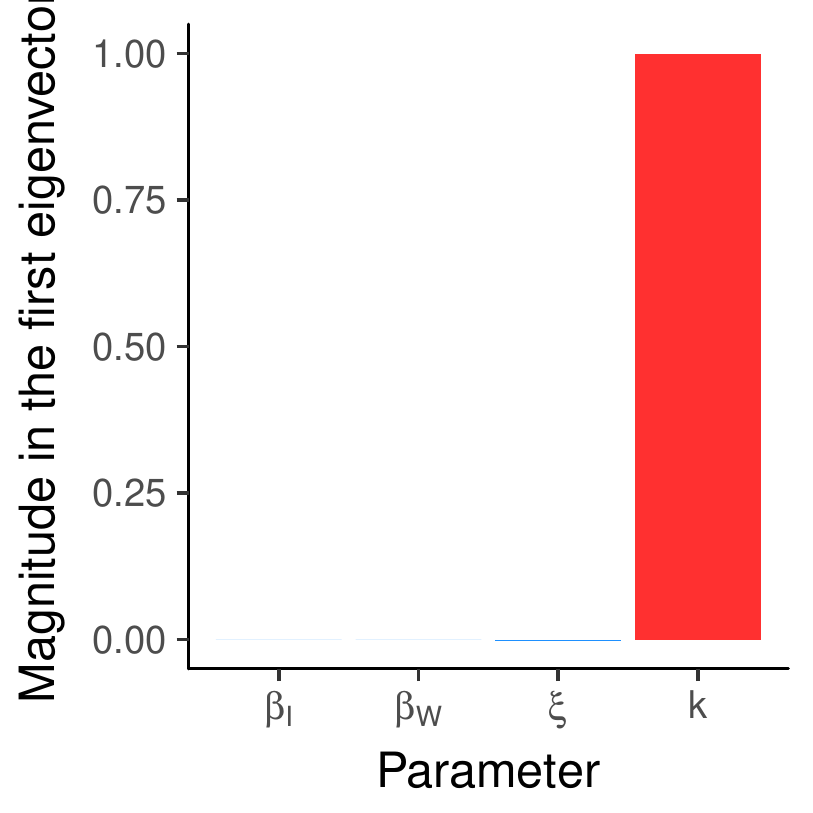}
	\end{subfigure}
	\begin{subfigure}[b]{0.22\textwidth}
	\includegraphics[width = \textwidth]{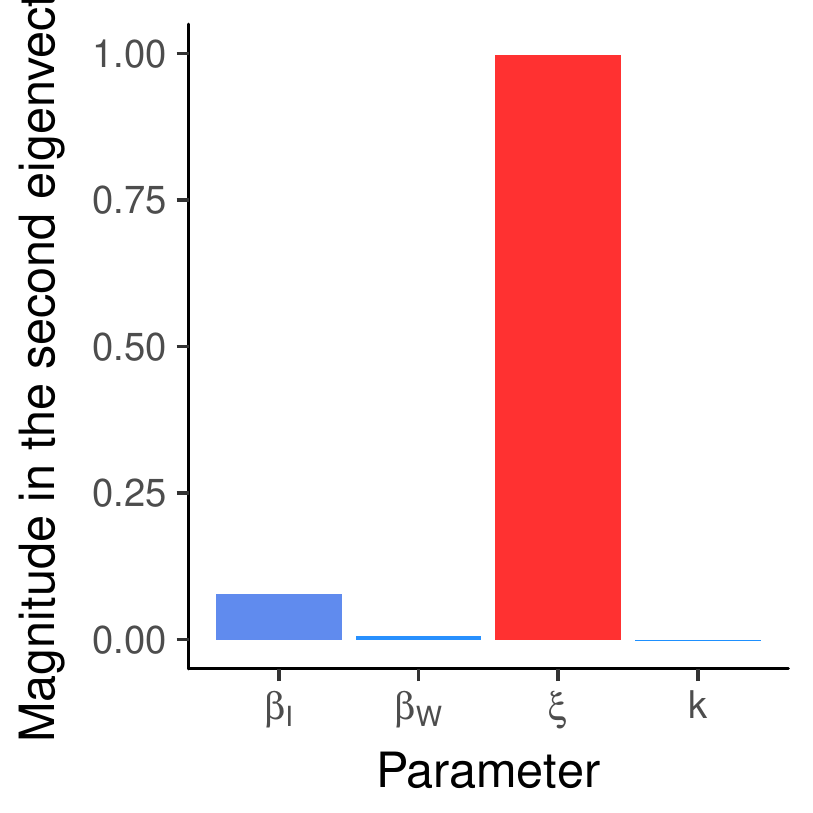}
%		\caption{}
	\end{subfigure}
	\begin{subfigure}[b]{0.22\textwidth}	
	\includegraphics[width = \textwidth]{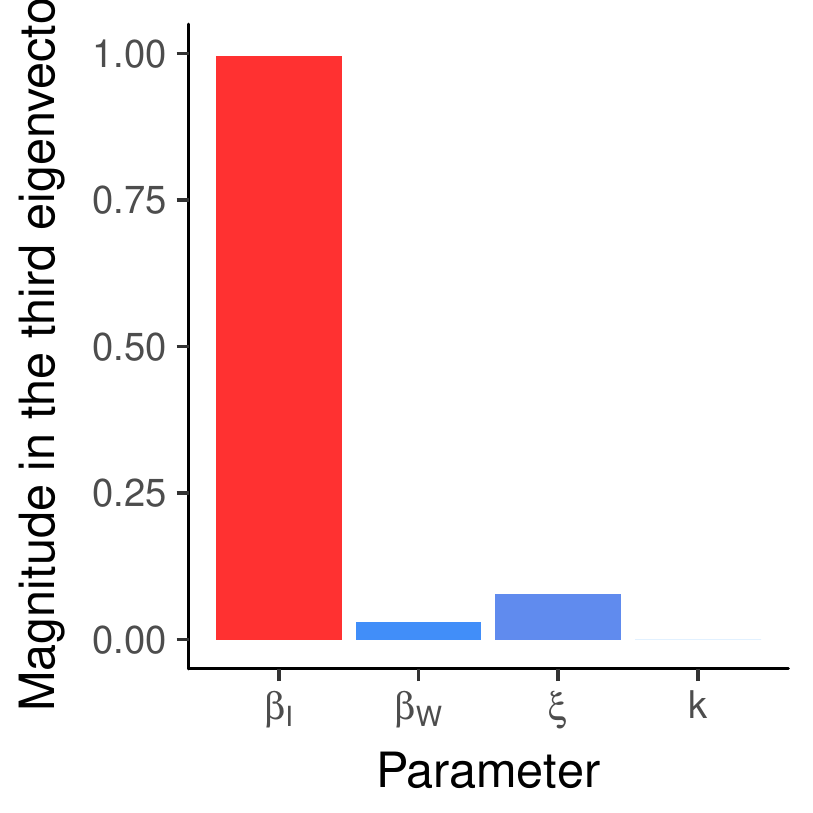}
	\end{subfigure}
	\begin{subfigure}[b]{0.22\textwidth}	
	\includegraphics[width = \textwidth]{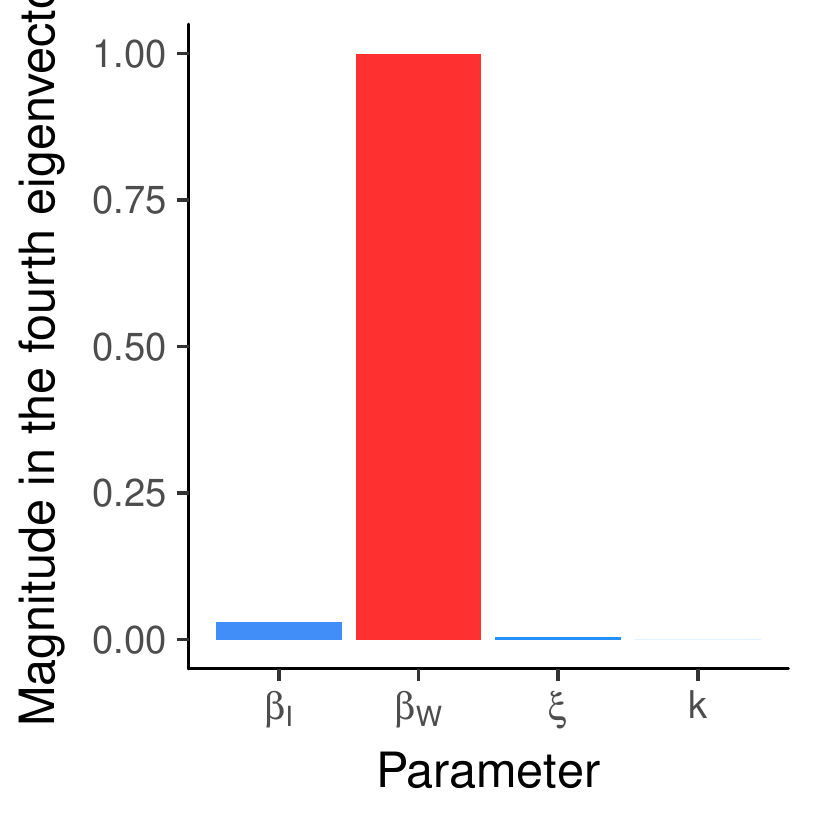}
%		\caption{}
	\end{subfigure}
\\
$C_q$\\
	\begin{subfigure}[b]{0.09\textwidth}	
	\includegraphics[height=1.5in]{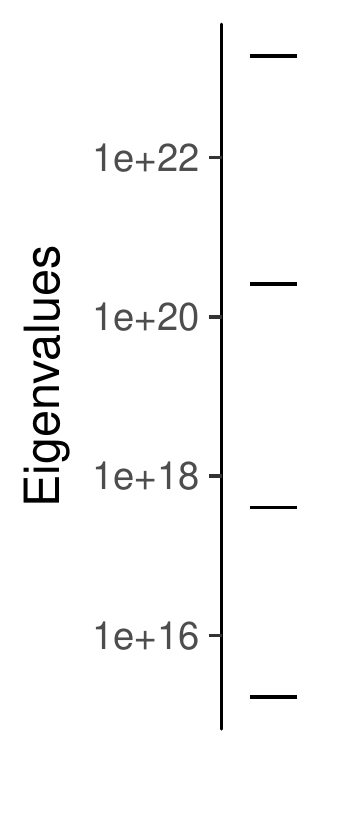}
%		\caption{}
	\end{subfigure}
	\begin{subfigure}[b]{0.22\textwidth}	
	\includegraphics[width = \textwidth]{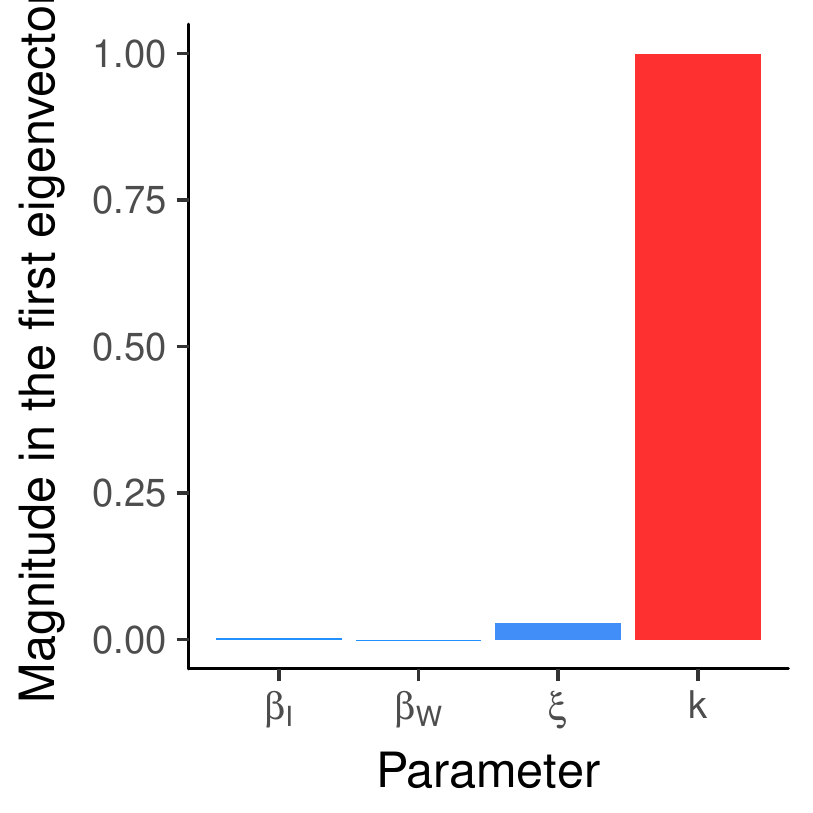}
	\end{subfigure}
	\begin{subfigure}[b]{0.22\textwidth}
	\includegraphics[width = \textwidth]{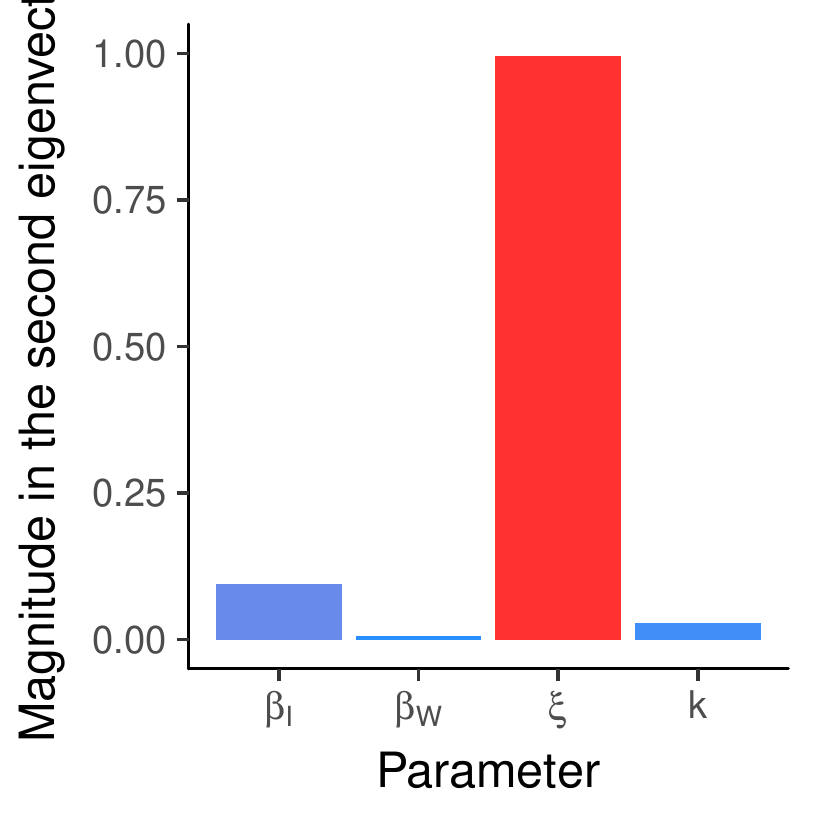}
%		\caption{}
	\end{subfigure}
	\begin{subfigure}[b]{0.22\textwidth}	
	\includegraphics[width = \textwidth]{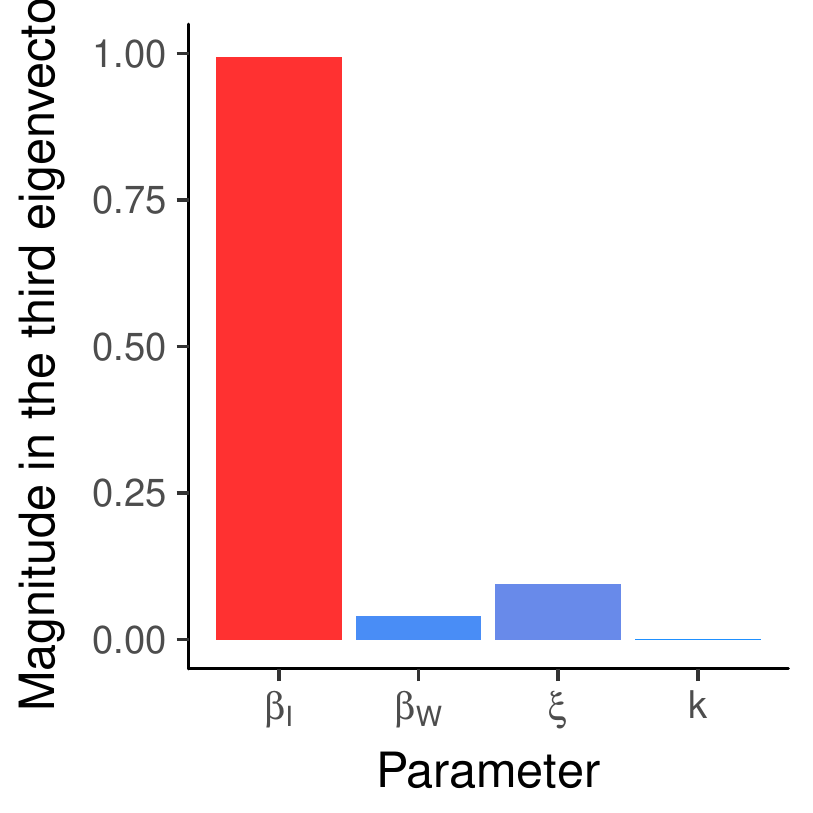}
	\end{subfigure}
	\begin{subfigure}[b]{0.22\textwidth}	
	\includegraphics[width = \textwidth]{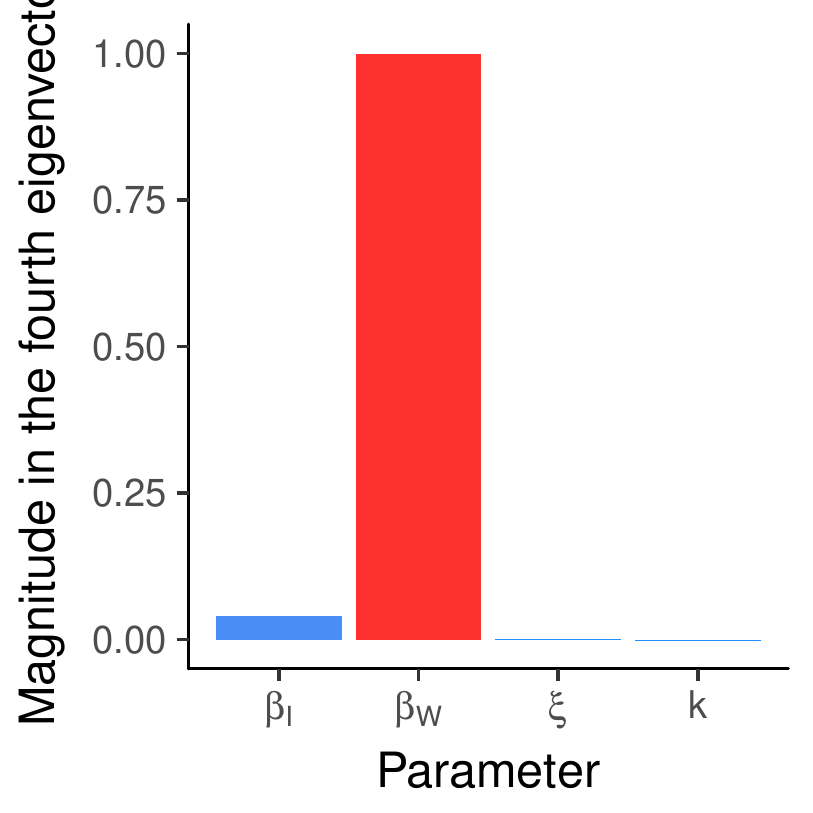}
%		\caption{}
	\end{subfigure}
\caption{Eigenvalues (leftmost column) and eigenvector magnitudes for the four sFIM-based quantities. For all four versions, the four eigenvalues were fairly evenly spaced apart. Note that for $F_q$, the second and third eigenvalues were very close and so appear as one bar in the eigenvalue plot (second row). Eigenvectors for all four sFIMs each captured one main parameter, with largely the same order, except for $F_q$, where the third eigenvector was primarily $\xi$ and the second was $\beta_I$. Because the two associated eigenvalues were so close, for clarity we swapped the order of the two eigenvectors in the plot.}
\label{fig:normalxi}
\end{figure}

We first take $\hat{\theta}$ to be the default parameters described above. From these default parameters and our LH sample, we calculated four sFIM-related quantities, two local and two global:
\begin{itemize}
\item $F_y$, the sFIM  at $\hat{\theta}$ calculated using our QOI $y$ (time series vector)
\item $F_q$, the sFIM at $\hat{\theta}$ calculated using the summary QOI $q$ (sum of squares using $y(\hat{\theta})$ as the `data') 
\item $C_y$, the average sFIM calculated using our QOI $y$ (time series vector)
\item $C_q$, the average sFIM calculated using the summary QOI $q$ (sum of squares using $y(\hat{\theta})$ as the `data')
\end{itemize}
Figure~\ref{fig:normalxi} shows the eigenvalues and eigenvector component magnitudes for each sFIM type. For all four quantities, the four eigenvalues were fairly evenly spaced apart. In terms of eigenvectors, the four sFIM versions agreed quite closely, with each eigenvector corresponding largely to a single parameter. The sufficient summary plots do not show any cohesive patterns (not shown), consistent with the profile likelihood results that the model is identifiable in this region of parameter space. 

\begin{figure}%[h!]
\centering
%$C_y$\\
	\begin{subfigure}[b]{0.09\textwidth}	
	\includegraphics[height=1.5in]{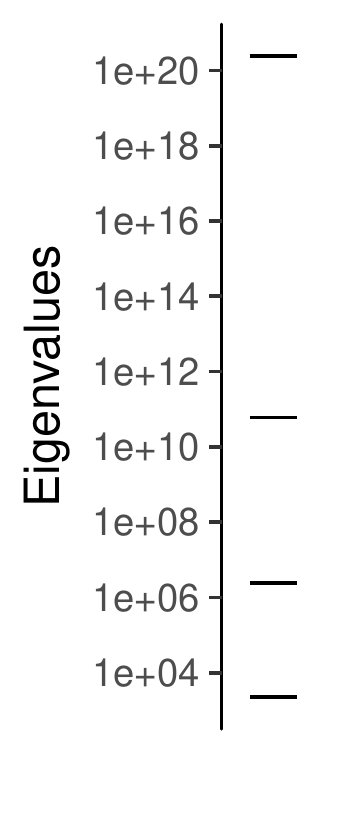}
%		\caption{}
	\end{subfigure}
	\begin{subfigure}[b]{0.22\textwidth}	
	\includegraphics[width = \textwidth]{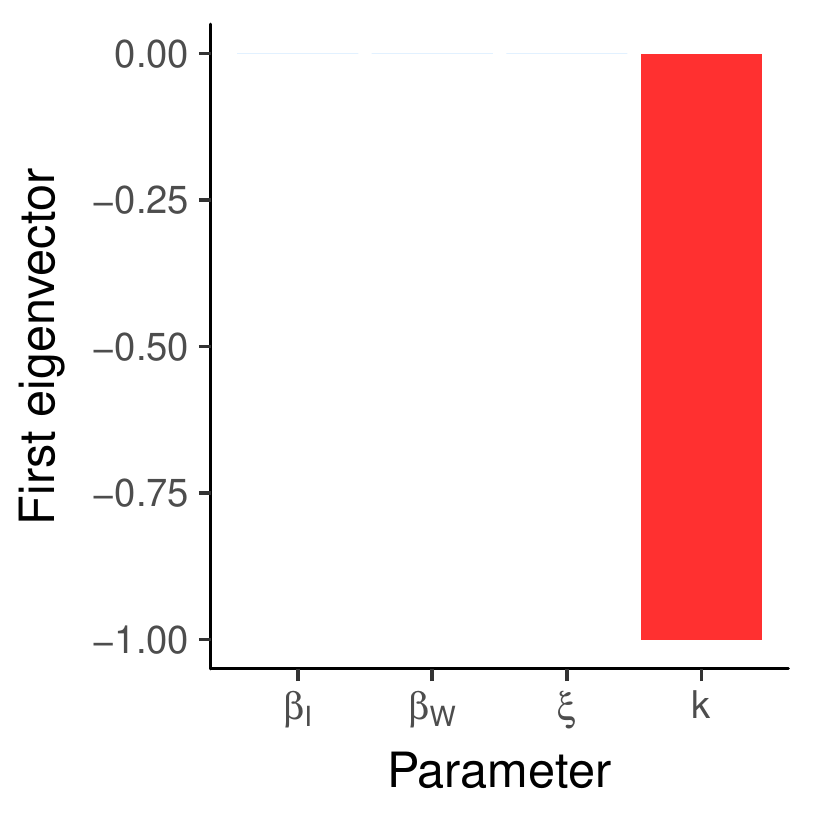}
	\end{subfigure}
	\begin{subfigure}[b]{0.22\textwidth}
	\includegraphics[width = \textwidth]{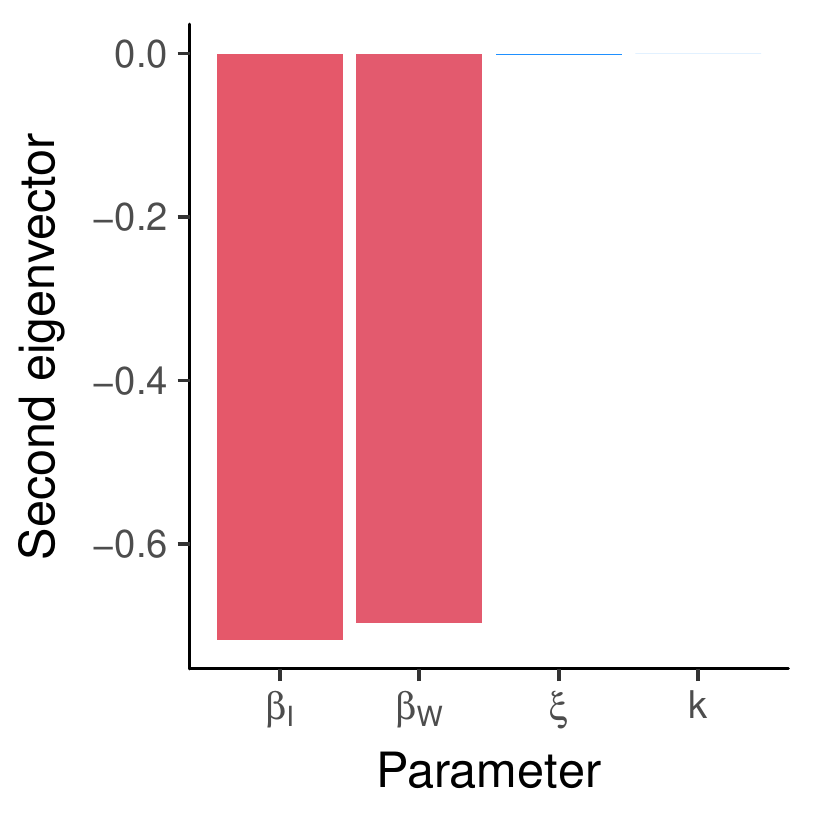}
%		\caption{}
	\end{subfigure}
	\begin{subfigure}[b]{0.22\textwidth}	
	\includegraphics[width = \textwidth]{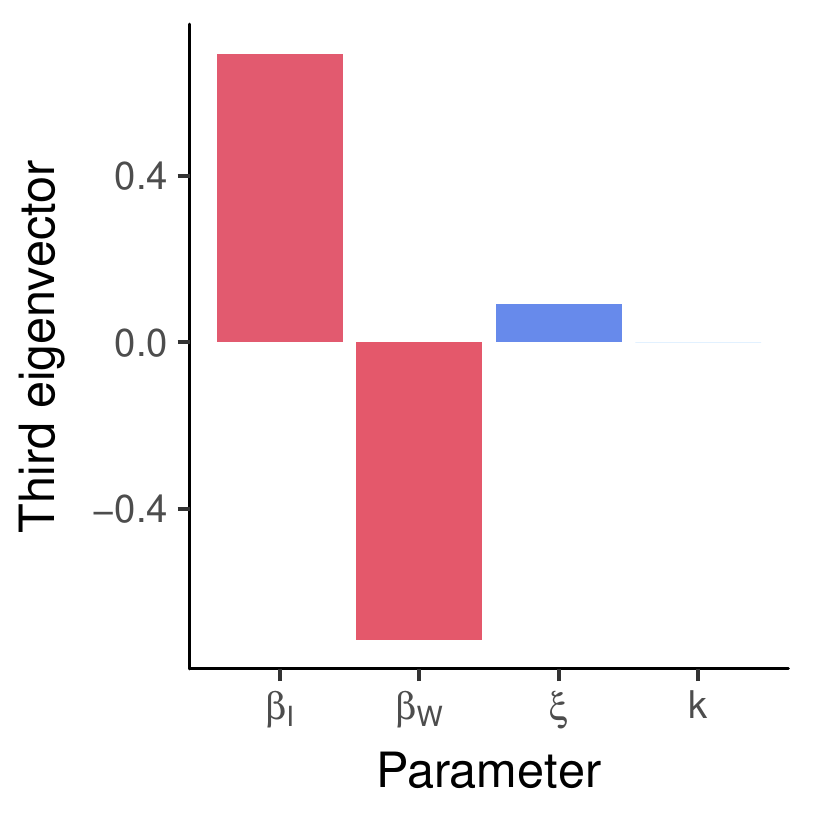}
	\end{subfigure}
	\begin{subfigure}[b]{0.22\textwidth}	
	\includegraphics[width = \textwidth]{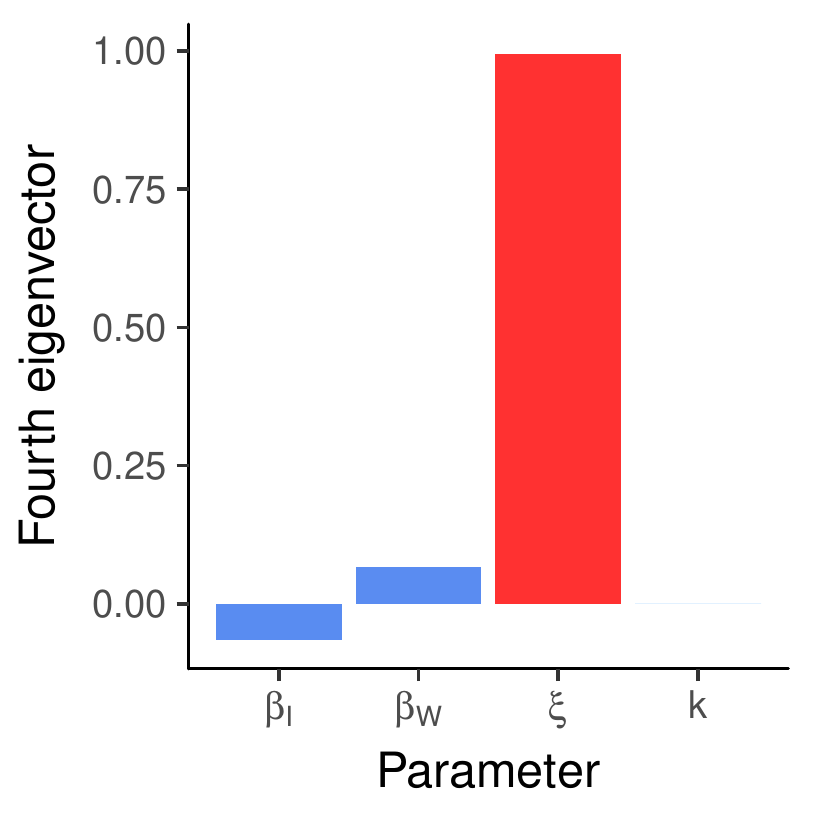}
%		\caption{}
	\end{subfigure}
\caption{Eigenvalues (leftmost panel) and eigenvectors for the four $C_y$ (the average sFIM, evaluated for the full vector QOI $y$), in the fast $\xi$ case. The first two eigenvectors capture the active directions corresponding to identifiable combinations $k$ and $\beta_W+\beta_I$. The third and fourth eigenvectors correspond to the inactive or unidentifiable directions, corresponding to compensation between $\beta_W$ and $\beta_I$, and $\xi$.}
\label{fig:fastxi}
\end{figure}

\begin{figure}[h!]
\centering
$F_y$\\
	\begin{subfigure}[b]{0.09\textwidth}	
	\includegraphics[height=1.5in]{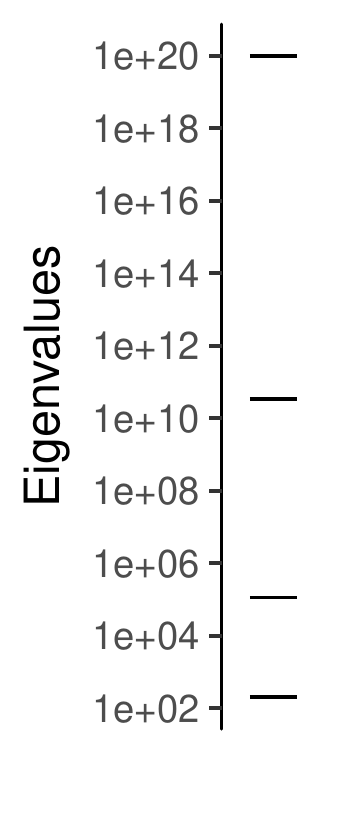}
%		\caption{}
	\end{subfigure}
	\begin{subfigure}[b]{0.22\textwidth}	
	\includegraphics[width = \textwidth]{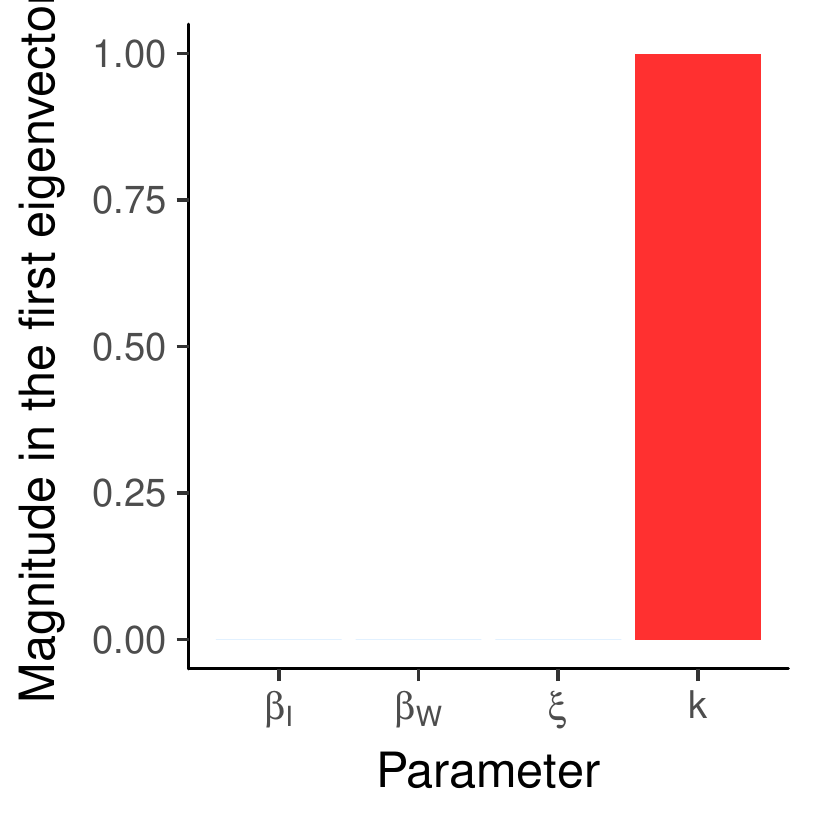}
	\end{subfigure}
	\begin{subfigure}[b]{0.22\textwidth}
	\includegraphics[width = \textwidth]{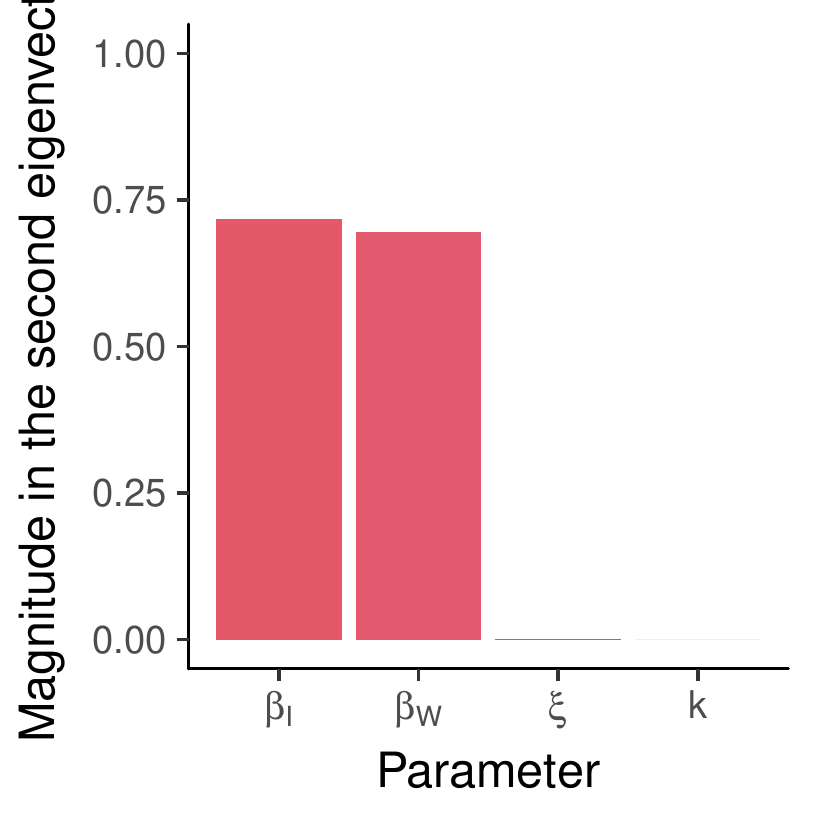}
%		\caption{}
	\end{subfigure}
	\begin{subfigure}[b]{0.22\textwidth}	
	\includegraphics[width = \textwidth]{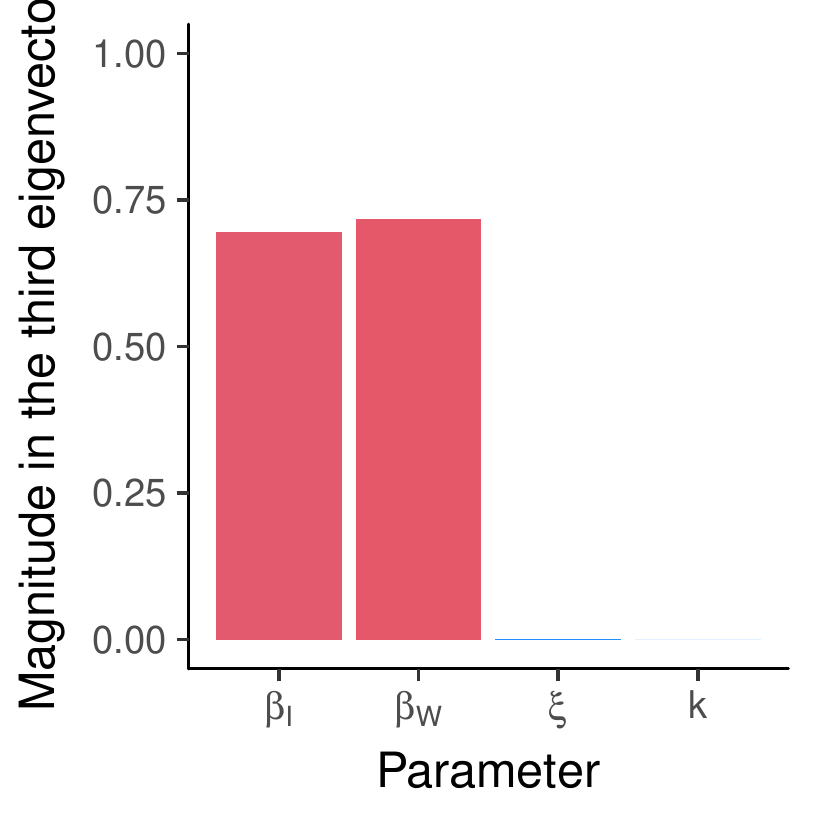}
	\end{subfigure}
	\begin{subfigure}[b]{0.22\textwidth}	
	\includegraphics[width = \textwidth]{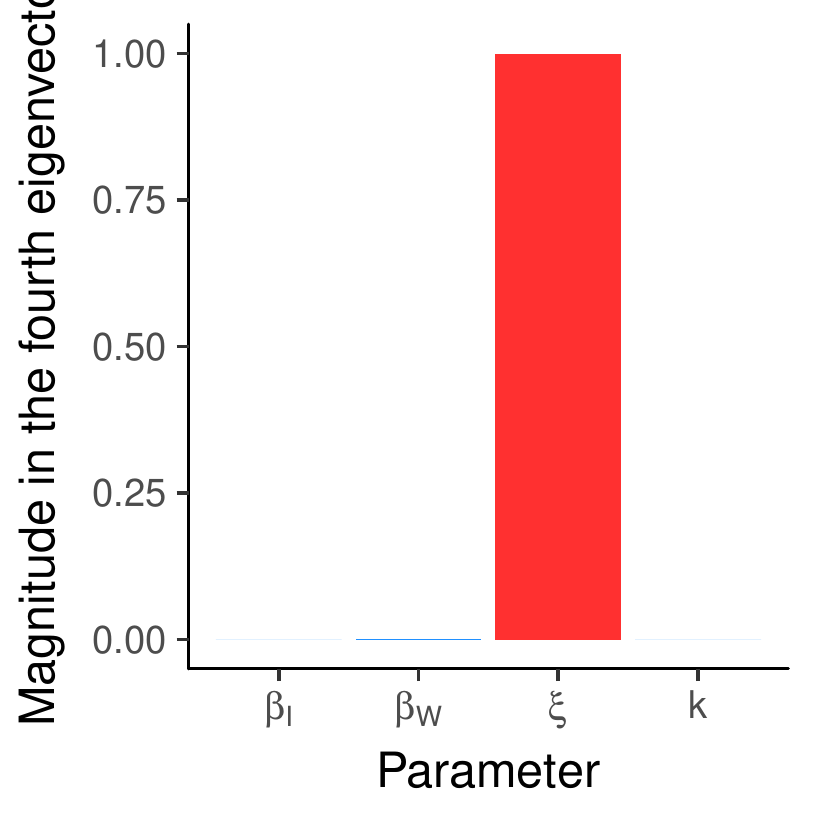}
%		\caption{}
	\end{subfigure}
\\
$F_q$\\
	\begin{subfigure}[b]{0.09\textwidth}	
	\includegraphics[height=1.5in]{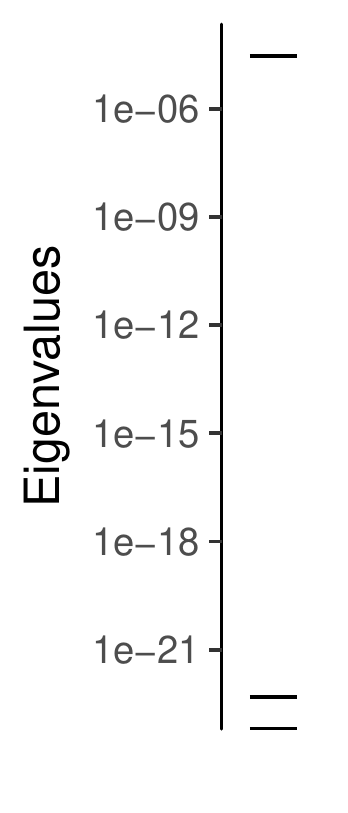}
%		\caption{}
	\end{subfigure}
	\begin{subfigure}[b]{0.22\textwidth}	
	\includegraphics[width = \textwidth]{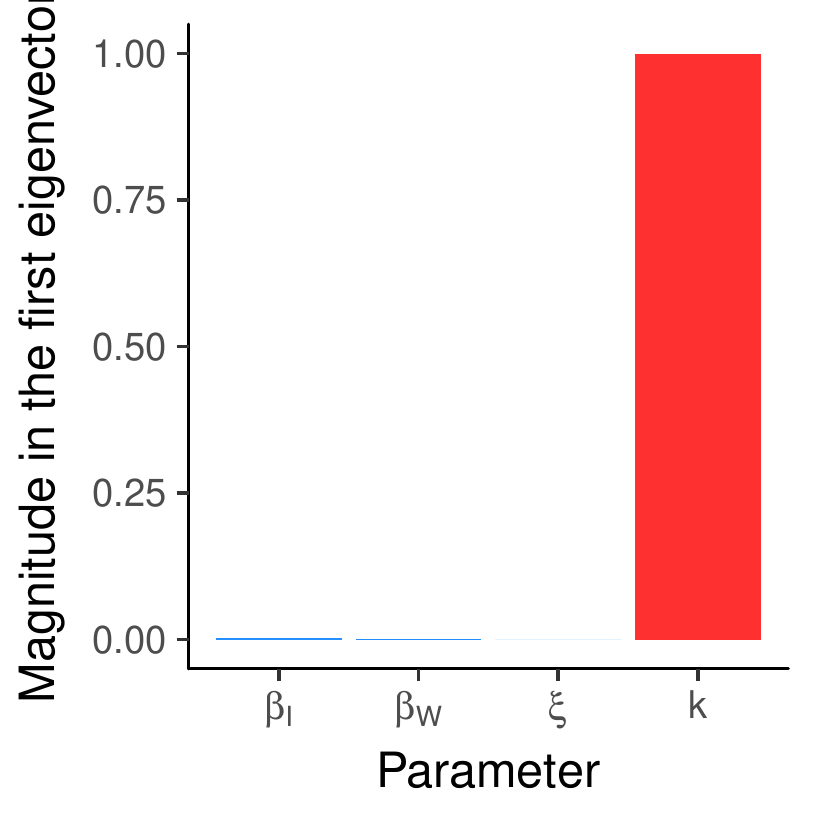}
	\end{subfigure}
	\begin{subfigure}[b]{0.22\textwidth}
	\includegraphics[width = \textwidth]{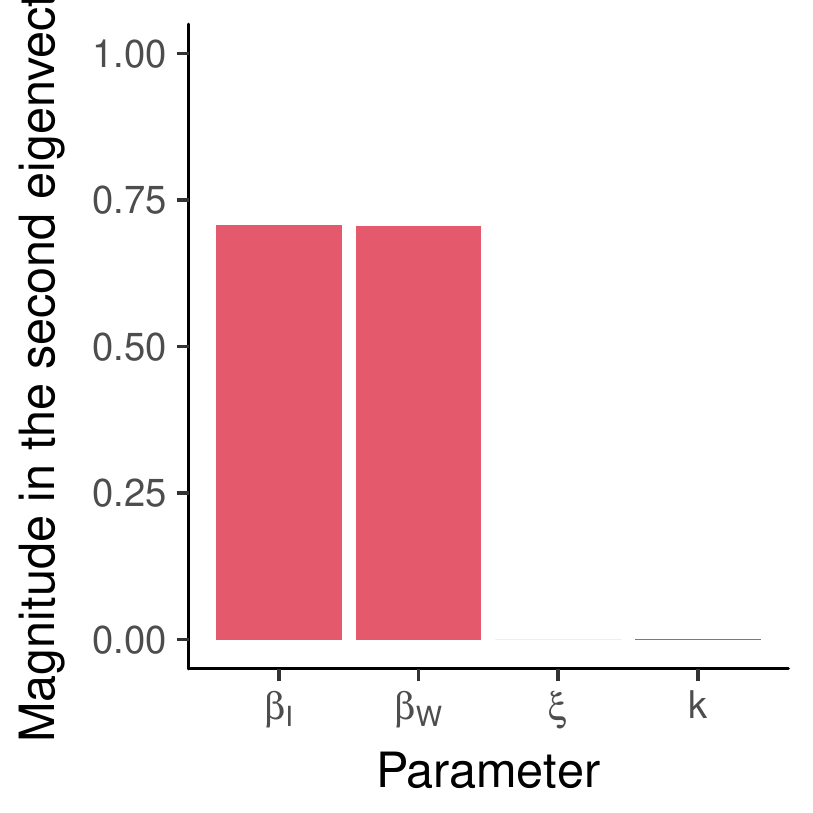}
%		\caption{}
	\end{subfigure}
	\begin{subfigure}[b]{0.22\textwidth}	
	\includegraphics[width = \textwidth]{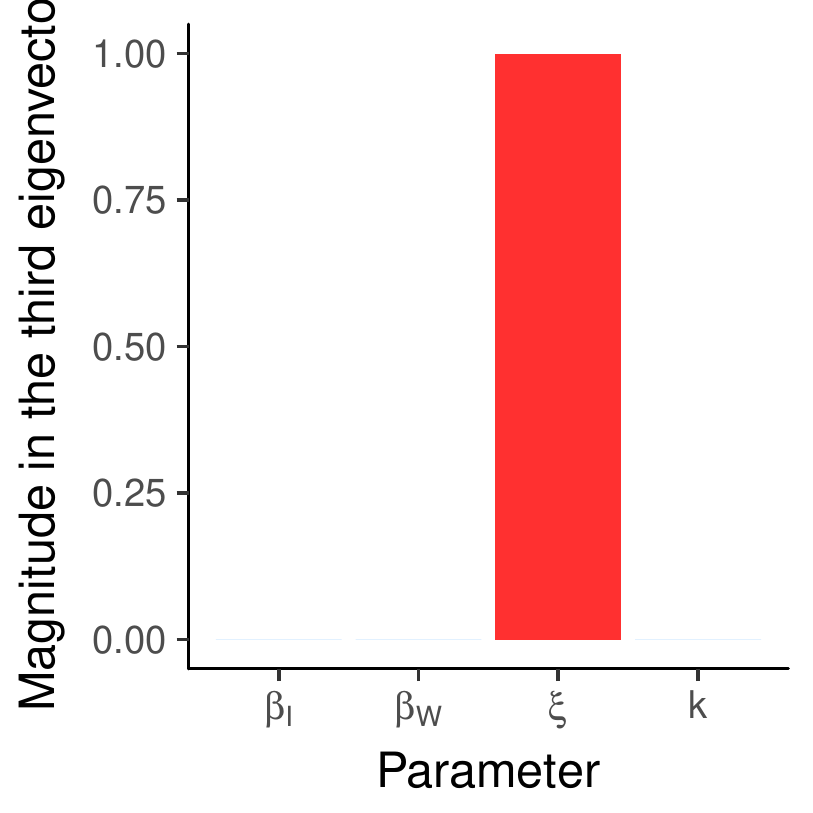}
	\end{subfigure}
	\begin{subfigure}[b]{0.22\textwidth}	
	\includegraphics[width = \textwidth]{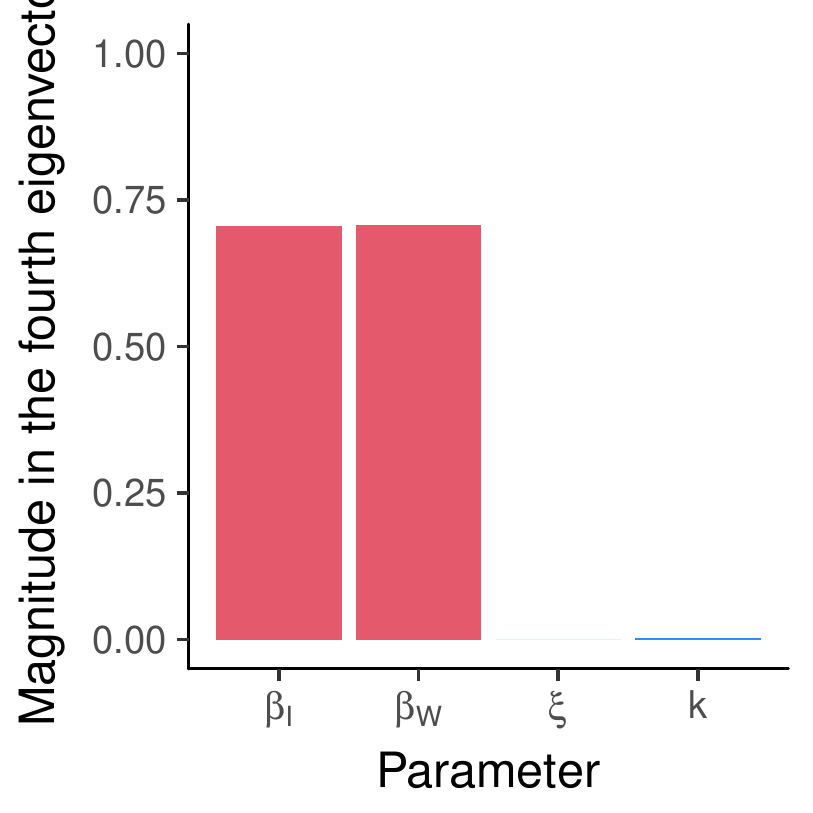}
%		\caption{}
	\end{subfigure}
	\\
$C_y$\\
	\begin{subfigure}[b]{0.09\textwidth}	
	\includegraphics[height=1.5in]{fastxi-non-scaled-50percent/eigplot-Cvec}
%		\caption{}
	\end{subfigure}
	\begin{subfigure}[b]{0.22\textwidth}	
	\includegraphics[width = \textwidth]{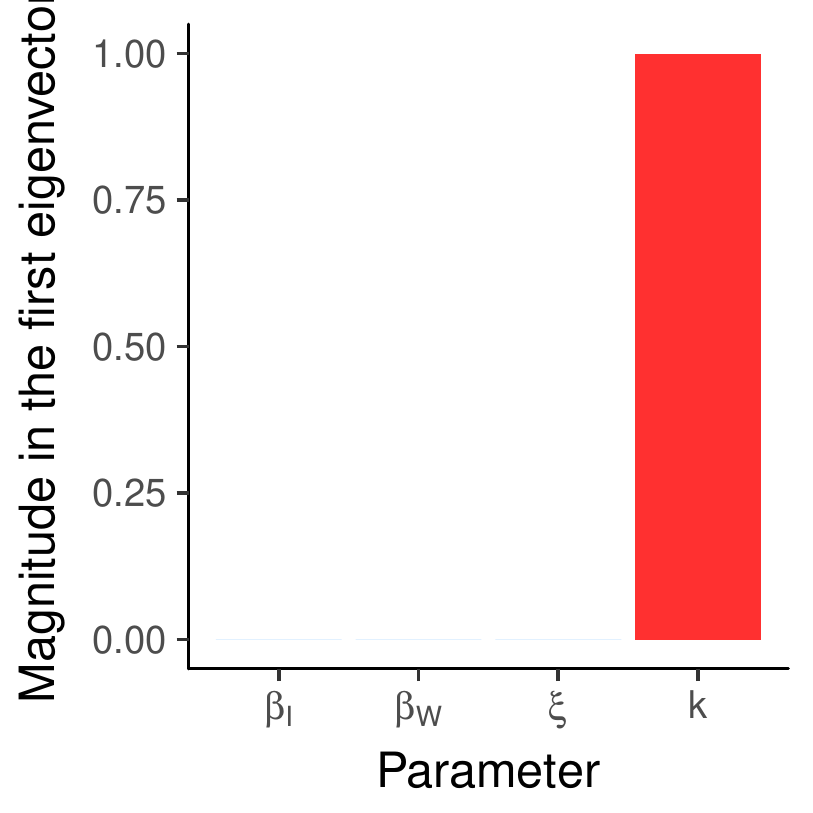}
	\end{subfigure}
	\begin{subfigure}[b]{0.22\textwidth}
	\includegraphics[width = \textwidth]{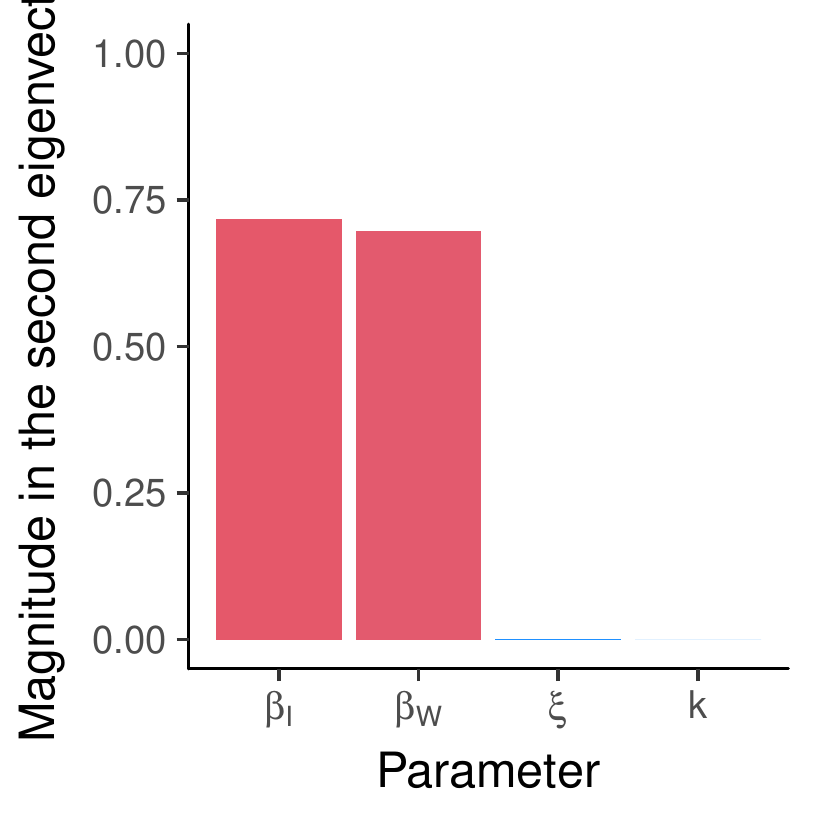}
%		\caption{}
	\end{subfigure}
	\begin{subfigure}[b]{0.22\textwidth}	
	\includegraphics[width = \textwidth]{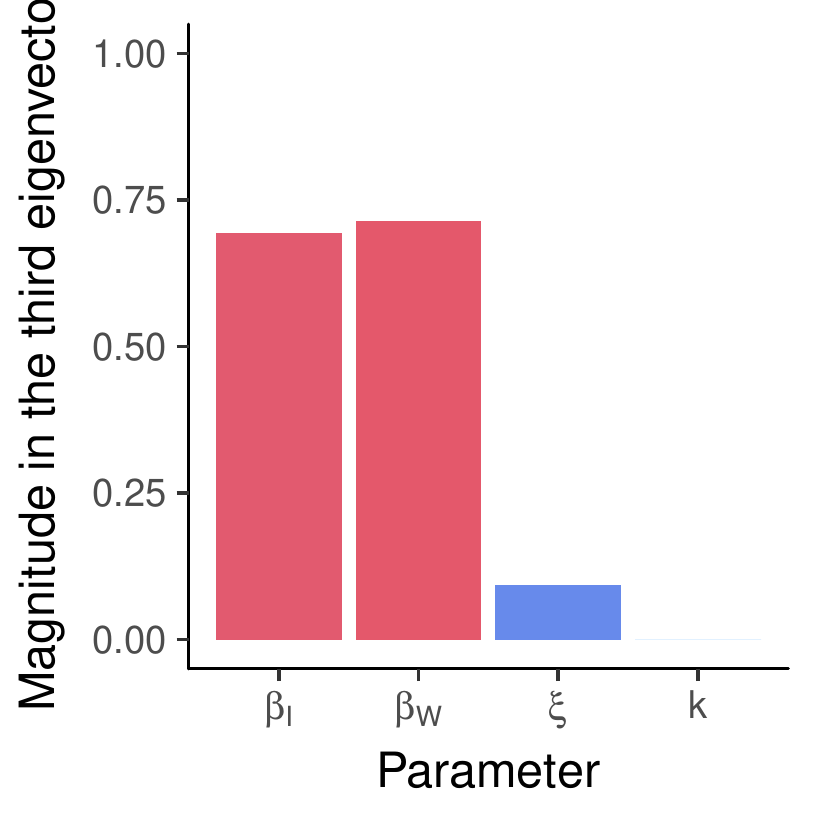}
	\end{subfigure}
	\begin{subfigure}[b]{0.22\textwidth}	
	\includegraphics[width = \textwidth]{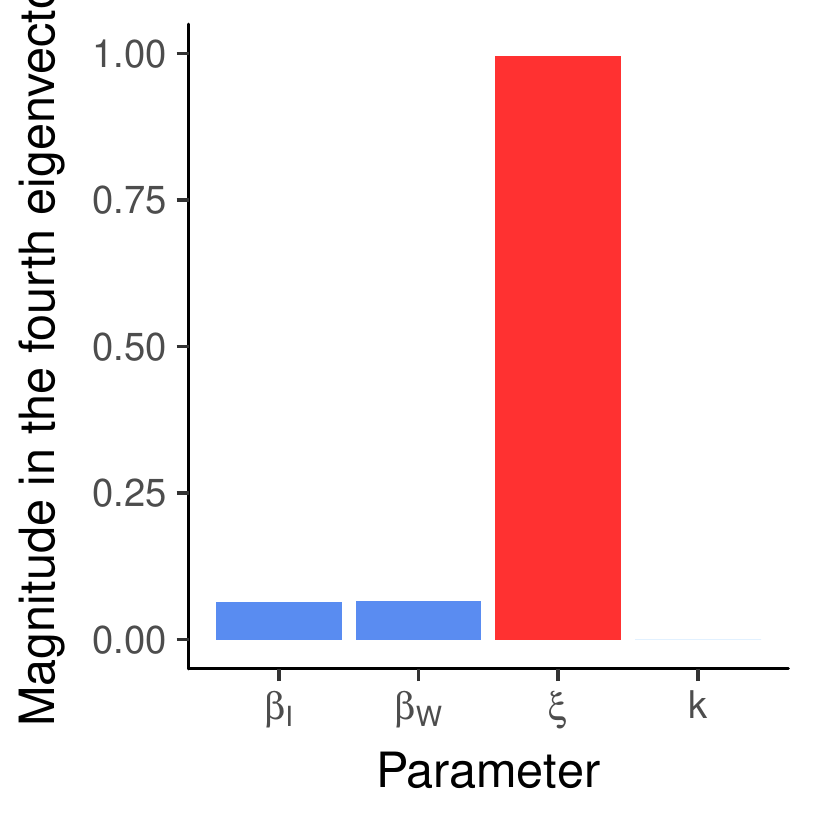}
%		\caption{}
	\end{subfigure}
\\
$C_q$\\
	\begin{subfigure}[b]{0.09\textwidth}	
	\includegraphics[height=1.5in]{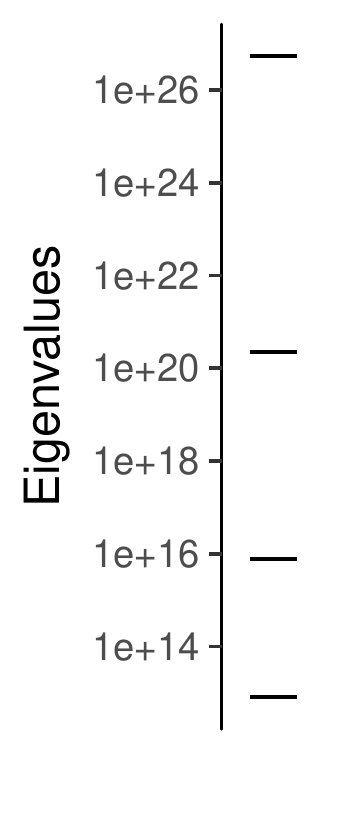}
%		\caption{}
	\end{subfigure}
	\begin{subfigure}[b]{0.22\textwidth}	
	\includegraphics[width = \textwidth]{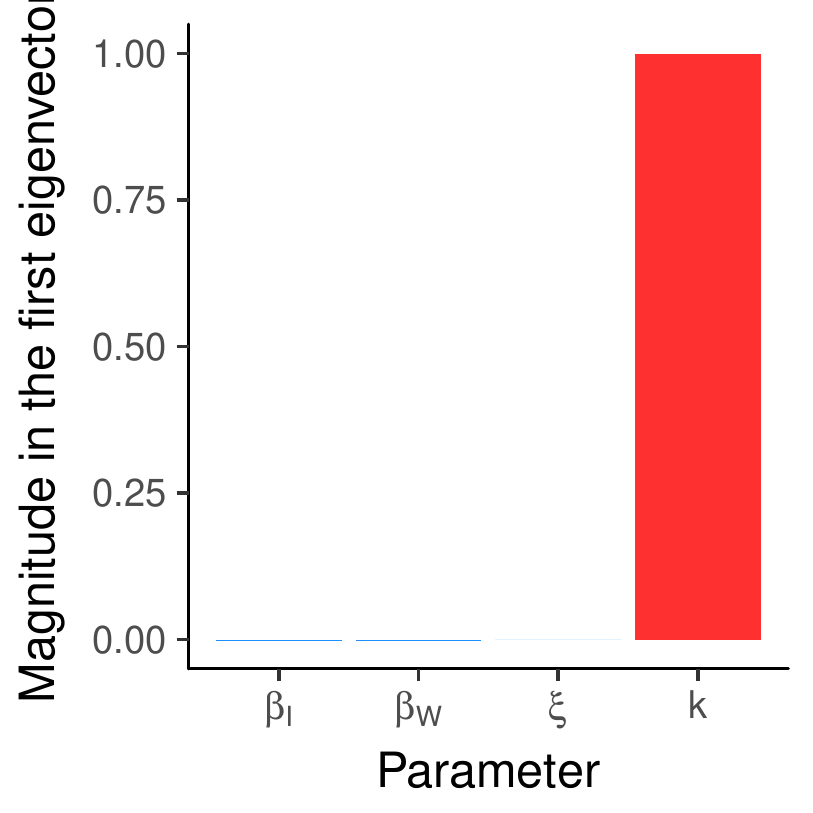}
	\end{subfigure}
	\begin{subfigure}[b]{0.22\textwidth}
	\includegraphics[width = \textwidth]{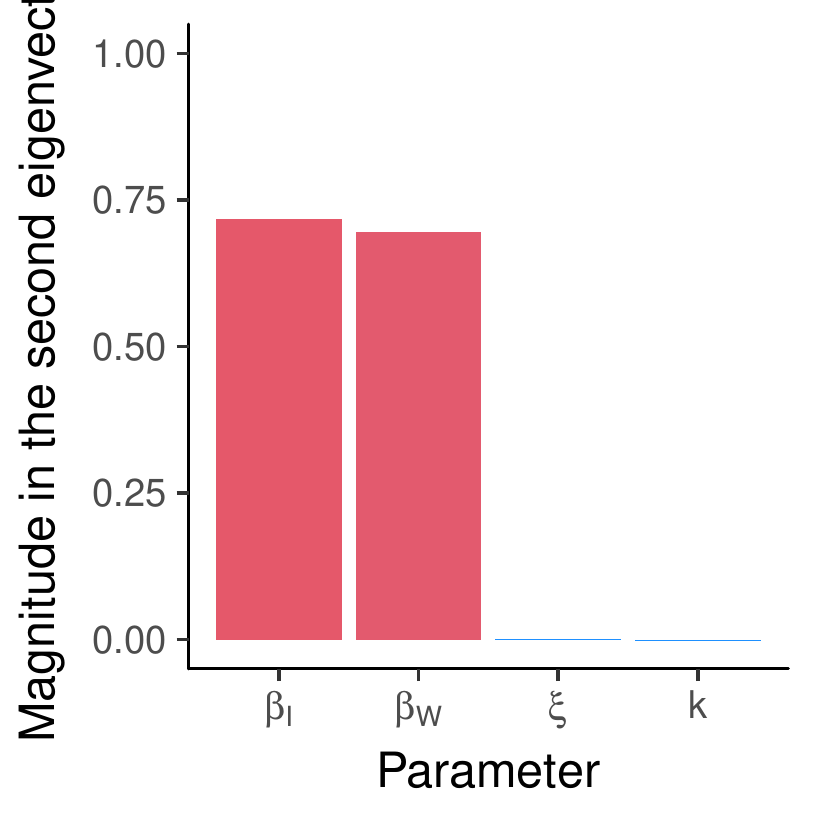}
%		\caption{}
	\end{subfigure}
	\begin{subfigure}[b]{0.22\textwidth}	
	\includegraphics[width = \textwidth]{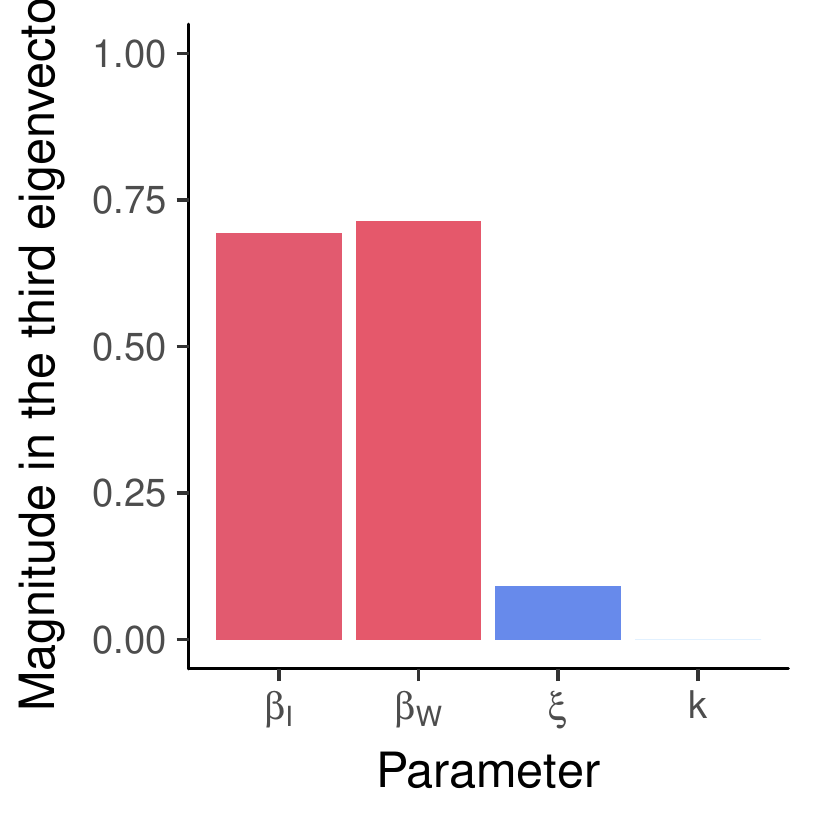}
	\end{subfigure}
	\begin{subfigure}[b]{0.22\textwidth}	
	\includegraphics[width = \textwidth]{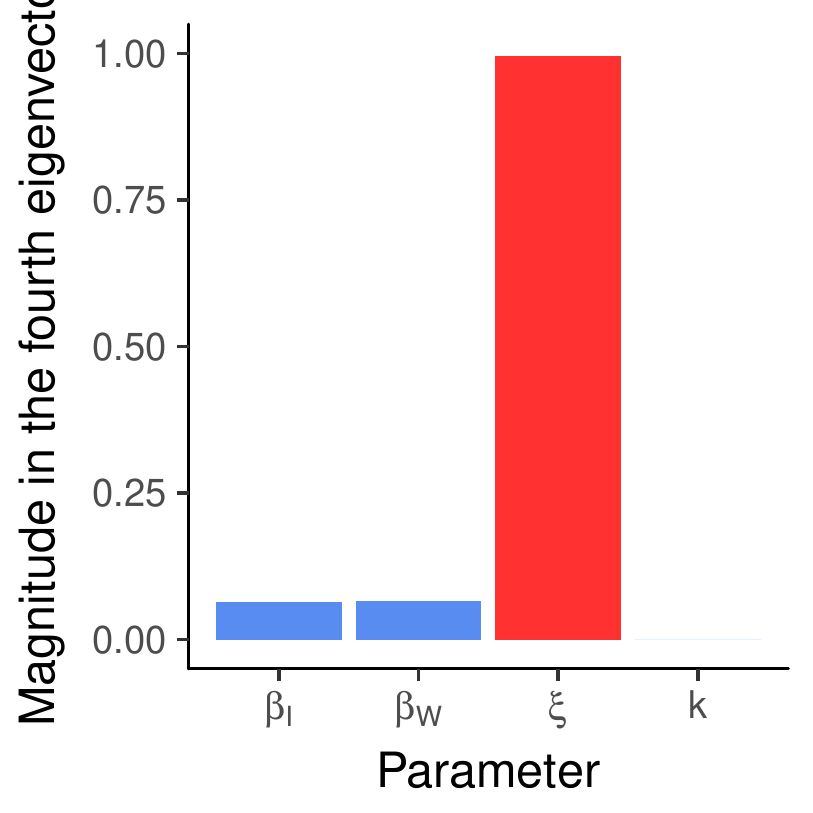}
%		\caption{}
	\end{subfigure}
\caption{Eigenvalues (leftmost column) and eigenvector loads for the four sFIM-based quantities, in the fast $\xi$ case. The first two eigenvectors capture the active directions corresponding to identifiable combinations $k$ and $\beta_W+\beta_I$. The lower two eigenvectors correspond to the inactive or unidentifiable directions, corresponding to compensation between $\beta_W$ and $\beta_I$, and $\xi$.}
\label{fig:fastxicompare}
\end{figure}

Next, we ran the same analyses, but with the larger value of $\xi = 37.8$ used in Figure~\ref{SIWR_profiles}c, shown in Figures~\ref{fig:fastxi} and~\ref{fig:fastxicompare}. The eigenvector directions are similarly consistent across all four sFIM quantities (Appendix Figure~\ref{fig:fastxicompare}). The first two eigenvectors capture the sensitive directions for the model, approximately matching $k$ and $\beta_W+\beta_I$, the main identifiable parameters. The lower two eigenvalues capture the inactive/unidentifiable directions representing compensation between the two transmission routes and $\xi$. 

For both the fast and slow $\xi$ regimes, the similarity of eigenvectors with all four quantities (local/global, vector/scalar QOIs) highlights how active subspaces can be used to explore different regions of parameter space, and also how for more regional analyses, a cost-function based scalar QOI can be useful as a way of summarizing a vector or time series of QOI's.

\begin{figure}%[h!]
\centering
	\begin{subfigure}[b]{0.45\textwidth}	
	\includegraphics[width = \textwidth]{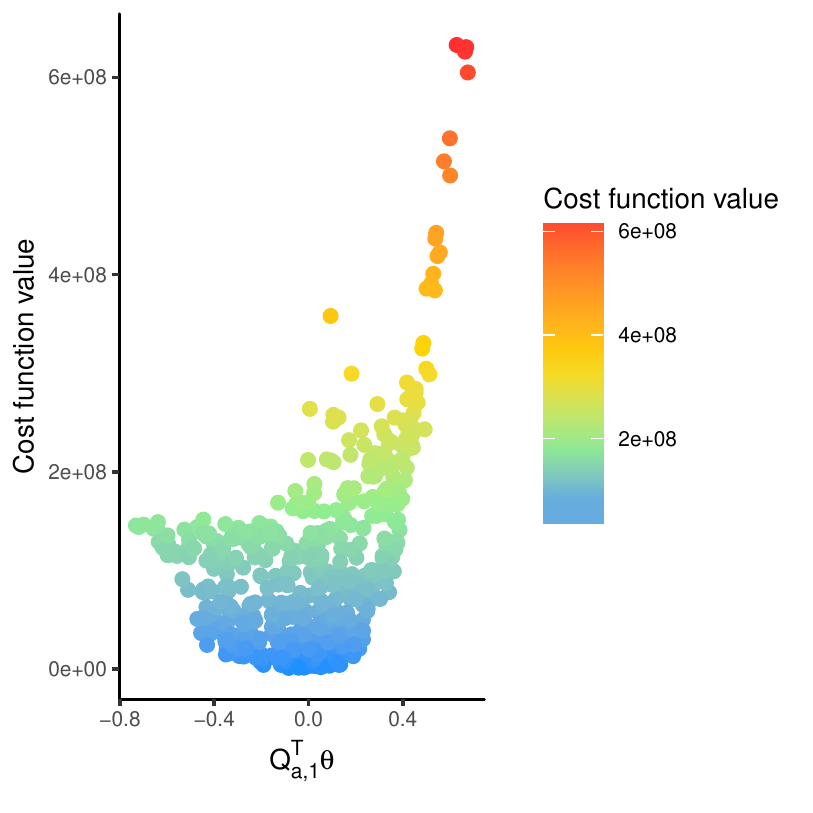}
		\caption{$C_q$ - normal $\xi$}
	\end{subfigure}
	\begin{subfigure}[b]{0.45\textwidth}
	\includegraphics[width = \textwidth]{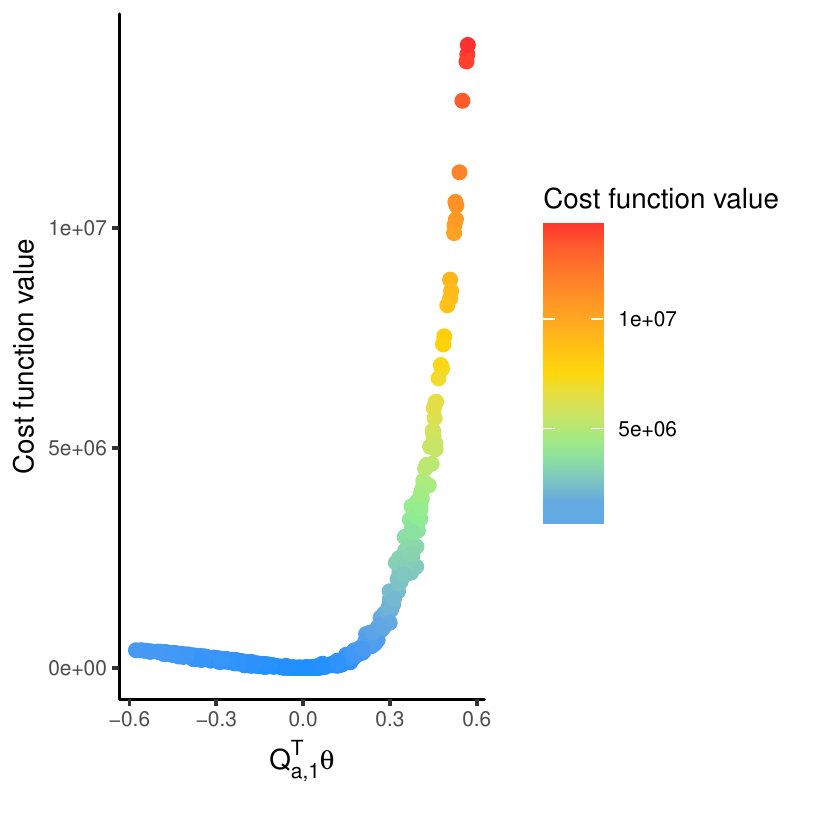}
		\caption{$C_q$ - fast $\xi$}
	\end{subfigure}\\
	\begin{subfigure}[b]{0.45\textwidth}	
	\includegraphics[width = \textwidth]{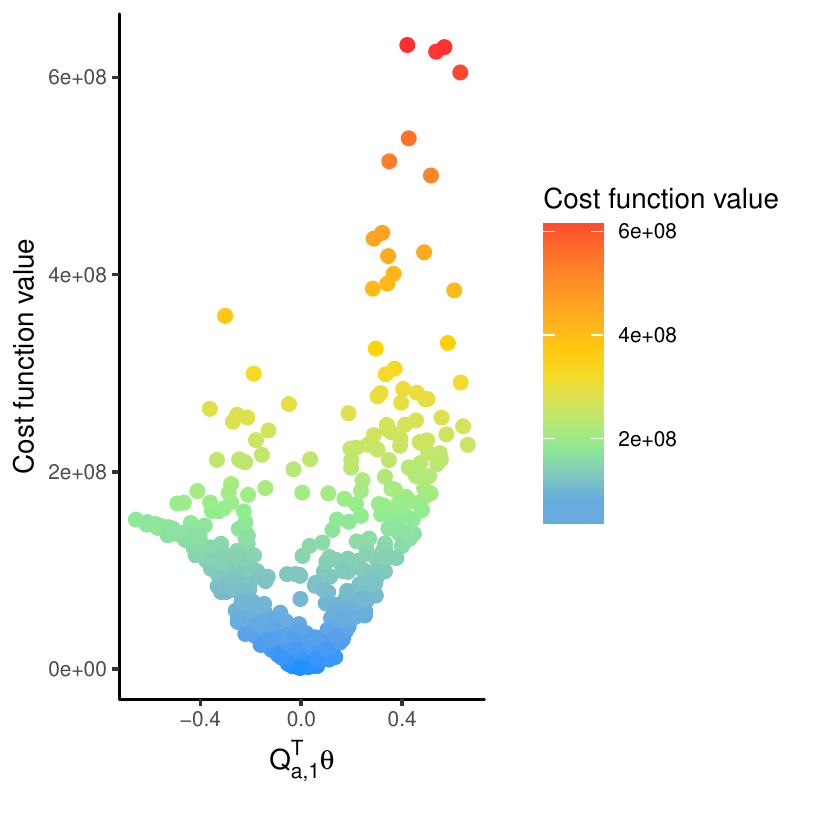}
		\caption{$C_y$ - normal $\xi$}
	\end{subfigure}
	\begin{subfigure}[b]{0.45\textwidth}	
	\includegraphics[width = \textwidth]{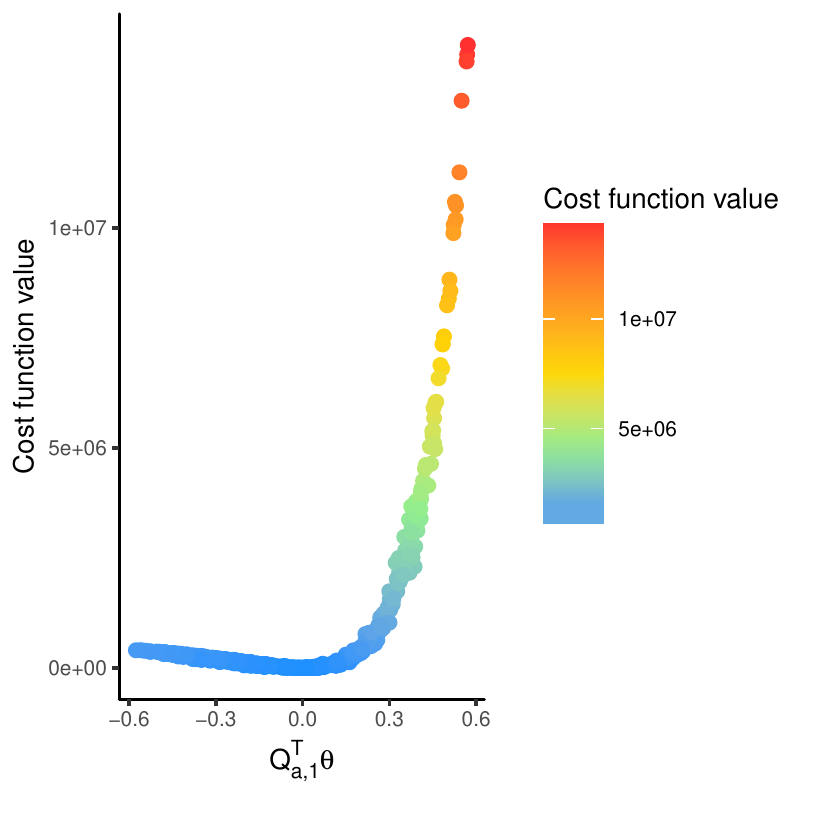}
		\caption{$C_y$ - fast $\xi$}
	\end{subfigure}
\caption{Sufficient summary plots for the first active eigenvector, with the average sFIM $C$ calculated using either the scalar QOI $q$ (top row) or the vector QOI $y$ (bottom row). In both cases, the cost function scalar QOI $q$ is used as the $y$-value in the plot. Left column panels show the normal $\xi$ case, which shows a somewhat unclear relationship between the first eigenvector and the cost function $q$, while the right column panels show the fast $\xi$ case, where the relationship with the cost function $q$ is close to one-dimensional.}
\label{fig:SIWRsum}
\end{figure}

Lastly, we ran the same four metrics, but with the parameters scaled and translated to be centered at zero (with a range of $(-1,1)$) in order to examine the sufficient summary plots \cite{constantine2015active}, shown in Figure~\ref{fig:SIWRsum}. The resulting plots show a much stronger univariate relationship in the unidentifiable, fast $\xi$ case than the default parameter case (consistent with the larger gap after the first eigenvalue in the fast $\xi$ case, shown in Appendix Figures~\ref{fig:SIWRscaledq} and~\ref{fig:SIWRscaledy}). Additionally, we generated sufficient summary plots by plotting $Q_a^T\theta$ for $C_y$ versus $q$ (our least squares cost function)---while $q$ was not the quantity used to generate $C_y$, the sufficient summary plots looked quite similar to those generated with $C_q$. Sufficient summary plots for the subsequent eigenvectors showed no clear trend in any of the cases (not shown). The associated eigenvalues and eigenvectors for the scaled versions of the parameters are given in Figures~\ref{fig:SIWRscaledq} and~\ref{fig:SIWRscaledy}.

\clearpage
\section*{Conclusions and future directions}
To conclude, in this paper we have examined the relationships between identifiability, active subspaces, and sloppiness using the sensitivity FIM as a common framework across each approach. By framing the active subspaces quantity $C$ as the average sensitivity FIM over the parameter space of interest, we were able to examine how local/global and linear/nonlinear identifiability and parameter reduction tools can each generate useful insights into a range of real-world applications. The framing of these parameter reduction tools in a parameter estimation context also let us examine the potential use of cost functions (such as least squares and likelihood functions), as a summary QOI when dealing with vector-QOIs (such as for time series). We hope that the sensitivity FIM-based framework developed here will facilitate further cross-talk between different areas of identifiability, uncertainty quantification, and parameter space reduction.

\clearpage
\beginsupplement
\section*{Appendix}

\begin{figure}[h!]
\centering
$C_q$ - normal $\xi$ \\
	\begin{subfigure}[b]{0.09\textwidth}	
	\includegraphics[height=1.5in]{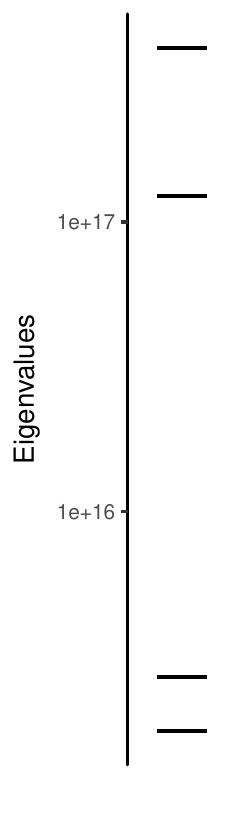}
%		\caption{}
	\end{subfigure}
	\begin{subfigure}[b]{0.22\textwidth}	
	\includegraphics[width = \textwidth]{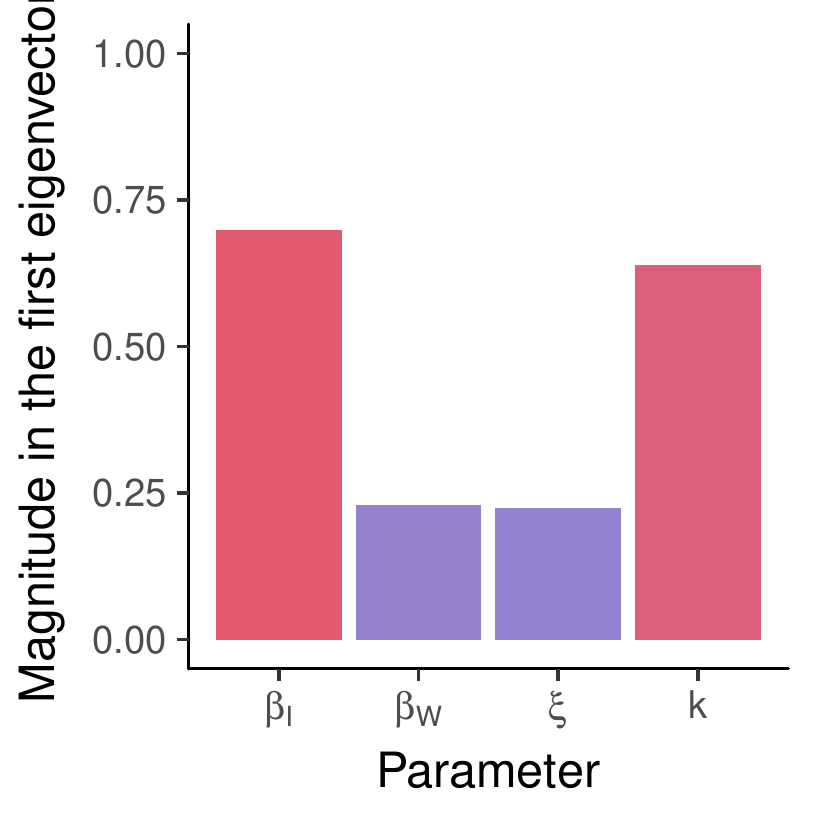}
	\end{subfigure}
	\begin{subfigure}[b]{0.22\textwidth}
	\includegraphics[width = \textwidth]{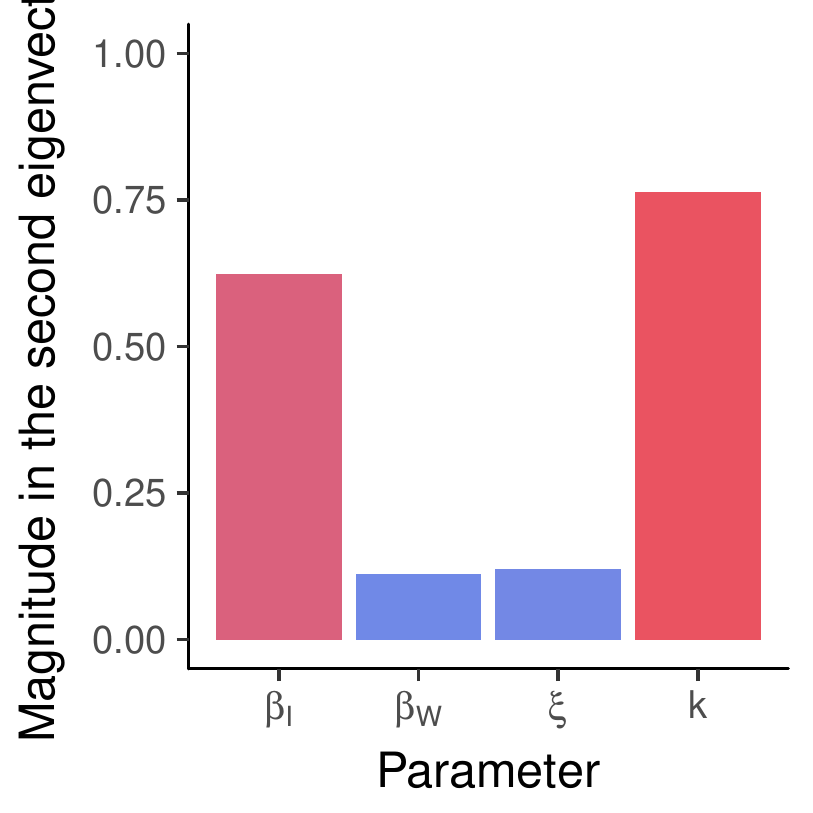}
%		\caption{}
	\end{subfigure}
	\begin{subfigure}[b]{0.22\textwidth}	
	\includegraphics[width = \textwidth]{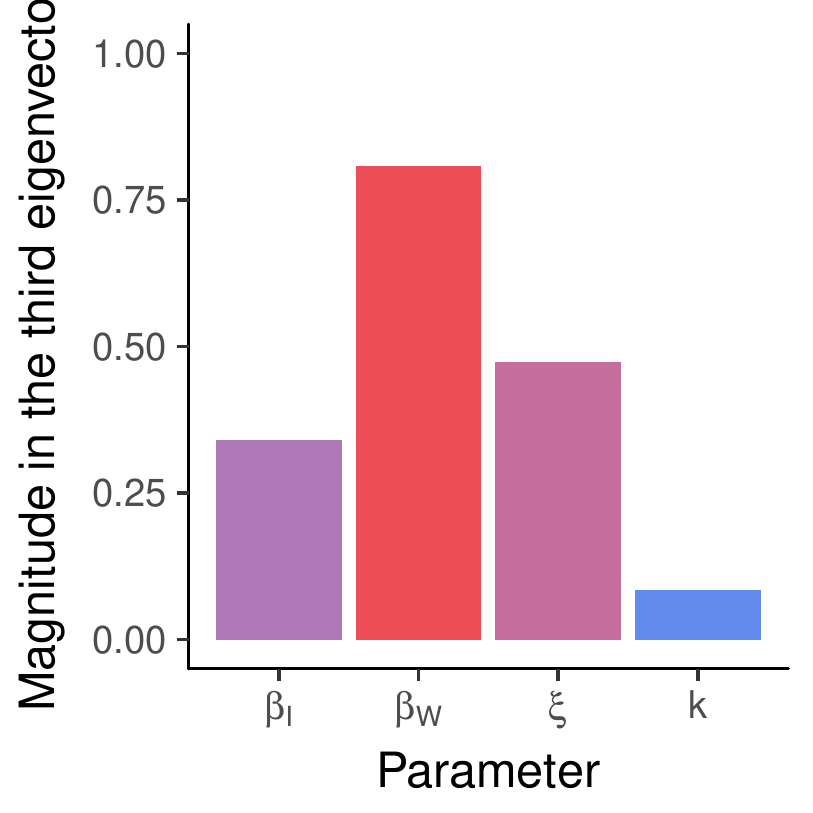}
	\end{subfigure}
	\begin{subfigure}[b]{0.22\textwidth}	
	\includegraphics[width = \textwidth]{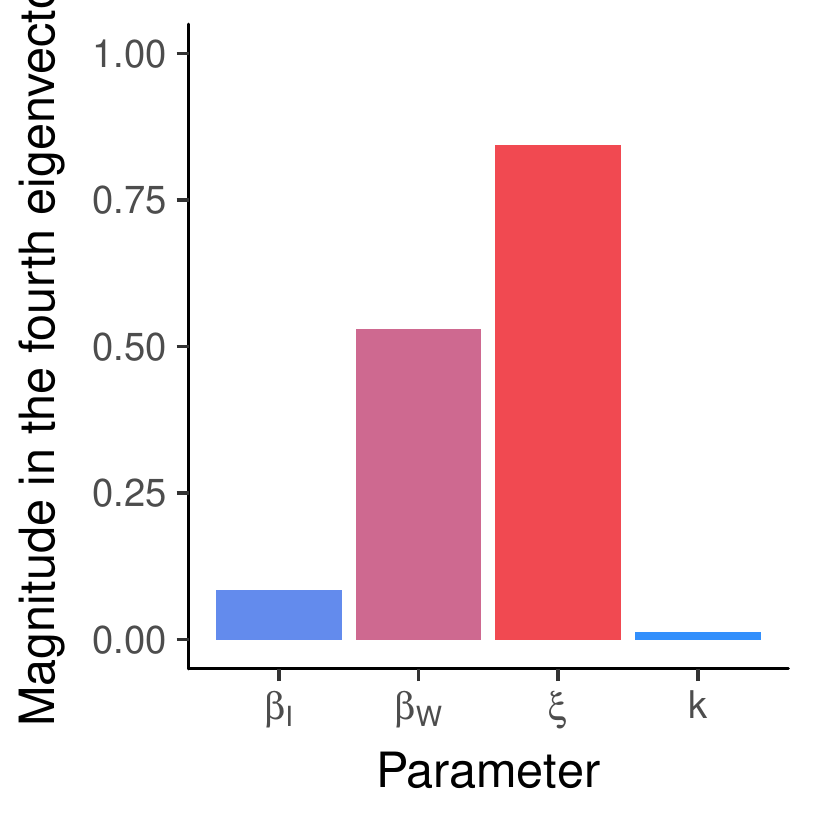}
%		\caption{}
	\end{subfigure}
\\
$C_q$ - fast $\xi$ \\
	\begin{subfigure}[b]{0.09\textwidth}	
	\includegraphics[height=1.5in]{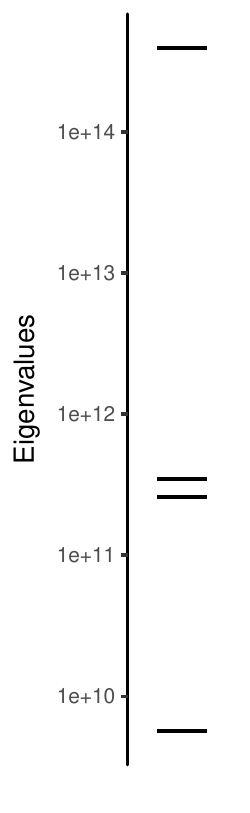}
%		\caption{}
	\end{subfigure}
	\begin{subfigure}[b]{0.22\textwidth}	
	\includegraphics[width = \textwidth]{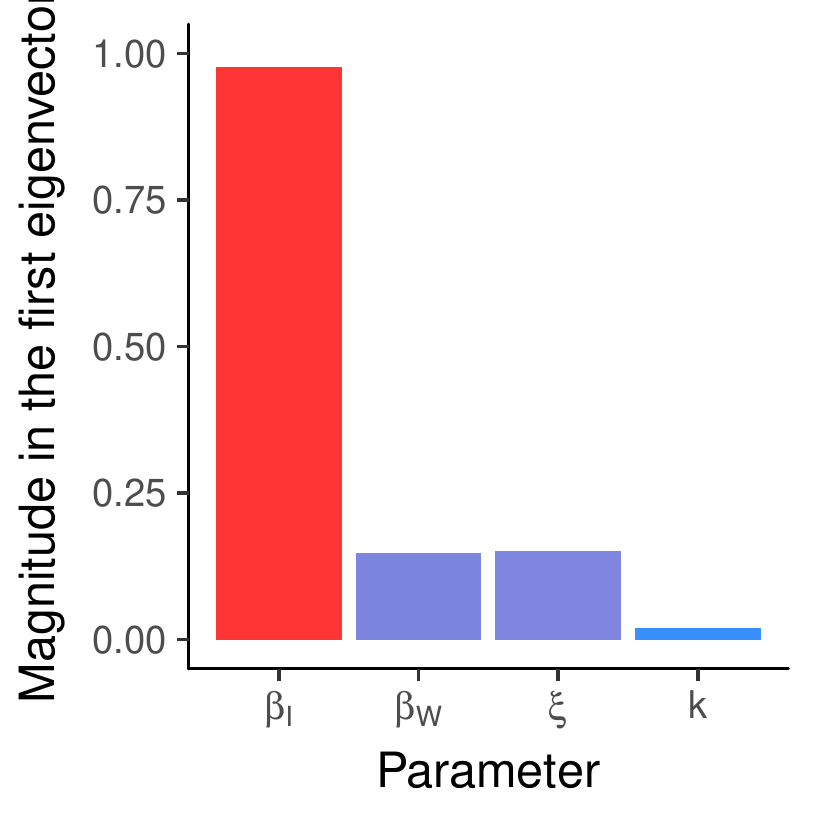}
	\end{subfigure}
	\begin{subfigure}[b]{0.22\textwidth}
	\includegraphics[width = \textwidth]{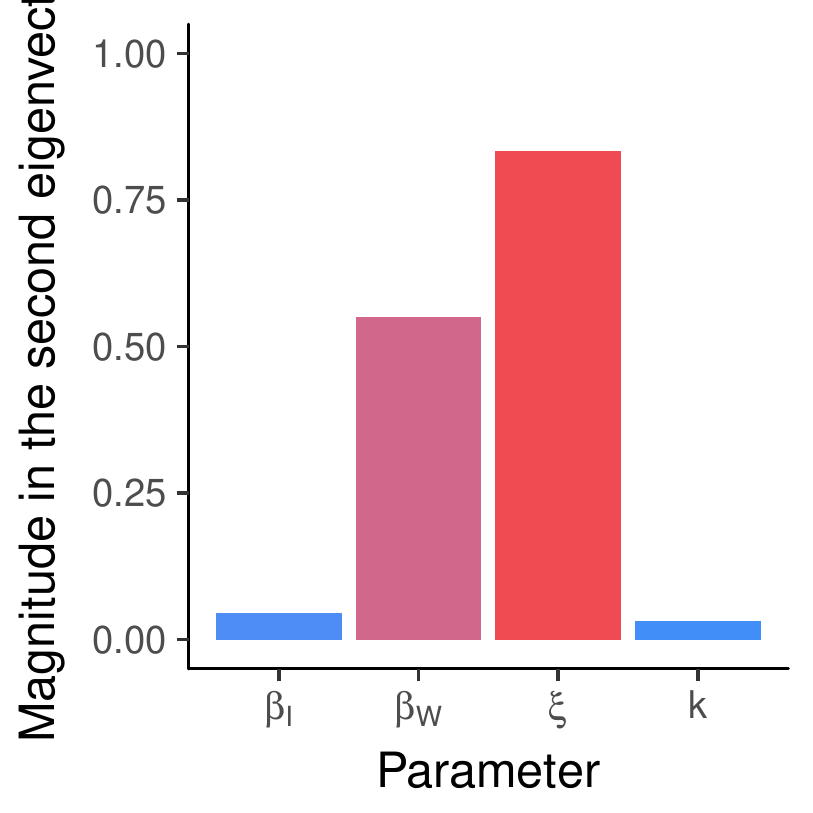}
%		\caption{}
	\end{subfigure}
	\begin{subfigure}[b]{0.22\textwidth}	
	\includegraphics[width = \textwidth]{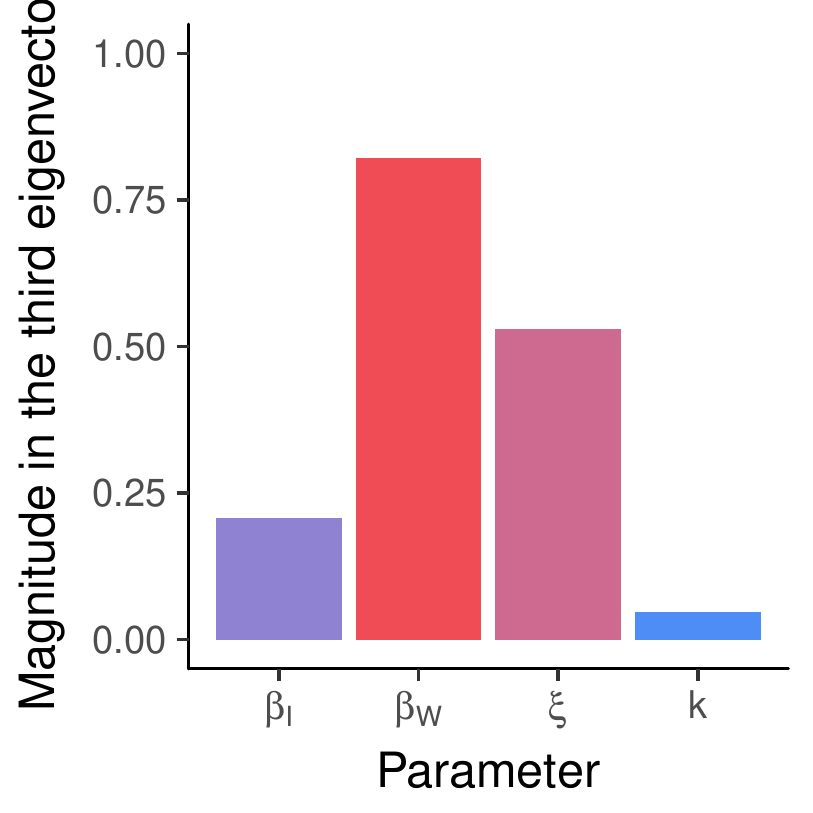}
	\end{subfigure}
	\begin{subfigure}[b]{0.22\textwidth}	
	\includegraphics[width = \textwidth]{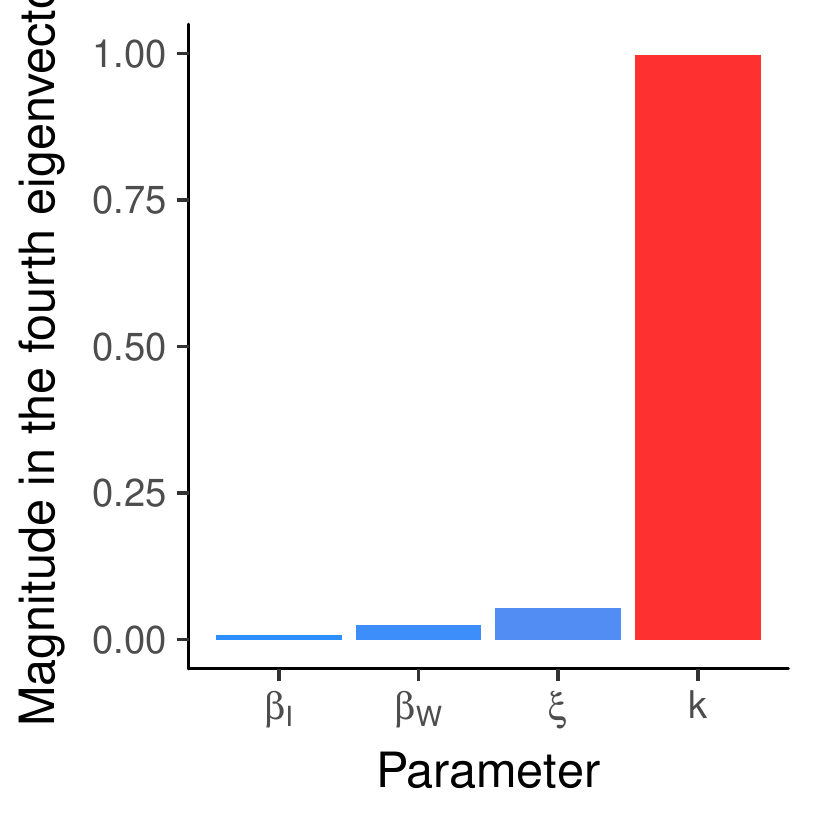}
%		\caption{}
	\end{subfigure}
\caption{Eigenvalues (leftmost column) and eigenvector component magnitudes for the average sFIM in both the normal and fast $\xi$ cases, using the least squares cost function $q$ as the QOI, and with parameters translated and scaled to be within $(-1,1)$.}
\label{fig:SIWRscaledq}
\end{figure}

\begin{figure}[h!]
\centering
$C_y$ - normal $\xi$ \\
	\begin{subfigure}[b]{0.09\textwidth}	
	\includegraphics[height=1.5in]{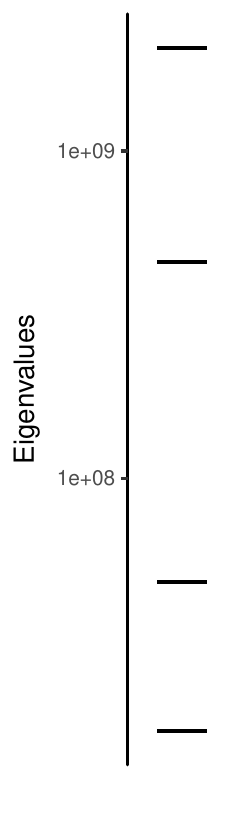}
%		\caption{}
	\end{subfigure}
	\begin{subfigure}[b]{0.22\textwidth}	
	\includegraphics[width = \textwidth]{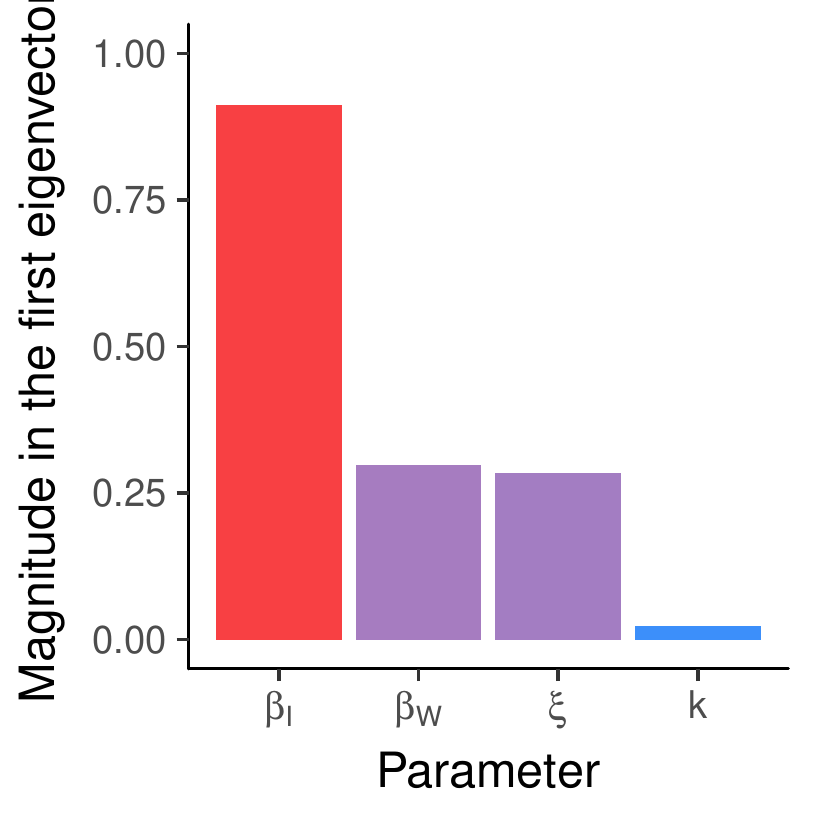}
	\end{subfigure}
	\begin{subfigure}[b]{0.22\textwidth}
	\includegraphics[width = \textwidth]{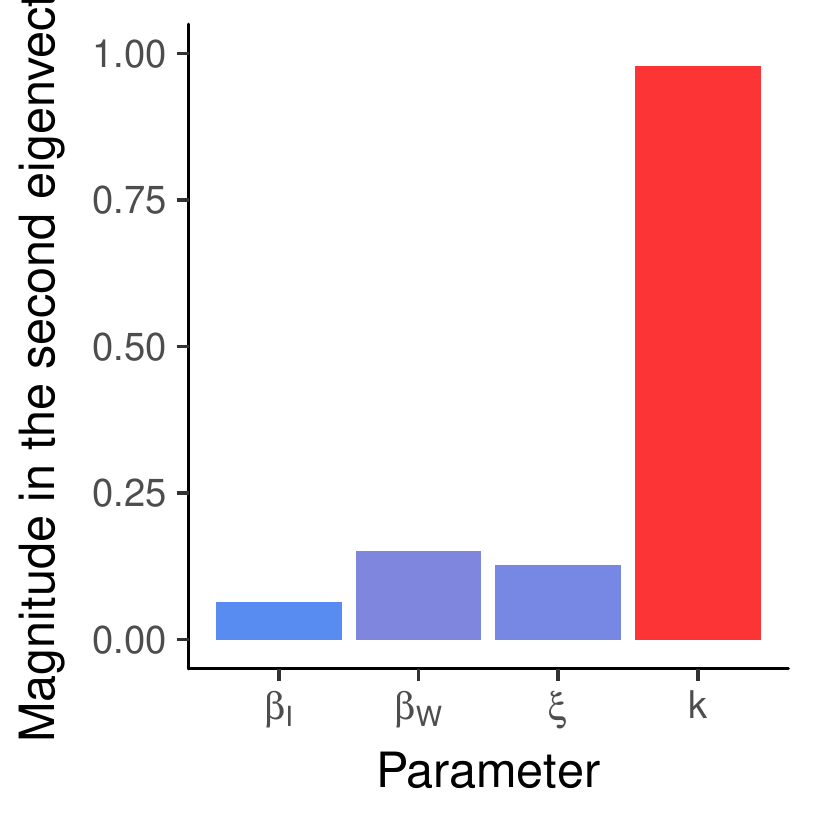}
%		\caption{}
	\end{subfigure}
	\begin{subfigure}[b]{0.22\textwidth}	
	\includegraphics[width = \textwidth]{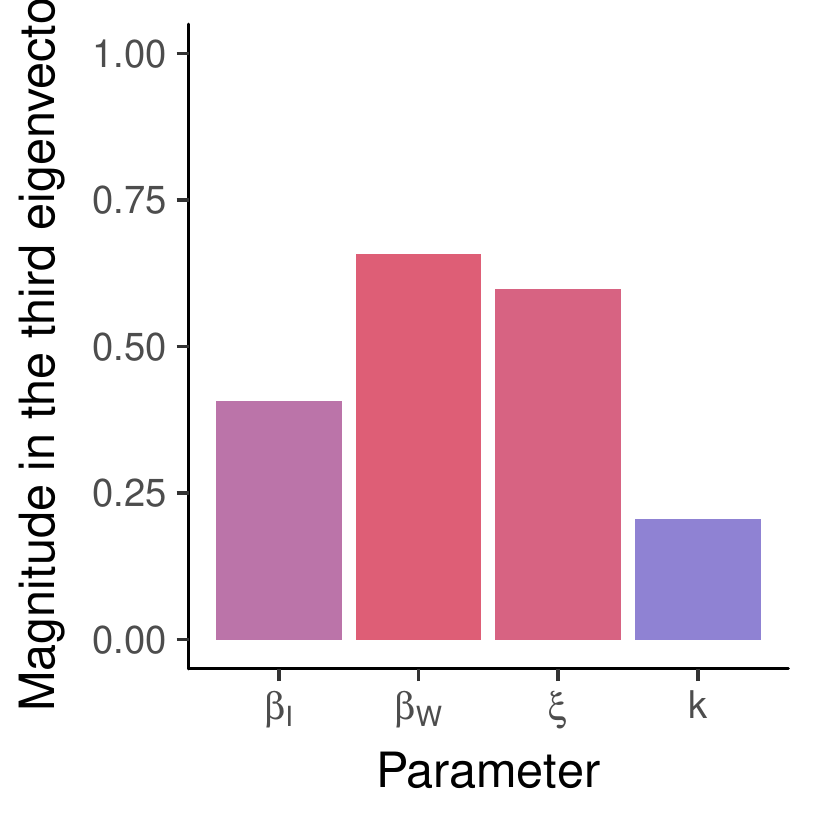}
	\end{subfigure}
	\begin{subfigure}[b]{0.22\textwidth}	
	\includegraphics[width = \textwidth]{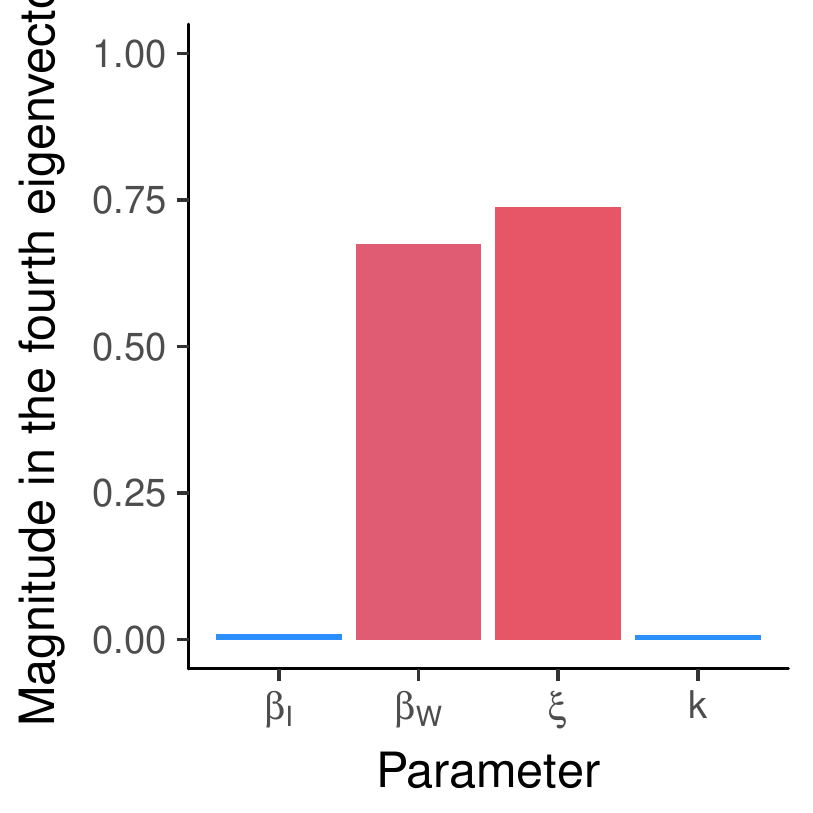}
%		\caption{}
	\end{subfigure}
\\
$C_y$ - fast $\xi$ \\
	\begin{subfigure}[b]{0.09\textwidth}	
	\includegraphics[height=1.5in]{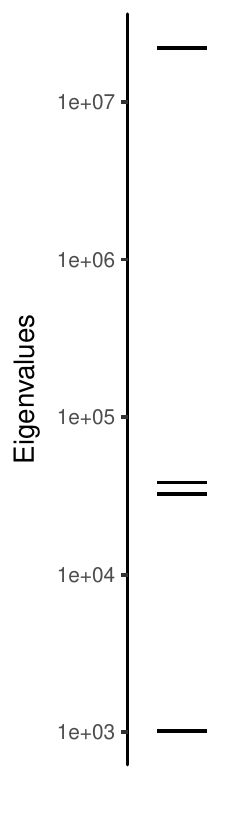}
%		\caption{}
	\end{subfigure}
	\begin{subfigure}[b]{0.22\textwidth}	
	\includegraphics[width = \textwidth]{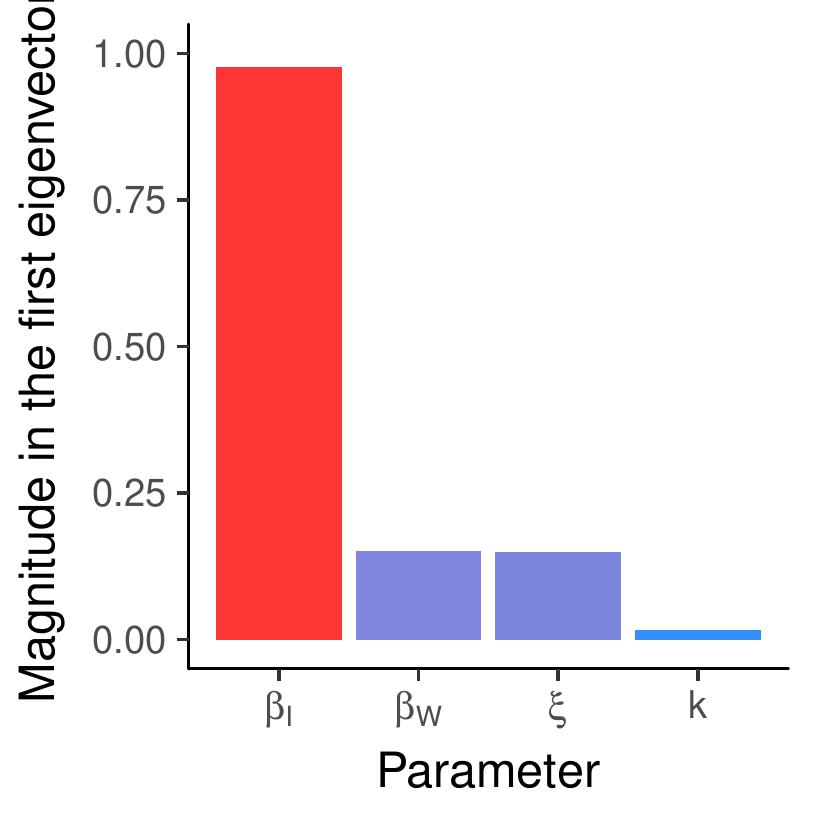}
	\end{subfigure}
	\begin{subfigure}[b]{0.22\textwidth}
	\includegraphics[width = \textwidth]{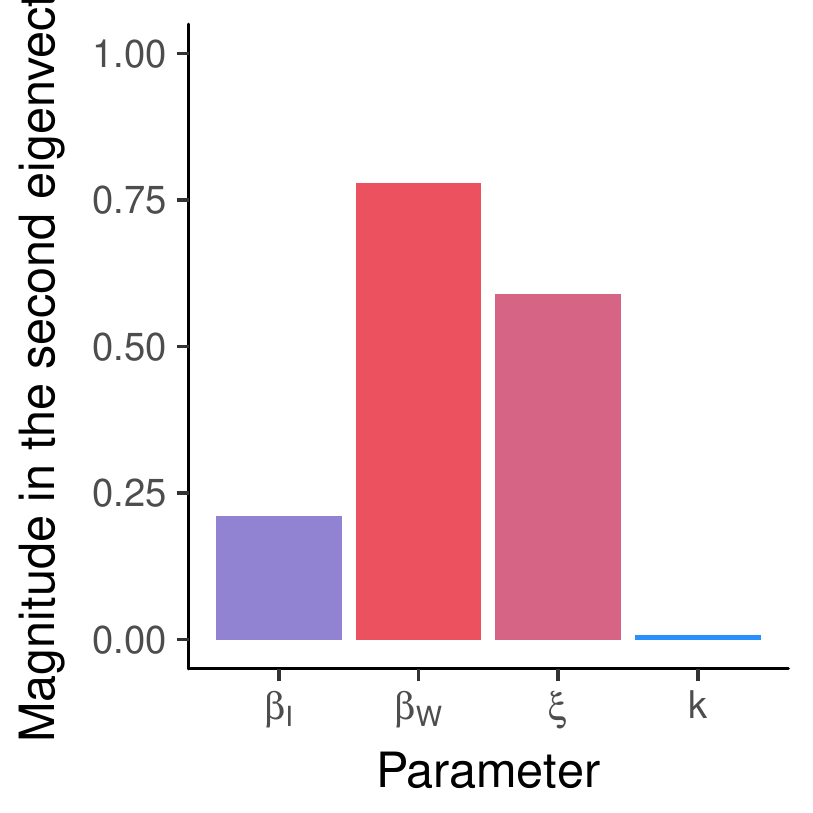}
%		\caption{}
	\end{subfigure}
	\begin{subfigure}[b]{0.22\textwidth}	
	\includegraphics[width = \textwidth]{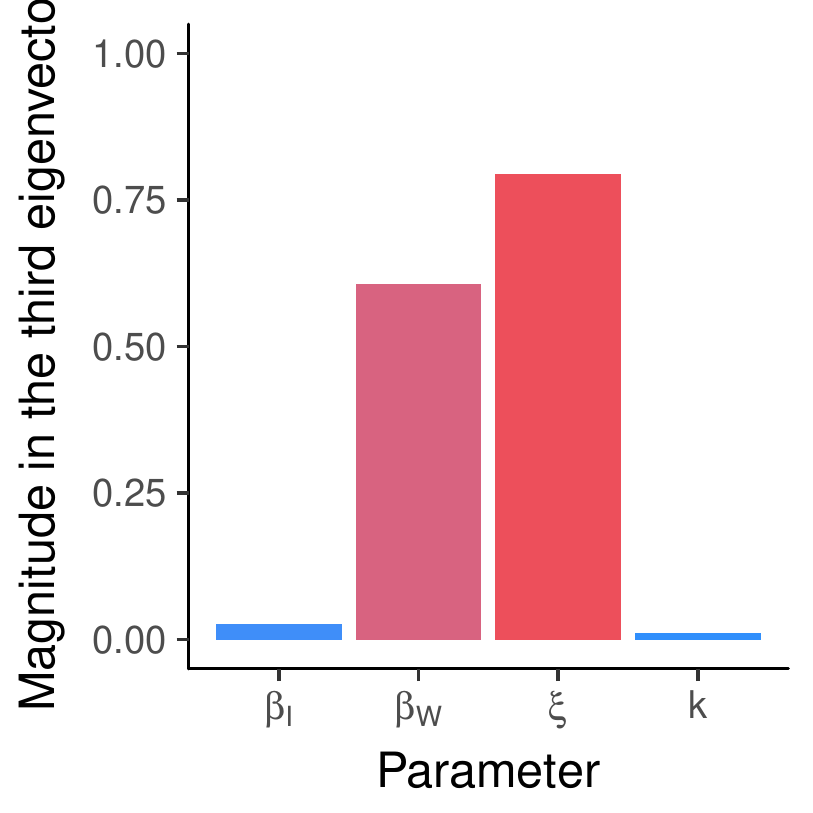}
	\end{subfigure}
	\begin{subfigure}[b]{0.22\textwidth}	
	\includegraphics[width = \textwidth]{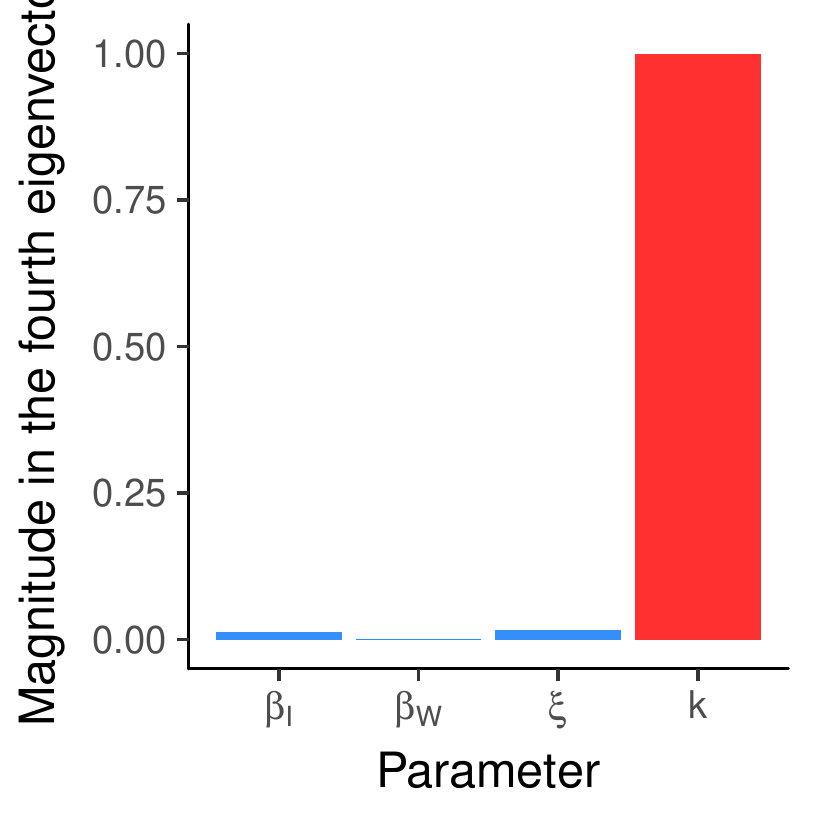}
%		\caption{}
	\end{subfigure}
\caption{Eigenvalues (leftmost column) and eigenvector component magnitudes for the average sFIM in both the normal and fast $\xi$ cases, using the time series vector $y$ as the QOI, and with parameters translated and scaled to be within $(-1,1)$.}
\label{fig:SIWRscaledy}
\end{figure}

\clearpage
%\small
\bibliographystyle{siamplain}
\bibliography{references}

\end{document}